\documentclass[a4paper,12pt]{amsart}
\numberwithin{equation}{section}

\usepackage[T1]{fontenc}
\usepackage[utf8]{inputenc}
\usepackage[english]{babel}
\usepackage{euler}
\usepackage{palatino}
\DeclareUnicodeCharacter{00B7}{\cdot}

\usepackage{microtype}

\usepackage[margin=2.5cm]{geometry}
\usepackage[linktocpage,hyperindex,breaklinks]{hyperref}
\usepackage[alphabetic]{amsrefs}
\usepackage{stmaryrd,amssymb}
\usepackage[arrow,matrix]{xy}
\usepackage{graphicx}

\ifx\texorpdfstring\undefined
  \def\texorpdfstring#1#2{#1}
\fi

\def\leq{\leqslant}
\def\geq{\geqslant}
\def\eps{\varepsilon}
\def\phi{\varphi}
\def\drond{\partial}

\renewcommand\emptyset\varnothing

\newcommand\Z{\mathbb{Z}}
\newcommand\Q{\mathbb{Q}}
\newcommand\R{\mathbb{R}}
\newcommand\bbC{\mathbb{C}}

\newcommand\into{\hookrightarrow}

\newcommand\set[1]{\ensuremath{\left\{#1\right\}}}

\newcommand\pair[1]{\ensuremath{\left\langle#1\right\rangle}}

\DeclareMathOperator{\Tr}{Tr}
\DeclareMathOperator{\Hom}{Hom}
\DeclareMathOperator{\Ext}{Ext}
\DeclareMathOperator{\Extrond}{\underline{\mathbf{Ext}}}

\DeclareMathOperator{\Sym}{Sym}

\DeclareMathOperator{\SL}{SL}
\DeclareMathOperator{\GL}{GL}
\DeclareMathOperator{\Gr}{Gr}
\def\MGITexp[#1]{{\mathcal M}^{\mathrm{#1}}_{\mathrm{GIT}}}
\def\MGIT{{\mathcal M}_{\mathrm{GIT}}}
\makeatletter
\def\Mgit{\@ifnextchar[{\MGITexp}{\MGIT}}
\makeatother

\newcommand{\Ocal}{\mathcal O}

\DeclareMathOperator{\Proj}{Proj}
\DeclareMathOperator{\Pic}{Pic}
\DeclareMathOperator{\Res}{Res}

\DeclareMathOperator{\Aut}{Aut}

\DeclareMathOperator{\id}{id}

\def\arXiv#1{%
    \href{http://arxiv.org/abs/#1}{arXiv:#1}%
}

\theoremstyle{plain}
\newtheorem{thm}{Theorem}[section]
\newtheorem*{thm*}{Theorem}

\newtheorem{prop}[thm]{Proposition}
\newtheorem{lemma}[thm]{Lemma}
\newtheorem{corol}[thm]{Corollary}

\newtheorem*{defn}{Definition}

\theoremstyle{remark}

\title{Periods of an arrangement of six lines and Campedelli surfaces}
\author{Rémy Oudompheng}
\date\today

\begin{document}

\begin{abstract}
  We define a period map for classical Campedelli surfaces, using a
  covering trick as in the case of Enriques surfaces: the period map
  is shown to come from a family of Enriques surfaces,
  obtained as quotients of the Campedelli surface by an involution.

  The period map realises an isomorphism between a projective variety
  obtained by invariant theory, and the Baily-Borel compactification
  of an arithmetic quotient, in the same fashion as in the work of
  Matsumoto, Sasaki and Yoshida. The result is proved from scratch
  using traditional methods.

  As another consequence we determine properties of the monodromy of
  Campedelli surfaces with a choice of double cover.
\end{abstract}

\maketitle

\setcounter{tocdepth}{1}
\tableofcontents

\newpage

\section*{Introduction}

This work aims at investigating the periods of Campedelli surfaces,
and their monodromy. Campedelli surfaces were among the first examples
of general type surfaces with $p_g=0$, meaning that there exists no
holomorphic differential 2-form on them: they do not have periods in the
traditional sense. However, they do have nonzero 2-forms with values
in a non trivial local system. Much like in the case of Enriques
surfaces, this allows to define a period vector, which corresponds to
periods of an étale double cover.

\subsection*{Campedelli surfaces}

The construction of Campedelli surfaces is very explicit: consider the
projective space $\mathbb P^6$ and the linear system of quadrics generated
by the forms $x_i^2$, where $x_1, \dots, x_7$ are a system of projective
coordinates. Selecting a subsystem generated by four general equations
defines a surface in $\mathbb P^6$, which is a smooth complete
intersection.

Consider also the finite group $G = (\mathbb Z/2)^3$ and its seven
non-trivial characters with values in $\set{\pm 1}$: it defines an
action on $\mathbb P^6$, when given a mapping between these characters
and the seven coordinates $x_i$. In generic situations, the action of
$G$ is free, and Campedelli surfaces are obtained as quotients under
this action. A less general choice of equations yields surfaces with
rational double points: if the action of $G$ is still free, the
quotient surface is the canonical model of a smooth surface containing
rational $(-2)$-curves, which still satisfies $p_g=q=0$, and can be
considered an appropriate extension of the definition. In this
construction, $G$ can be identified with the topological fundamental
group of Campedelli surfaces.

Campedelli surfaces have a natural projective parameter space
\cite{Miyaoka-campedelli} which is the Grassmann variety $\Gr(4,7)$,
corresponding to the choice of its equations inside the vector space
generated by the $x_i^2$ polynomials.
The parameter point of a Campedelli surface $X$ is only defined
up to a $G$-equivariant change of coordinates (action of the diagonal
torus $T$). It also depends on an identification between $G$ and $(\mathbb
Z/2)^3$, which determines the labelling of the seven coordinates.
It is thus natural to consider the variety $(\Gr(4,7) \sslash T) /
\GL_3(\mathbb F_2)$ as a (compactification of a) moduli variety for
Campedelli surfaces. The finite group $\GL_3(\mathbb F_2)$ represents
the coordinate permutations that do not change the resulting surface.

In order to study a period map, we are required to fix a non trivial
character $\kappa$ of $G$ (in order to define the period vector): the
natural coarse moduli variety for pairs $(X,\kappa)$ is then
$(\Gr(4,7) \sslash T) / \mathfrak S_4$. The finite group which appears
here is the affine group of $(\mathbb F_2)^2$, which is isomorphic to
$\mathfrak S_4$: it is the subgroup of $\GL_3(\mathbb F_2)$ which
fixes the chosen character $\kappa$. It also acts on the coordinates
of $\mathbb P^6$, permuting the six other coordinates (in the same way
as it permutes the six pairs of numbers among $\set{1,2,3,4}$).

\subsection*{Configurations of lines}

The various involutions $x_i \mapsto -x_i$ of $\mathbb P^6$ act on
intersections of diagonal quadrics: in particular, each of them
descends to an involution $s_i$ of Campedelli surfaces as defined
above, and generate a subgroup $(\Z/2)^6/G \simeq (\Z/2)^3$ of their automorphism group.
A geometric way of realising this quotient is obtained by squaring all
coordinates, which gives a finite morphism $\mathbb P^6 \to \mathbb
P^6$. This morphism is also well-defined on a Campedelli surface,
and the image of the induced morphism $X \to \mathbb P^6$ is a
plane (since the squared coordinates of points of $X$ satisfy linear relations).
This exhibits $X$ as an abelian cover of $\mathbb P^2$ ramified
over a configuration of seven lines \cites{AlexeevPardini,Pardini-covers}.
Another moduli variety for Campedelli surfaces is then the GIT
quotient $((\mathbb P^2)^7 \sslash PSL_3) / \GL_3(\mathbb F_2)$, where
$\GL_3(\mathbb F_2)$ is seen as a subgroup of $\mathfrak S_7$.

We will find that the differential of the period map is never
injective: indeed, an étale double cover of a Campedelli surface is a
special case of a Todorov surface, whose period map is known not to be
injective. More precisely, the period map of Campedelli surfaces can
be factored through that of a family of lattice-polarised Enriques
surfaces, whose natural parameter space is $\Gr(3,6)$, up to commuting
actions of ${\Z/2 \wr \mathfrak S_3} \simeq \Z/2 \rtimes \mathfrak
S_4 \subset \mathfrak S_6$ and an involution $Q$ which is described in section
\ref{ssec:involution}. The correspondence between Campedelli surfaces
and Enriques surfaces is given by the rational map $\Gr(4,7)
\dashrightarrow \Gr(3,6)$ which maps a (generic) 4-dimensional
subspace of $\bbC^7$ to its 3-dimensional intersection with a given
hyperplane $\bbC^6 \subset \bbC^7$. In terms of configurations of
seven lines, it translates to the fact that periods with values in the
local system $\Z_\kappa$ do not depend on the position of the line
labelled by $\kappa$.

This correspondance gives insight about the relationship between
Campedelli surfaces and their periods: the isomorphism classes of
Enriques surfaces are determined by their periods according to the
work of Horikawa, and the correspondance between Campedelli surfaces
and Enriques surfaces described above is given by geometry.

The subgroups of $\mathfrak S_6$ which appear are related in the
following way: $\mathfrak S_4$ acts on the six 2-element subsets
of $\set{1,2,3,4}$. Since $\mathfrak
S_4$ acts on the 3 decompositions of $\set{1,2,3,4}$ as complementary
pairs, we can write it as a semi-direct product $V_4 \rtimes \mathfrak
S_3$ where $V_4$ is the group of double transpositions. We have
inclusions
\[ \mathfrak S_4 \simeq (\Z/2)^2 \rtimes \mathfrak S_3 \subset
(\Z/2) \wr \mathfrak S_3 \subset \mathfrak S_6 \]
where $(\Z/2) \wr \mathfrak S_3$ is the \emph{wreath product} of
$\Z/2$ by $\mathfrak S_3$, which also acts naturally on three pairs
of objects. All of these groups act on $\bbC^6$ by
permutation of coordinates, and on $\bbC^7$ by adding a fixed
coordinate.

\subsection*{Enriques surfaces polarised by a $D_6$ lattice}

The Enriques surfaces appearing as quotients of Campedelli surfaces
are also naturally parametrised by their period space: geometric
arguments show that their $H^{1,1}(\Z)$ contains a distinguished copy
of the $D_6$ lattice. Their generic transcendental lattice
(orthogonal complement to the $D_6$ lattice) is then isomorphic to $L
= \Z^2(2) \oplus \Z^4(-1)$, which defines a bounded symmetric domain
$\mathcal D_L$ and an arithmetic quotient $\mathcal X_L = \mathcal D_L
/ O(L)$. The question is then: is it true that all (or almost all)
such Enriques surfaces are obtained from Campedelli surfaces (and what
are the geometric properties of the correspondance)?

We can obtain the following structure:
\begin{thm*}
  The structure of the period map of Campedelli surfaces
  $\Gr(4,7) \dashrightarrow \mathcal X_L$ can be
  described by the following diagram
  \[ \xymatrix{
    \Gr(4,7) \ar@{-->}[r] \ar[d] & \Gr(3,6) \ar[d] & \\
    \Gr(4,7) / \mathfrak S_4 \ar@{-->}[r] &
    \Gr(3,6) / (\Z/2)^2 \rtimes \mathfrak S_3 \ar[r] &
    \frac{\Gr(3,6)}{(\Z/2 \wr \mathfrak S_3) \times \pair{Q}}
    \simeq \mathcal X_L^\mathrm{BB} \\
  } \]
  where $\mathcal X_L^\mathrm{BB}$ is the Baily-Borel compactification
  of $\mathcal X_L$. In the leftmost square, the various arrows are
  the rational maps coming from the geomtric correspondances described
  earlier.
\end{thm*}

This statement includes the fact that the natural GIT moduli space for
Enriques surfaces polarised by the lattice $D_6$, whose presentation
is deduced from the moduli space we chose for Campedelli surfaces, is isomorphic to a
Baily-Borel compactification by means of the period map. The bulk of
this paper consists in a proof of this fact, which is to be related to
the work of Matsumoto, Sasaki, Yoshida \cite{MSY}, who prove a similar
statement for K3 surfaces which are double planes ramified over six
lines, with different techniques. Since indeed, a $D_6$-polarised
Enriques surfaces are bidouble covers of the plane, also ramified
over six lines: from both the geometric and lattice-theoretic pont of
view, the moduli space appearing here is a 15-fold cover of
theirs, but this aspect is not studied here.

Another way of summarising the present work is the following:
\begin{thm*}
  The bidouble covers of the plane, ramified over three pairs of lines
  have a natural GIT-theoretic moduli space, given by $\Gr(3,6) / (\Z/2
  \wr \mathfrak S_3)$. The period map of Enriques surfaces realises
  this space as a degree 2 cover of the natural Baily-Borel
  compactification of the associated period space.

  The isomorphism classes of Enriques surfaces which appear are polarised by a lattice
  $D_6$, and correspond generically to exactly two nonequivalent
  configurations of lines, related by a Cremona transformation.
  Such a generic Enriques surface gives rise to four deformation types of
  Campedelli surfaces with choice of a local system.
\end{thm*}

Note that given a generic Enriques surface as above, there are two
choices of linear systems (corresponding to the inequivalent polarisations
of degree four which describes it as a bidouble cover of $\mathbb
P^2$), such that a double cover, with ramification locus chosen in
these linear systems is a Campedelli surface. The additional factor
two arises from the fact that there are two ways of constructing
double covers with given ramification locus, since an Enriques surface
is not simply connected.

I express my gratitude to my PhD adviser Arnaud Beauville for
all the wise remarks he made during the progression of this work,
as well as Carlos Simpson, who kindly answered several of my questions.
I also warmfully thank Olivier Debarre, Christoph Sorger, and Claire
Voisin, for their remarks and comments on this work.

\subsection*{Notations}

In the table below we gather a list of common notations used
throughout the paper.
\begin{center}
\begin{tabular}{l p{.8\textwidth}}
$X$ & a Campedelli surface or its canonical model \\
$\tilde X$ or $Y$ & the universal cover of $X$ \\
$G$ & the group $(\Z/2\Z)^3$ \\
$\kappa$ & a character of $\pi_1(X)$ or $G$ \\
$\Z_\kappa$, $\R_\kappa$, $\bbC_\kappa$ & the local systems of
  $\Z$-modules, $\R$-vector spaces, $\bbC$-vector spaces attached to
  $\kappa$ \\
$\mathcal L_\kappa$ & the 2-torsion holomorphic line bundle associated
  to $\kappa$ \\
$X_\kappa$ & the étale double cover of $X$ associated to $\kappa$ \\
$s_\kappa$ & the involution of $X$ associated to $\kappa$ \\
$S_\kappa$ & the quotient of $X$ by $s_\kappa$ \\
$T_\kappa$ & the quotient of $X_\kappa$ by the natural lift of $s_\kappa$ \\
$S$ & an Enriques surface \\
$T$ & the K3 universal over of $S$ \\
$\Z^{p,q}$ & the odd unimodular integral lattice of signature $(p,q)$
  \\
$D_{p,q}$ & the even sublattice of $\Z^{p,q}$ \\
$D_n$ & the even sublattice of the Euclidean lattice $\Z^n$ \\
$L$ & the lattice $\Z^2(2) \oplus \Z^4(-1)$ \\
\end{tabular}
\end{center}

\section{Campedelli surfaces}

Campedelli surfaces (and numerical Campedelli surfaces which share the
same numerical invariants) have been thoroughly studied
\cites{Campedelli,Miyaoka-campedelli,Pardini-involution}. A number of useful results are
contained in the unpublished manuscript of M. Reid \cite{Reid-k22}.

In this section we give the basic definitions of Campedelli surfaces,
and their first properties. We are mainly concerned about involutions
that are induced by coordinate reflections in $\mathbb P^6$. In the
generic case, the quotient of a Campedelli surface $X$ by such a
reflection is an Enriques surface whose periods determine those of $X$
(proposition \ref{prop:campedelli-quotient-is-enriques}). This
Enriques surface usually has six nodes corresponding to isolated fixed
points on $X$.

We then give local properties of the period map. Its differential has
rank 4 at points of the moduli space parameterising smooth surfaces
(proposition \ref{prop:target-dim-4}). We will define later a natural
4-dimensional period domain for Campedelli surfaces. We also prove
in section \ref{ssec:period-double-pt},
along the lines of \cite{Voisin-4fold}, a property which is needed
later: when crossing the discriminant hypersurface, the period of the
associated vanishing cycle has non-zero derivative, on the double
cover ramified over the (corresponding irreducible component of the)
discriminant.

\subsection{Description and general properties}

\begin{defn}
  A \emph{numerical Campedelli surface} is a minimal smooth projective
  surface $X$ with numerical invariants $p_g = q = 0$ and $K_X^2 = 2$.

  A (classical) \emph{Campedelli surface} is a numerical Campedelli
  surface whose fundamental group is isomorphic to $(\Z/2)^{\oplus 3}$.
\end{defn}

Many topological invariants of Campedelli surfaces can be calculated
in terms of these numbers: since $\chi(\Ocal_X)=1$, by Noether's formula, the topological Euler
characteristic is $e(X) = c_2(X) = 12 - K_X^2 = 10$, thus the nonzero
Hodge numbers of $X$ are $h^{0,0} = h^{2,2} = 1$ and $h^{1,1} =
8$. The signature of $X$ is $\tau(X) = (K_X^2 - 2e(X))/3 = -6$, hence
the torsion-free quotient of $H^2(X, \Z)$ (denoted by $H^2(X,
\Z)_\mathrm{num}$), which is a $\Z$-module of rank 8 with a
unimodular quadratic form, is isomorphic to $\Z^{1,7}$, the standard
Lorentzian lattice.





\index{surface!de Campedelli}
Campedelli surfaces have a six-dimensional, unirational moduli
variety. This is implied by the following structure theorem:
\begin{thm}[see \cite{Miyaoka} or \cite{Reid-k22}]
  The universal cover $\tilde X$ of a (classical) Campedelli surface
  $X$ is birational to a complete intersection of 4 diagonal quadrics
  in $\mathbb P^6$, where $\pi_1(X)$ acts by its 7 distinct nonzero
  characters. Moreover, this complete intersection is the canonical
  model of $\tilde X$, and its quotient is the canonical model of $X$.
\end{thm}
Because of this simple description, we give the name of
\emph{canonical Campedelli surface} to the canonical model itself.

An abstract approach to this property is the fact that the action of
$G = \pi_1(X)$ on $H^0(\widetilde{X}, K_{\widetilde{X}})$ is
decomposed as a sum of eigenspaces for each character $\kappa \in \hat
G$, which can be identified with $H^0(X, K_X \otimes L_\kappa)$,
where $L_\kappa$ is the (flat) line bundle arising from the
representation of $\pi_1(X)$ given by $\kappa$. Each of these spaces
has dimension one, except for $\kappa=0$ (since $p_g(X) = h^0(X, K_X) = 0$).
Up to homothety, it is then possible to choose canonical coordinates
$x_\kappa$ (where $0 \neq \kappa \in \hat G$) for the embedding of $\tilde X$
in $|K_{\tilde X}|^\vee \simeq \mathbb P^6$.

Since diagonal quadratic equations in $\mathbb P^6$ form a vector
space of rank 7, the choices of a linear system spanned by four
elements can be identified with points in the Grassmann variety $\Gr(4,7)$: this
gives an identification between the moduli space of triples:
\[ (X, f: G = (\Z/2\Z)^{\oplus 3} \simeq \pi_1(X), g: \mathbb P^6
\simeq |K_{\tilde X}|^\ast) \]
consisting of a Campedelli surface with a framing of its fundamental
group and a $G$-equivariant linear isomorphism of $|K_{\tilde X}|^\ast$ with
$\mathbb P^6$, and an open Zariski subset
$\mathcal M_u$ of $\Gr(4,7)$ which parametrises
linear systems of diagonal quadrics whose base locus is a normal
surface with at worst ordinary double points as singularities.

The smoothness of such a complete intersection is easy to detect:
\begin{prop}
  \label{prop:campedelli-smoothness}
  A complete intersection of four diagonal quadrics in $\mathbb P^6$ is smooth if and
  only if its linear system does not contain a quadric of rank three.
  Under this assumption no point has more than two vanishing
  coordinates.
\end{prop}

\begin{proof}
  Let $\tilde X \subset \mathbb P^6$ denote the surface defined by
  the linear system. If it contains a quadric of rank three, e.g.
  $x_5^2 + x_6^2 + x_7^2$, there is a point of $\tilde X$ such that
  $x_5 = x_6 = x_7 = 0$. Such a point is indeed defined by
  three quadratic equations inside the space $\set{x_5=x_6=x_7=0}
  \simeq \mathbb P^3$. By the Jacobian criterion,
  the complete intersection cannot be smooth at this point, which is a
  vertex of $\set{x_5^2 + x_6^2 + x_7^2 = 0}$.

  Conversely, if the quadrics have the form, $Q_j = \sum q_{ij}
  x_i^2$, the Jacobian matrix of the linear system
  is given by $(q_{ij} x_i)_{i,j}$. Assume this matrix has not full
  rank at a point $R$. If $R$ has four nonzero coordinates, the
  corresponding minor should vanish, and this implies that
  some linear combination of the $Q_j$ has rank three.
  If $R$ has four zero coordinates (say $x_1=\dots=x_4=0$), since $R$ is solution to the
  equations, the matrix $(q_{ij})_{i=5,6,7}$ has rank at most two.
  Thus there is a two dimensional subsystem of $\langle Q_j \rangle$
  which is contained in $\pair{x_1^2, \dots, x_4^2}$ and one of them
  has rank three.
\end{proof}

\begin{defn}
  A \emph{framed Campedelli surface} is a pair $(X, \phi: G \to
  \pi_1(X))$ where $X$ is a Campedelli surface, and $\phi$ is an
  isomorphism between $G=(\Z/2\Z)^3$ and $\pi_1(X)$.
\end{defn}

A natural moduli space for framed Campedelli surfaces is the GIT
quotient of $\Gr(4,7)$ under action of the diagonal torus $T$ in $\SL_7$.
Indeed, the set of $G$-equivariant isomorphisms of $H^0({\tilde X}, K_{\tilde X})^\vee$
with $\bbC^7$ is a torsor under $T$: any such isomorphism must be
compatible with the splitting of these spaces into a direct sum of
one-dimensional eigenspaces.

Our main object of interest is the following:
\begin{defn}
  A \emph{marked Campedelli surface} is a pair $\kappa = (X,
  X_\kappa \to X)$ where $X$ is a Campedelli surface and $X_\kappa \to
  X$ is an étale connected double cover of $X$.
\end{defn}

The datum of a marked Campedelli surface amounts to give, along
with $X$, one of the following (equivalent) structures:
\begin{itemize}
\item a connected étale double cover $X_\kappa \to X$;
\item a non-trivial rank one local system of integral coefficients
  $\Z_\kappa$ (whose square is the constant sheaf);
\item a 2-torsion non trivial holomorphic line bundle on $X$,
  denoted by $\mathcal L_\kappa$;
\item a non trivial character of $\pi_1(X)$, $\kappa: \pi_1(X) \to \Z/2\Z$;
\item an effective divisor numerically equivalent to $K_X$
  (which must have the form $H_\kappa \in |K_X+\mathcal L_\kappa|$).
\end{itemize}

\begin{prop}
  The moduli space of framed Campedelli surfaces is connected and can
  be written as a $\mathfrak S_4$-Galois cover of the moduli space of
  marked Campedelli surfaces.
\end{prop}

\begin{proof}
  The set of isomorphisms between $G=(\Z/2\Z)^{\oplus 3}$ and
  $\pi_1(X)$ mapping a chosen character $\kappa$ to the character
  $(b_1, b_2, b_3) \mapsto b_1$ is naturally equipped with a free
  transitive action of the affine linear group $\mathrm{GA}_2(\mathbb F_2)$ of the plane
  $\mathbb A^2(\mathbb F_2)$, and $\mathrm{GA}_2(\mathbb F_2)$
  is canonically isomorphic to the permutation group of the
  subset $\set{g \in G \text{ such that } \kappa(g) = 1}$.
\end{proof}

The \emph{periods} of a marked Campedelli surface, with a given
holomorphic 2-form with values in $\bbC_\kappa$ (equivalently,
an element $\omega \in H^0(K_X+\mathcal L_\kappa)$), are the various
integrals of $\omega$ along cycles of $H_2(X, \Z_\kappa)$: they are
fully determined by the associated element $\omega \in H^2(X,
\bbC_\kappa)$.

\begin{prop}
  The lattice $H^2(X, \Z_\kappa)$ (modulo torsion) is equipped with
  the quadratic form defined by the cup product with values in $H^4(X,
  \Z_\kappa \otimes \Z_\kappa) \simeq \Z$. As such, it is unimodular
  and isomorphic to $\Z^{2,8}$, the standard odd quadratic lattice with
  signature $(2,8)$.
\end{prop}
\begin{proof}
  The cohomology of $\Z_\kappa$ has the same numerical invariants as
  the usual cohomology of $X$: $e(X) = 10$, and $\tau(X) = 6$. Its
  only nonzero Betti number is $b_2(X, \Z_\kappa) = 10$, and the value
  of $\tau(X)$ forces the signature to be $(2,8)$.

  By Poincaré-Verdier duality, $H^2(X, \Z_\kappa)$ is unimodular,
  and since we know it is indefinite and odd (the signature is
  not a multiple of 8), its isomorphism class is uniquely determined
  (see \cite{Serre}).
\end{proof}


\subsection{Involutions of a Campedelli surface}

Let $X$ be a \emph{canonical} Campedelli surface, and $G = (\Z/2\Z)^{\oplus 3}$ as before.
Let $\tilde X$ be the universal cover of $X$. Then $X$ is isomorphic to
$\tilde X / G$ and $\tilde X$ can be written as
a complete intersection of 4 diagonal quadrics in $\mathbb P^6$.

Let $\Gamma$ be the projection in $PGL_7$ of the
diagonal group $(\pm 1)^7$, acting on $\tilde X$: $G$ is naturally embedded
in $\Gamma$, its action being given by its seven non-trivial characters.


Note that $\tilde X$ is a $\Gamma$-invariant subvariety of $\mathbb
P^6$, and that squaring coordinates gives a Galois cover $\mathrm{sq}:
\mathbb P^6 \to \mathbb P^6$ with group $\Gamma$. Since
$\mathrm{sq}(\tilde X)$ is defined in $\mathbb P^6$ by four linear
equations, it is isomorphic to a plane, and the map $\tilde X \to
\mathrm{sq}(\tilde X)$ is also a $\Gamma$-Galois cover. This proves
the following proposition:
\begin{prop}
  \label{prop:campedelli-galois-cover}
  A Campedelli surface is a Galois cover of $\mathbb P^2$
  with group $\Gamma/G$. \qed
\end{prop}

The group $\Gamma/G$ will be identified with a set of reflections
in $\mathbb P^6$, giving a simple description of the group of
automorphisms of a generic Campedelli surface. The elements
of $\Gamma$ can be classified by \emph{weight}: this notion will provide
a convenient vocabulary for the rest of the paper.
\begin{defn}
  The \emph{weight} of an element of $\Gamma$ (represented by
  a diagonal matrix $g$ in $\set{\pm 1}^7 \subset \GL_7$)
  is defined as $|\Tr g|$ (the
  difference between the number of $+1$ and $-1$ coefficients
  in $g$).
\end{defn}


There are in $\Gamma$
\begin{itemize}
\item one element with weight 7: the identity;
\item 7 elements with weight 5: the reflections;
\item 21 elements with weight 3;
\item 35 elements with weight 1, seven of them being the nonzero
  elements of $G$.
\end{itemize}

\begin{prop}
  The projection $\Gamma \to \Gamma/G$ induces a bijection between the
  elements $s_i$ (reflections across coordinate hyperplanes) and
  nonzero elements of $\Gamma/G$.

  Moreover, if $\chi$ is a nonzero character of $G$, and $s_\chi$ is
  the reflection across the associated hyperplane, the map
  $\chi \mapsto [s_\chi] \in \Gamma/G$ is a group isomorphism.
\end{prop}

\begin{proof}
  Let $s_i$ and $s_j$ be different reflections in $\Gamma$,
  and suppose $s_i - s_j$ lies in $G$: then there would be an element of
  $G$ with weight 3, which is impossible. Thus $\Gamma \to \Gamma/G$
  is injective on the $s_i$'s.

  Moreover, let $i$, $j$, $k=i+j$ be nonzero elements in $\hat G$:
  then $G$ contains a unique nonzero element $g$ annihilated by $i$, $j$
  and $k$. By definition, $\pm (s_i s_j s_k)$ is the
  element of $\Gamma$ associated to $g$, hence $s_k = s_i + s_j$ in
  $\Gamma/G$.
\end{proof}

The character group $\hat G$ is thus realised as a subgroup of
$\Aut(X)$: for a generic $X$, this inclusion is even an equality.
This results from the fact that the bicanonical map must be
equivariant under the action of $\Aut(X)$, and from the fact that a
general configuration of 7 lines in the plane has no automorphisms.

Another interpretation of these facts is that the quadratic form over
$(\Z/2\Z)^7$ defined by
\[ q(x_1, \dots, x_7) = \sum_i x_i^2 + \sum_{i<j} x_i x_j \]
has a polar symplectic form
\[ b(u,v) = \sum_{i \neq j} u_i v_j = q(u+v) - q(u) - q(v) \]
whose kernel is the vector
\[ \begin{pmatrix} 1 & 1 & 1 & 1 & 1 & 1 & 1 \\ \end{pmatrix} \]
It thus defines a non-degenerate quadratic form on $\Gamma$,
for which nonzero isotropic vectors are vectors of weight $1$,
the weight being defined as the difference between the number of zero
and nonzero coordinates. If $x$ has $n$ nonzero coordinates, $q(x)
\equiv n(n+1)/2$, which is zero iff $n=0, 3, 4, 7$.


Then $G$ is an isotropic subspace of $\Gamma$ where the quadratic form
vanishes, and $b$ defines a non-degenerate pairing
between $G$ and $\Gamma/G$. Taking $G$ to be generated by the lines
of the matrix
\[ \begin{pmatrix}
  0 & 0 & 0 & 1 & 1 & 1 & 1 \\
  0 & 1 & 1 & 0 & 0 & 1 & 1 \\
  1 & 0 & 1 & 0 & 1 & 0 & 1 \\
\end{pmatrix} \]
for which coordinates are the 7 characters of $G$, it is easy to check
that the pairings $b(e_i, -)$ with the basis vectors coincide with
the nonzero element of $\hat G$ given by the $i$-th column of the
matrix.

\subsection{Quotients of Campedelli surfaces under involutions}
\label{sec:quotient}

Let $(X, \kappa)$ be a (possibly nodal) marked Campedelli surface, and $s_\kappa$
the involution of $X$ canonically associated to the character
$\kappa$, which is ramified over $D_\kappa$, the unique effective
divisor in $|K_X+L_\kappa|$ (a curve of arithmetic genus 3).
We refer to \cite{Pardini-involution} for theoretical results
about the possible quotients of numerical Campedelli surfaces
under involutions, since we are dealing here with a particular case
of those.

\begin{prop}
  \label{prop:reflection-fixed-locus}
  When $X$ is smooth, the fixed point set of $s_\kappa$ consists of
  $D_\kappa$ and six isolated points.
\end{prop}

\begin{proof}
  A point $x \in X$ is fixed under the action of $s_\kappa$ if and
  only if there exists $g \in G$ such that $s_\kappa·\tilde x =
  g·\tilde x$ for some lifting of $x$ to $\tilde X$. If $g$
  is the identity, this means that $x$ lies on $D_\kappa$.

  Otherwise, if $g+s_\kappa$ has weight 1 ($\kappa(g) \neq 0$), $x$
  has at least three vanishing coordinates: this is impossible.

  If $g+s_\kappa$ has weight 3, then $\kappa(g) = 0$ and $x$ should
  have only two vanishing coordinates (corresponding to characters
  $\chi, \chi+\kappa$ such that $\chi(g) = 0$ and $\chi \neq \kappa,
  0$). We find two fixed points for each $\chi$ corresponding to $D_\chi
  \cap D_{\chi+\kappa}$, so there are six isolated fixed points for
  $s_\kappa$.
\end{proof}

\begin{prop}
  If $x \in X$ is a node of a Campedelli surface, $x$ is the
  intersection of three divisors $D_{\kappa_1}$, $D_{\kappa_2}$,
  $D_{\kappa_3}$, such that the $\kappa_i$ generate $\hat G$, and $x$
  is fixed by any involution in $\Gamma/G$ (notably $s_\kappa$).
\end{prop}


\begin{proof}
  This follows from the description of the Campedelli surface as an
  octuple plane, and the classification of singularities which can be
  found at \cite{AlexeevPardini}, see also proposition
  \ref{prop:campedelli-smoothness}.
\end{proof}

In this description, if $\kappa$ is one of the $\kappa_i$'s, the node
lies on $D_\kappa$, and is part of the fixed locus described in proposition
\ref{prop:reflection-fixed-locus}. If $\kappa$ is a sum of two $\kappa_i$'s,
then the node is a fixed point of the type $D_\chi \cap D_{\chi+\kappa}$
also described in prop. \ref{prop:reflection-fixed-locus}. Nodes that add
new fixed points are defined by $D_{\kappa_1} \cap D_{\kappa_2} \cap
D_{\kappa_3}$ where $\sum \kappa_i = \kappa$.

There are at most four such nodes: the number of bases of $\hat G$
(up to permutation of vectors) whose sum is $\kappa$ is 4 (there
are $7 \cdot 6 \cdot 4 / 3! = 28$ bases, and the number of sums is
7, giving four bases for each possible sum).

\begin{prop}
  \label{prop:campedelli-quotient-is-enriques}
  Suppose $X$ is smooth. The quotient of $X$ by $s_\kappa$ is an
  Enriques surface $S_\kappa$ with six ordinary double points. If
  $\tilde X$ is the blow-up of the isolated fixed points, and $\tilde
  S_\kappa$ is the minimal resolution of $S_\kappa$, the morphism
  $\tilde X \to \tilde S_\kappa$ is a double cover ramified over
  $D_\kappa + \sum E_i$ where $E_i$ are the exceptional curves,
  where $D_\kappa$ is now a genus 3 curve on $S_\kappa$.
\end{prop}

\begin{proof}
  Let $X_\kappa$ be the connected étale double cover of $X$ associated
  to the character $\kappa$. Since $X_\kappa$ can be written as
  $\tilde X / H_\kappa$, where $H_\kappa \subset G$ is the kernel of
  $\kappa$, $s_\kappa$ has a natural lift to $X_\kappa$ associated to
  the reflection of $\mathbb P^6$ across the hyperplane
  $\set{x_\kappa=0}$.

  The previous argument still applies and shows that $s_\kappa$ fixes
  the image of the hyperplane $D_\kappa$ and 12 isolated points.
  Drawing a diagram for this situation (here $T_\kappa = X_\kappa /
  s_\kappa$),
  \[ \xymatrix{
    X_\kappa \ar[r] \ar[d] & X \ar[d] \\
    T_\kappa \ar[r] & S_\kappa \\
  } \]
  we note that the unique nonzero section of $K_X+L_\kappa$ lifts to a
  holomorphic 2-form on $T_\kappa$, which vanishes only along
  $D_\kappa$: considering an expression of this form in local
  coordinates along $D_\kappa$ shows that it is necessarily invariant
  under $s_\kappa$, and descends to a non-vanishing 2-form on
  $T_\kappa$. It follows that $T_\kappa$ is a K3 surface with 12
  isolated double points.

  It is now easy to check that the action of $G/H_\kappa$, defining an
  involution of $T_\kappa$, has no fixed point (it is represented by
  an element $g$ of $\Gamma$, with weight one, such that $\kappa(g) =
  1$). This proves that $S_\kappa$ is an Enriques surface.
\end{proof}


Note that if $s_\kappa$ has $\nu$ fixed nodes outside the usual fixed
locus, the quotient $S_\kappa$ is a rational surface with $K_X^2 =
-\nu$. We have seen that $\nu \leq 4$.

\subsection{Infinitesimal variation of periods}

Let $\kappa \in \hat G$ be a fixed character: following the results of
Griffiths, the infinitesimal variation of $\omega \in H^2(X,
\bbC_\kappa)$ is described by the \emph{infinitesimal
  $\kappa$-periods map}:
\[ H^1(X, TX) \to \Hom(H^0(X, K_X+L_\kappa), H^1(X,
\Omega^1_X(L_\kappa))) \]
It is indeed related to the standard situation, in the following way:
if $X_\kappa$ is the étale cover of $X$ associated to $\kappa$,
whose Galois group is identified with $\Z/2\Z$, $H^1(X, TX)$ is
identified to the invariant part of $H^1(X_\kappa, TX_\kappa)$,
and $H^0(K_{X_\kappa})$, $H^1(X_\kappa, \Omega^1_\kappa)$ can be split
into the direct sum of an invariant and anti-invariant part. The
infinitesimal variation of $\kappa$-periods is then identified
to an eigenspace of the usual variation of Hodge structure of
$X_\kappa$ (which is $\Z/2\Z$-equivariant for all elements of
$H^1(X,TX)$).

Since $|K_X+L_\kappa|$ contains only $D_\kappa$, we can use the
isomorphism $TX \simeq \Omega^1_X \otimes (-K_X)$ to identify the
infinitesimal period map with the linear map
\[ H^1(X, TX) \to H^1(X, TX(D_\kappa)) \]
(recall that the first space has dimension 6 while the target space
has dimension 8), induced by the product by a section vanishing on
$D_\kappa$.

Following the method of Konno \cite{Konno}, it will be shown in the
next section that the kernel of this map is the space of deformations such
that the family of maps $Q_t: X_t \to |2K_X|^\vee \simeq \mathbb P^2$
moves the lines $Q_t(D_\kappa)$ but not $Q_t(D_\chi)$ for $\chi \neq
\kappa$. These deformations can be described as a deformation of
coordinates
\[ x_\kappa^2 \leadsto x_\kappa^2 + \sum_{i \neq \kappa}
\eps_i x_i^2. \]

Let $\mathcal T$ be the subsheaf of $TX$ consisting of vector fields
which are tangent to $D$ (we drop $\kappa$ subscripts starting from
here) along $D$. It is defined by a Cartesian square, which induces
short exact sequences:
\[ \xymatrix{
  TX(-D) \ar[r] \ar@{=}[d] & \mathcal T \ar[r] \ar[d] & TD \ar[d] \\
  TX(-D) \ar[r] & TX \ar[r] \ar[d] & TX_{\upharpoonleft D} \ar[d] \\
  & N_D \ar@{=}[r] & N_D \\
} \]

The long exact sequences arising from the previous diagram can be used
to build a commutative diagram:
\[ \xymatrix{
  H^0(D, TX(D)) \ar@{^(-->}[r] \ar[d] &
  H^0(D, N_D(D)) \ar[d] \ar@{-->}[rd] \\ 
  H^1(X, TX) \ar@{^(->}[r] \ar@{=}[d] & H^1(X, \mathcal T(D)) \ar[r]
  \ar@{->>}[d]
  & H^1(D, TD(D)) \ar@{-->>}[d] \\
  H^1(X, TX) \ar[r] & H^1(X, TX(D)) \ar[r] \ar[d]
  & H^1(D, TX(D)) \ar@{-->}[d] \\
  & H^1(D, N_D(D)) \ar@{=}[r] & H^1(D, N_D(D)) = 0\\
} \]
where all lines, the 2nd column and the dashed arrows are exact
sequences. Note that $H^0(X, TX(D)) = 0$ since as before, $TX(D)$ is
isomorphic to $\Omega^1_X \otimes L_\kappa$, whose sections are part
of the Hodge decomposition of $H^1(X, \bbC_\kappa) = 0$.
This implies that the maps in the top left square are injective.


Notice that $N_D(D)$ is a line bundle on $D$ of degree $2g_D-2 = 4$,
with $h^0=2$ (since $N_D(D) = \Ocal_D(2D)$, which is the line bundle
corresponding to the linear system of quadrics).


\begin{prop}
  The map $H^0(D, TX(D)) \to H^0(D, N_D(D))$ is an isomorphism.
\end{prop}

\begin{proof}
  Since the (meromorphic) vector field $\drond/\drond x_\kappa$ is
  tangent to $X$ along $D$ (all equations of $X$ have zero derivative at $x_\kappa=0$),
  any vector field
  \[ \sum_{i \neq \kappa} \alpha_i x_i^2
  \frac 1 {x_\kappa} \drond / \drond x_\kappa \]
  is a regular section of $TX(D)$ over $D$. This shows that the
  composite
  \[ \Ocal_D(2D) \to TX(D)_{\upharpoonleft D} \to N_D(D)
  \text{ where the first map is }
  s \mapsto \frac s {x_\kappa} \frac{\drond}{\drond x_\kappa} \]
  is an isomorphism.
\end{proof}


This implies that $H^0(D, N_D(D))$ lifts to a 2-dimensional subspace
of $H^1(X, TX)$ which is the kernel of the infinitesimal period map.
The corresponding deformations can be expressed informally by writing
\[  dx_\kappa = \frac{1}{x_\kappa}
\sum_{i \neq \kappa} \eps_i x_i^2 \iff
d(x_\kappa^2) = 2 \sum_{i \neq \kappa} \eps_i x_i^2 \]

\subsection{The Jacobian ring}
\label{sec:campedelli-jacobian}

The classical theory of Griffiths \cites{Griffiths,Voisin}
details how the infinitesimal variation of Hodge
structure of a smooth hypersurface in $\mathbb P^n$ (inside the moduli
space of hypersurfaces) can be recovered from the Jacobian ring of its
equation. A similar construction describes the variation of periods for
complete intersections. We refer to \cite{Konno} for a detailed
exposition of the theory.

Let $Y$ be a complete intersection of four diagonal quadrics in
$\mathbb P^6$, and $\Lambda = |\mathfrak I_Y(2)| \simeq \mathbb P^3$
be the linear system
of quadrics through $Y$. Then $\mathbb P^6 \times \Lambda$ carries a
“universal” divisor $\tilde Y$, whose fibre over $\ell \in \Lambda$
is the quadric defined by $\ell$. If $\pi$ is the projection
$\tilde Y \to \mathbb P^6 \times \Lambda \to \mathbb P^6$,
the fibres of $\pi$ over $Y$ are projective spaces, and
$\pi$ is a $\mathbb P^2$-bundle outside $\pi^{-1}(Y)$.

\begin{thm}[Konno]
  The variation of Hodge structure of $H^8(\tilde Y)$ is canonically
  identified (up to a shift in gradings) with the variation of Hodge
  structure of $H^2(Y)$, which is canonically embedded as
  $H^2(R^6 \pi_\star \Z_{\tilde Y}) \subset H^8(\tilde Y, \Z)$.
\end{thm}

\index{anneau jacobien}
Consider the bigraded ring $S^{\bullet,\bullet} = \bbC[x_i; y_j]$
with 11 variables, which is generated by projective coordinates
on $\mathcal P^6$ ($i = 1 \dots 7$) and $\Lambda$ ($j=1,2,3,4$).
Suppose
\[ f = y_1 Q_1 + y_2 Q_2 + y_3 Q_3 + y_4 Q_4 \]
(where the $Q_j$'s are diagonal quadratic forms
$Q_j = \sum q_{ij} x_i^2$) is a parametrisation
of the linear system $\Lambda$: it is an equation of bidegree
$(2,1)$ of $\tilde Y \subset \mathbb P^6 \times \mathbb P^3$.
Then the variation of Hodge structure for $\tilde Y$ can be elegantly
described by the Jacobian ring of $f$ \cite{Green}.

For the sake of consistency with \cite{Konno}, we also introduce the
notation
\[ R_{p,q} = (S/J)^{2p-7q,p-4q} \simeq
H^0(\mathbb P^6 \times \Lambda, p \tilde Y + q K) \]
where $K$ is the canonical class of $\mathbb P^6 \times \Lambda$.

The Jacobian ideal of $\tilde Y$ is
\[ J^{\bullet, \bullet} = S \langle \partial_{x_i} f, \partial_{y_j} f \rangle
= S \Big\langle \sum_j q_{ij} x_i y_j, Q_j \Big\rangle \text. \]

\begin{thm}[Konno]
  \label{thm:hodge-eq-jacobian}
  If $p+q=8$, the space $H^p(\tilde Y, \Omega^q)_\mathrm{prim}$ is isomorphic to
  \[ R_{p+1, 1} = (S/J)^{2p-5,p-3} \simeq H^0((p+1)\tilde Y +
  K_{\mathbb P^6 \times \Lambda}) \text, \]
  using Griffiths's techniques, which associate to $P \in R_{p+1,1}$
  the meromorphic volume form $P \Omega / f^{p+1}$, where $\Omega$
  is a standard volume form on $\mathbb P^6 \times \mathbb
  P^3$ with values in $\Ocal(7,4)$.
\end{thm}

In this setting, the first-order deformations of the linear system
associated to $f$ are given by the elements of $R_{1,0} =
(S/J)^{2,1}$, while the infinitesimal period map of $Y$ is given by
the natural morphism
\[ R_{1,0} \to \Hom(R_{3+1,1}, R_{4+1,1}) \]
induced by multiplication in the ring $S/J$.

\begin{prop}
  The infinitesimal period map for $(X = Y/G, \kappa \in \hat G)$,
  which the identification given by theorem \ref{thm:hodge-eq-jacobian},
  is proportional to the natural map given by multiplication
  \[ R_{1,0}^{(0)} \to \Hom(R_{3+1,1}^{(\kappa)},
  R_{4+1,1}^{(\kappa)}) \]
  where $V^{(\kappa)}$ is the isotypic component of $V$ where $G$ acts
  by the character $\kappa$. Moreover, $R_{4,1}^{(\kappa)}$ is one-dimensional,
  generated by $x_\kappa$, so the map
  \[ x_\kappa: R_{1,0}^{(0)} \to R_{4+1,1}^{(\kappa)} \]
  also describes the infinitesimal $\kappa$-period map.
\end{prop}

An explicit description of $R_{1,0}^{(0)}$ is given by the $G$-invariant
part of the vector space, generated by the monomials $x_i x_j
  y_k$ in $S/J$, which is isomorphic to
\[ \frac{\langle x_i^2 y_j \rangle}
   {\big\langle Q_k y_j, \sum_j q_{ij} x_i^2 y_j \big\rangle} \]
which has dimension $6 = 28 - 16 - 7 + 1$, since the only non trivial
relation is $\sum Q_j y_j = \sum_{i,j} q_{ij} x_i^2 y_j$. The brackets here
denote the linear span of the elements they enclose.

The space $R_{5,1}^{(\kappa)} = (S/J)^{3,1}_{(\kappa)}$ is
\[ \frac{\langle  x_p x_q x_r y_j  \rangle}
   {\big\langle   x_\kappa Q_l y_j,
                  \sum_j q_{ij} x_i x_s x_t y_j
    \big\rangle}
   \qquad
   \begin{aligned}
      \text{where } p+q+r=\kappa \\
      \text{where } i+s+t=\kappa \\
   \end{aligned} \]
which can be written more suggestively as
\[ \frac{\langle x_i^2 x_\kappa y_j \rangle}
   {\big\langle  x_\kappa Q_l y_j,
                 \sum_j q_{ij} x_i^2 x_\kappa y_j,
                 \sum_j q_{\kappa j} x_\kappa x_s^2 y_j \text{ for } s \neq \kappa
    \big\rangle}
  \bigoplus
  \frac{\langle x_p x_q x_r y_j \rangle}
       {\big\langle \sum_j q_{pj} x_p x_q x_r y_j \big\rangle} \]
(in the second summand, $p+q+r=\kappa$, and none of $p,q,r$ equals
$\kappa$). The first space has dimension $4$ as we see below, and the second
space has dimension $8-6=2$ (note there are two possible triples $(p,q,r)$).
These dimensions coincide with the dimension of eigenspaces of
$s_\kappa$ on $H^{1,1}(Y, \bbC_\kappa)$: the (twisted) differential
form $\Omega / f^{p+1}$ is anti-invariant under $s_\kappa$, as well as
$x_\kappa$.

The image of multiplication by $x_\kappa$ lies inside the first space,
and by simplifying out $x_\kappa$, we observe that the annihilator of
$x_\kappa$ can be analysed by looking at the quotient morphism,
\begin{equation}
 \label{eq:target-dim-4}
 \frac{\langle x_i^2 y_j \rangle}
   {\big\langle Q_k y_j, \sum_j q_{ij} x_i^2 y_j \big\rangle}
   \longrightarrow
   \frac{\langle x_i^2 y_j \rangle}
        {\big\langle Q_k y_j, \sum_j q_{ij} x_i^2 y_j,
          \sum_j q_{\kappa j} x_i^2 y_j \big\rangle}
\end{equation}

\begin{prop}
  \label{prop:target-dim-4}
  If $Q$ is the matrix of a linear system defining a smooth (universal cover of)
  Campedelli surface, the target of the morphism in equation \ref{eq:target-dim-4}
  has dimension 4.
\end{prop}

The space of linear combinations $\sum a_{ij} x_i^2 y_j$ is identified
with the space of matrices $(a_{ij})$ with size $7 \times 4$.
The subspace $\langle Q_k y_j \rangle$ is then the space of products $MQ$
where $M \in \bbC^{4 \times 4}$. Since $Q$ defines a surface,
it must have full rank, hence the space of products
$MQ$ has dimension 16. We use an intermediate lemma to prove proposition
\ref{prop:target-dim-4}.

\begin{lemma}
  If $Q = (q_{ij})$ is the matrix of the linear system defining a smooth
  Campedelli surface, written in the standard basis $(x_i^2)_{i \in \hat G}$,
  no set of four columns of $Q$ are linearly dependent.
\end{lemma}

\begin{proof}
  If there were such a set, then some element of the linear system
  would be a quadric of rank three: then there would be three
  linearly dependant monomials $x_i^2$ in the bicanonical linear
  system. The configuration of lines describing the Campedelli
  surface has a triple point, which creates a singularity.
\end{proof}


\begin{proof}[Proof of proposition \ref{prop:target-dim-4}]
  We work with the spaces of matrices described above:
  a matrix $MQ$ lies in $\big\langle \sum_j q_{\kappa j} x_i^2 y_j, \sum_j
    q_{\kappa j} x_i^2 y_j \big\rangle$, if $Mx$ is a linear combination of $x$
  and $K$ for any column $x$ of $Q$ ($K$ being the column $\kappa$).
  Then $K$ is an eigenvector of $M$, and we need to analyse
  its action on $\bbC^4/K$.

  The six columns of $Q$ (except $K$) should define 3 different
  eigenvectors of $M$ in $\bbC^4/K$. Using the previous lemma, we know
  that no three of them are linearly dependent. By standard arguments,
  this forces $M$ to act as a scalar multiplication on $\bbC^4/K$.
  Conversely, if $M$ acts as a scalar on $\bbC^4/K$, then $MQ$ lies in
  the given subspace: the space of such $MQ$ has dimension 5.

  The dimension of $\langle Q_k y_j, \sum_j q_{ij} x_i^2 y_j, \sum_j
    q_{\kappa j} x_i^2 y_j \rangle$ is then $16 + (7+7-1) - 5 = 24$, hence
  the result.
\end{proof}

\subsection{The period map around Campedelli surfaces with double points}
\label{ssec:period-double-pt}

We will need to examine the behaviour of the period map around a
Campedelli surface with a double point, following the method described by
C. Voisin in \cite{Voisin-4fold}. We are interested in the
following particular situation: let $\alpha$, $\beta$, $\gamma$
be a basis of $(\Z/2)^3$, and $\kappa = \alpha+\beta+\gamma$.

Let $X$ be a Campedelli surface (the associated configuration of 7
lines in the bicanonical plane should have at worst triple points) and
assume that the lines $D_\alpha$, $D_\beta$ and $D_\gamma$ are
concurrent in $|2K_X|^\vee$: then $x_\alpha^2 + x_\beta^2 +
x_\gamma^2$ belongs to the defining linear system, and there exists a
point $R$ such that $x_\alpha = x_\beta = x_\gamma = 0$. Note that
there exists an element of $\pi_1(X)$ whose action on coordinates
reverses the sign of $x_i$ for $i = \alpha, \beta, \gamma,
\kappa$. This implies that $R$ is fixed under the involution
$s_\kappa$. In the following, we denote $(\kappa, \alpha, \beta,
\gamma)$ by numbers 4, 5, 6, 7.

Note that $R$ cannot have a fourth vanishing coordinate, since this
would mean that a fourth branching line is concurrent with $D_\alpha$,
$D_\beta$, $D_\gamma$.

By the Picard-Lefschetz formula, the monodromy of the Gauss-Manin connection
around the hypersurface of configurations with three concurrent lines
has order two, and there cannot exist (even locally) a continuous
choice of trivialisation of the lattices $H^2(X_\kappa(t), \Z)$
(hence neither of $H^2(X_t,\Z_-)$) in any neighbourhood of $X$ in the moduli space.

However, taking the double cover ramified along this hypersurface
trivialises the monodromy and allows to get a simultaneous resolution
of singularities \cite{Atiyah-resolution}, since the hypersurface
parametrises surfaces acquiring double points. It becomes possible
to define a local period map with values
in a type IV domain $\mathcal D$, which is then holomorphic.

Let $t$ be a complex parameter and consider the deformation of
$X$ given by the matrices
\[ \begin{pmatrix}
  a_1 & a_2 & a_3 & a_4 & a_5 & a_6 & a_7 \\
  b_1 & b_2 & b_3 & b_4 & b_5 & b_6 & b_7 \\
  c_1 & c_2 & c_3 & c_4 & c_5 & c_6 & c_7 \\
  0 & 0 & 0 & -t^2 & 1 & 1 & 1 \\
\end{pmatrix} \]
which corresponds to the double cover of a one-parameter family which
is transverse to the hypersurface defined by the minor of the first
columns.

We make the following regularity hypothesis:
\[ \text{the determinant } \begin{vmatrix}
  a_1 & a_2 & a_3 \\
  b_1 & b_2 & b_3 \\
  c_1 & c_2 & c_3 \\
\end{vmatrix} \text{ does not vanish.} \]

\begin{lemma}
  This hypothesis is equivalent to requiring that $R$ has no fourth
  vanishing coordinate, or to the assumption that no set of four lines
  are concurrent. \qed
\end{lemma}

As before, the period map of this family can be identified (locally)
with the period map of the associated family of hypersurfaces of
bidegree $(1,2)$ in $\mathbb P^3 \times \mathbb P^6$, with equations
$F_t = \sum u_j q_j$ where $q_1 = \sum a_i x_i^2$ and so on.

The singular locus of $F_t$ is defined by the equations $q_j$,
and the vanishing of the Jacobian matrix of $(q_1, q_2, q_3, q_4)$
on $(u_i)$. For generic choices of $q_1$, $q_2$, $q_3$ (described
by our regularity hypothesis), there is an isolated fixed point at
$([0:0:0:1], R)$. Choosing $u_4 = 1$, $x_4 = 1$ as a local chart, the
equation of the family is given by
\[ u_1 q_1 + u_2 q_2 + u_3 q_3 - t^2 + x_5^2 + x_6^2 + x_7^2 = 0 \]
which defines an isolated ordinary double point: in other words,
$\set{q_1, q_2, q_3, x_5, x_6, x_7}$ is a
set of local coordinates on $\mathbb P^6$. The
action of $s_\kappa$ lifts to these local coordinates by changing
signs in $x_5$, $x_6$, $x_7$.

If $\tilde X_t$ is the complete intersection of quadrics associated to the
parameter $t$, and $Y_t$ is the associated $8$-fold in $\mathbb P^3
\times \mathbb P^6$, the Hodge structure on $H^8(Y_t)$ is described
by $H^{3,5}(Y_t) = \mathbb C \tilde \omega$, where
\[ \tilde \omega = \Res_{Y_t} \frac{x_\kappa \Omega}{F_t^4} \]
where $\Omega$ is a non-vanishing generator of $\Omega^9(4,7)$.

By the works of Griffiths and Konno, we know that the variation of
Hodge structure along this family will be described by the following
type $(4,4)$ form:
\[ \tilde \omega' = \Res_{Y_t} \frac{t u_4 x_4^2 x_\kappa \Omega}{F_t^5} \]
up to some multiplicative constant.

In the local chart described above $(u_4 = x_4 = 1)$, we can write
\[ \tilde \omega' = \Res_{Y_t} \frac{t x_\kappa \Omega}{F_t^5} \]
where $\Omega$ is the canonical holomorphic volume form
on $\bbC^3 \times \bbC^6$.

If $\kappa$ were not $\alpha+\beta+\gamma$, then $R$ would not be a
fixed point of the involution: the regularity of the period map
follows from a blow-up at the singular point of the family $Y_t$.
In the affine chart where $x_i = tX_i$, $q_i = t Q_i$, $u_i = tU_i$
(note that $x_\kappa$ is nonzero):
\[ \tilde \omega'(t) = \Res_{Y_t}
\frac{t^{10} x_\kappa \Omega_\mathrm{bl}}{t^{10} (F^\mathrm{bl}_t)^5} \]
\[ \text{where } F_t^\mathrm{bl} =
U_1 Q_1 + U_2 Q_2 + U_3 Q_3 + X_5^2 + X_6^2 + X_7^2 - 1 \]
A change of variables turns $F_t$ into $\sum Y_i^2 - 1$,
and the period over the corresponding vanishing sphere
has a finite limit (since the given residue is the way of expressing
the primitive cohomology of the quadric).

We now need to know how this differential form integrates in the
case where $R$ is invariant under $s_\kappa$: in this case some
cohomology class is anti-invariant under $s_\kappa$ and its period
vanishes for $t=0$. Such a situation has already been studied by Horikawa
for Enriques surfaces \cite{Horikawa2}.

The key fact is that the projective quadric in $\mathbb P^9$, which is
the normal cone to the singular point in the 9-fold $\mathcal Y =
(Y_t)$, has two base classes in $H^{4,4}$, corresponding to classes
of maximal isotropic subspaces, $\sigma_1$ and $\sigma_2$, which
form a hyperbolic plane. They are such that $h^4 = \sigma_1 + \sigma_2$
is the class of a 4-dimensional linear section, and $\rho = \sigma_1 -
\sigma_2$ is the class of the real sphere $\mathbb S^8$. Note that
$\rho^2 = -2$.

Since $s_\kappa$ acts by changing signs of three coordinates, it
exchanges $\sigma_1$ and $\sigma_2$ (it is not in $SO_{10}$), hence
the real sphere is anti-invariant under this transformation
($s_\kappa$ reverses the orientation of the associated real manifold),
and by the same argument the associated period has a nonzero
derivative at $t=0$.

\section{Todorov surfaces and double covers of Enriques surfaces}

Todorov surfaces were introduced by Todorov to give
examples of surfaces whose infinitesimal period map is nontrivial and
not injective \cite{Todorov}. A systematic study of these surfaces is done in
\cite{Morrison}, which we quote for most of the properties stated
below.

Campedelli surfaces are double covers of Enriques surfaces, and their
étale double covers (whose periods are what we actually study) are
Todorov surfaces. In section \ref{ssec:todorov-enriques} we carry out
a study of basic properties of double covers of Enriques surfaces,
similar to \cite{Morrison}. If $X \to S$ is such a double cover,
the geometry of $S$ usually gives good information about the
transcendental part of $H^2(S, \Z_-)$, and in
section \ref{ssec:index-comp} we compute a formula giving
the index of the embedding $H^2(S, \Z_-) \to H^2(X, \Z_\kappa)$
(theorem \ref{thm:index-pullback}) between (twisted) cohomology
lattices. It will be used to compare the global period maps for $S$
and $X$.

\subsection{Classical Todorov surfaces}

\begin{defn}[Todorov surface]
  \index{surface!de Todorov}
  A \emph{Todorov surface} is a surface $Z$ with canonical singularities, and an
  involution $j$, such that $S = Z/j$ is a K3 surface with rational
  double points and $\chi(\Ocal_Z) = 2$.
\end{defn}

If $\tilde Z$ is the minimal resolution of $Z$, the natural morphism
$\tilde j: \tilde Z \times_{Z,j} Z \to \tilde Z$ lifts $j$ to an involution
of $\tilde Z$. The quotient $\Sigma = \tilde Z / \tilde j$ is again a
K3 surface, which is a partial desingularisation of $S$.

\begin{defn}[Fundamental invariants]
  Let $(Z,j)$ be a Todorov surface. Then $\Sigma$ as above is a nodal
  K3 surface: the lattice generated by its nodes has index $2^\alpha$
  inside its primitive saturation in $H^2(\Sigma, \Z)$. The number of
  nodes of $\Sigma$ is denoted by $k$. The \emph{fundamental
    invariants} of $Z$ are the integers $(\alpha,k)$.
\end{defn}

\begin{thm}[Todorov]
  Several topological invariants of $\tilde Z$ can be calculated using
  the fundamental invariants. The order of the 2-torsion subgroup of
  $H^2(\tilde Z, \Z)$ is $2^\alpha$. The divisorial part of the
  ramification locus $B \subset \Sigma$ satisfies $B^2 = 2k-16$. The
  integer $k$ can also be expressed as $c_1(\tilde Z)^2 + 8$.
\end{thm}


\begin{prop}
  If $X$ is a Campedelli surface without $(-2)$-curves, $X_\kappa$ is
  a Todorov surface with invariants $(2,12)$.
\end{prop}

\begin{proof}
  As seen in section \ref{sec:quotient}, the involution $s_\kappa$ on
  $X$ lifts to $X_\kappa$, which is a projective surface with rational
  double points, and the quotient $X_\kappa / s_\kappa$ is the
  universal cover of the Enriques surface $X/s_\kappa$. We recover
  $\alpha=2$ from the order of the 2-torsion subgroup of
  $H^2(X_\kappa,\Z)$ (which is identified with the character group of
  $\pi_1(X_\kappa)$), and $k=12$ from the number of double points of
  $T_\kappa = X_\kappa / s_\kappa$ (we know that $S_\kappa$ has six
  double points).
\end{proof}


Another invariant of Todorov surfaces is a natural sublattice of the
Picard group of the underlying K3 surface:
\begin{defn}[Todorov lattice \cite{Morrison}]
  Let $Z$ be a smooth Todorov surface and $\Sigma = Z/j$ be the
  associated nodal K3 surface. The resolution of singularities of
  $\Sigma$ is denoted by $\tilde \Sigma$: each double point of $\Sigma$
  is resolved to a $(-2)$-curve $E_i$ on $\tilde \Sigma$.

  The \emph{Todorov lattice} associated to $Z$ is the primitive
  saturation $L_T(Z)$ of the sublattice $\pair{B, E_i}$ of
  $H^2(\tilde \Sigma, \Z)$ generated by the ramification divisor $B$
  and the classes $E_i$.
\end{defn}


The Todorov lattice is the Picard lattice of K3 surfaces associated to
generic Todorov surfaces with given invariants $(\alpha,k)$. Such K3
surfaces are parametrised by a period domain: a moduli space for
Todorov surfaces can thus be constructed using their period map
\cite{Morrison}.






\subsection{Todorov-Enriques surfaces}
\label{ssec:todorov-enriques}

We define a class of surfaces inspired by the definition of Todorov
surfaces. This class includes Keum-Naie surfaces \cite{Naie} and the
construction by Mendes Lopes and Pardini \cite{MLP} of surfaces such
that $p_g=q=0$. We gather here a collection of results which fit
Campedelli surfaces in this class.

\begin{defn}
  \label{defn:todorov-enriques}
  A \emph{Todorov-Enriques} surface is a pair $(X,j)$ where $X$ is a
  canonical surface with an involution $j: X \to X$, such that
  $\chi(\Ocal_X) = 1$, and $X/j$ is an Enriques surface with at
  worst rational double points.
\end{defn}

Replacing if necessary $X$ by its minimal desingularisation, we will
assume that $X$ is smooth, and that the double points of $X/j$ come from
the isolated fixed points of $j$. The smooth Enriques surface obtained
by blowing-up the double points of $X/j$ is denoted by $S$.

The double points of $X/j$ define $(-2)$-classes $E_1, \ldots, E_k$ in
$\Pic(S)$. The \emph{invariants} of $(X,j)$ are $k$ (the number of
fixed points of $j$, which is also the number of double points on
$X/j$ and $\alpha$, which is the dimension of the kernel of
the natural map
\[ \Z/2\Z \pair{E_1, \dots, E_k} \to \frac{\Pic S}{2 \Pic S} \]
(the space of “even sets” made of $E_1, \dots, E_k$).


\begin{prop}
  \label{prop:todorov-selfint-B}
  The ramification divisor $B \subset X/j$ has self-intersection $2K_X^2 = 2k-8$,
  and $h^0(S, \Ocal_S(B)) = k-3$.
\end{prop}

\begin{proof}
  Since $K_X = K_S + \pi^\ast B/2$, $K_X^2 = B^2/2$.
  Now, if $X_b$ is the blowup of the $k$ fixed points of $j$,
  \[ e(X_b) - e(B) - 2k = 2(e(S) - e(B) - 2k) \]
  hence $e(X_b) = e(X) + k = 2 \cdot 12 + B^2 - 2k$. By Noether's formula,
  $12 = K_X^2 + e(X)$, giving $12 = 24 + 3B^2/2 - 3k$ and the
  equality $B^2=2k-8$.
\end{proof}

\begin{prop}
  The inequality $2\alpha \leq k \leq \alpha + 5$ holds, or
  equivalently $k-5 \leq \alpha \leq k/2$.
\end{prop}

\begin{proof}
  Let $N_S$ be the primitive saturation of the lattice generated by
  the nodal classes $E_i$ in $\Pic(S)$ (which is a rank 10 unimodular
  lattice).

  The double point lattice $N_S$ of $S$ has rank $k$ and discriminant
  $2^{k-2\alpha}$ (hence $k \geq 2\alpha$), but its orthogonal
  complement has rank $10-k$. Since the discriminant group of
  $N_S^\perp$ is a $\mathbb F_2$-vector space of rank at most
  $10-k$, which is isomorphic to the discriminant group of $N_S$,
  $k-2\alpha \leq 10-k$.
\end{proof}

Using this inequality (and the fact that $K_X^2 = k-4 > 0$),
the possible values of $(\alpha, k)$ are
\[ (0,5) \quad (1,5) \quad (2,5) \quad (1,6) \quad (2,6) \quad (3,6)
\quad (2,7) \quad (3,7) \quad (3,8) \quad (4,8) \]

Remember that an even set of nodes on an Enriques surface is made of 4
or 8 nodes. If $k=5$, the existence of two distinct even sets of nodes
would imply the existence of an even set of two nodes, which is
impossible.

If $k=8$, the universal cover of $S$ is necessarily a Kummer surface:
the study of even sets on a Kummer surface (see \cite{BHPV} for
example) tells us that in this case $\alpha=4$.

The pair $(3,6)$ is also impossible: since any even set of nodes has
four elements, two distinct even sets must share exactly two nodes,
in other words, the complements of distinct even sets are disjoint,
so there cannot be more than three non trivial even sets.

The possible values of $(\alpha, k)$ are now
\[ (0,5) \quad (1,5) \quad (1,6) \quad (2,6) \quad (2,7)
\quad (3,7) \quad (4,8) \]

\begin{prop}
  The 2-torsion subgroup of $\Pic(X)$ is an extension of
  $(\Z/2\Z)^\alpha$ by the 2-torsion subgroup of $H^2(S,\Z)$
  (which is isomorphic to $\Z/2\Z$).
\end{prop}

\begin{proof}
  Let $\pi: X_b \to S$ be the projection, and $U$ be the
  complement in $S$ of the image of the fixed locus $F \subset X_b$ of $j$, and $u: U
  \into S$ be the standard inclusion. Consider the extensions
  of sheaves
  \[ \begin{aligned}
    e_U: &0 \to \Z/2\Z \to \pi_\star \Z/2\Z \to \Z/2\Z \to 0 \\
    e_S: &0 \to \Z/2\Z \to \pi_\star \Z/2\Z \to u_! \Z/2\Z \to 0 \\
    \chi: &0 \to u_! \Z/2\Z \to \Z/2\Z \to i_\ast \Z/2\Z \to 0 \\
  \end{aligned} \]
  which determine $e_U \in H^1(U, \Z/2)$, $e_S \in \Ext^1_S(u_! \Z/2, \Z/2)$
  and $\chi \in \Ext^1_S(i_\ast \Z/2, u_! \Z/2)$.


  The exact sequence of sheaves
  \[ 0 \to \Z/2\Z \to \pi_\star \Z/2\Z \to u_! \Z/2\Z \to 0 \]
  induces an exact sequence
  \[ 0 \to \Z/2\Z \simeq H^1(S, \Z/2\Z)
  \to H^1(X_b, \Z/2\Z) \to H^1_c(S \setminus F, \Z/2\Z)
  \to H^2(S, \Z/2\Z) \text. \]
  which expresses $H^1(X_b, \Z/2\Z)$ as an extension of
  $\ker e_S \subset H^1_c(S \setminus F, \Z/2\Z)$
  by $\Z/2$.

  The relative cohomology exact sequence now tells us that
  \[ H^0(S, \Z/2\Z) \to H^0(F, \Z/2\Z) \xrightarrow{\chi}
  H^1_c(S \setminus F, \Z/2\Z) \to H^1(S, \Z/2\Z) \to H^1(F, \Z/2\Z) \]
  is exact, but the last map being injective, $H^1_c(S \setminus F,
  \Z/2\Z)$ is isomorphic to $(\Z/2\Z)^f / e$, where $f$ is the number of
  components of $F$ and $e$ is the sum of elements of $F$. We must now compute
  \[ \ker e_S = \ker(e_S \chi: H^0(F, \Z/2) \to H^2(S, \Z/2))/e \]

  But $e_S \chi$ is an element of $\Ext^2_S(i_\ast \Z/2, \Z/2) \simeq
  H^0(\Extrond^2(i_\ast \Z/2, \Z/2))$ which is totally determined by
  looking at the sections of $\Extrond^1(i_\ast \Z/2, u_! \Z/2)$ and
  $\Extrond^1(u_!\Z/2, \Z/2)$ determined by $e_S$ and $\chi$. It is
  easy to check that it is exactly the Gysin map.

  Its kernel is identified with the code of even sets of nodes
  in $(\Z/2\Z)^k$.
\end{proof}

This shows that $(1,5)$ is impossible, since it would give a numerical
Godeaux surface such that $H^1(X, \Z/2\Z)$ contains $(\Z/2)^2$, which is
impossible \cite{Miyaoka}*{sec.~3}.

\begin{prop}
  A Todorov-Enriques surface admits a canonical étale double cover,
  induced by the universal cover of the quotient Enriques surface.
  This double cover is a Todorov surface. \qed
\end{prop}

Since Todorov surfaces are known to satisfy $q=0$ \cite{Morrison}, Todorov-Enriques
surfaces satisfy $p_g=q=0$.


If $X$ is a Todorov-Enriques surface, and $(\tilde \alpha, \tilde k)$
are the invariants of a Todorov surface $Y$ which is an étale double
cover of $X$, then $\tilde k = 2k$, and we have the following
possibilities.
\[ \begin{matrix}
(\alpha, k) &
(0,5) & (1,6) & (2,6) &
(2,7) & (3,7) & (4,8) \\
(\tilde \alpha, \tilde k) &
(0,10); (1,10) & (1,12); (2,12) & (2,12) &
(3,14) & (3,14) & (5,16) \\
\end{matrix} \]

\index{surface!de Keum-Naie}
In \cite{Naie}, D. Naie constructed examples of such surfaces for the
invariants $(0,5)$, $(1,6)$, $(2,6)$, $(2,7)$ and $(4,8)$, starting
with an Enriques surface with 8 nodes (the value of $\alpha$
can be deduced from the observation that the 2-torsion groups have order
$2$ for $k=5$, $4$ or $8$ for $k=6$, $2^3$ for $k=7$). The general
description of surfaces of type $(2,7)$ can be found in \cite{MLP}.
The case $(3,7)$ is actually impossible \cite{MLP}*{4.4}: the Enriques
surface would be realised as a surface with seven nodes in $\mathbb
P^3$, the nodes being aligned like the seven points of $\mathbb
P^2(\mathbb F_2)$, which is impossible.

\begin{prop}
  \label{prop:campedelli-are-todorov-enriques}
  The Campedelli surfaces are examples of Todorov-Enriques surfaces
  with invariants $(\alpha,k) = (2,6)$.
\end{prop}

\begin{proof}
  Let $(X, \kappa)$ be a marked (smooth) Campedelli surface and $S =
  X/s_\kappa$ be the associated Enriques surface with six nodes. Here
  $B$ is the (reduced) image of $D_\kappa$ in $S$ ($B^2 =
  2k-8 = 4$ by proposition \ref{prop:todorov-selfint-B}). This gives $k=6$,
  and explained after definition \ref{defn:todorov-enriques}, $\alpha=2$
  means there are three sets of even nodes among the six nodes on $S$.

  Note that any twisted canonical divisor $D_i \neq D_\kappa$ goes through the
  two fixed points $\set{e_i, e'_i} = D_i \cap D_{i+\kappa}$ and no other fixed point
  (except those of $D_\kappa$): let $B_i$ be the image of $2D_i$ in $S$, which
  is a divisor with generic multiplicity two, which pulls back to
  $2K_X$, $B_i$ is an element of $|2B|$.

  There are two rational curves $E_i$, $E'_i$ in the desingularisation
  $\tilde S$ corresponding to the points $e_i$ and $e'_i$. The
  decomposition of $B_i$ into irreducible components is
  $B'_i + S_i + S'_i$ (since $S_i$ pulls back to $2E_i$),
  and $B'_i = 2B-S_i-S'_i$ is divisible by two. This implies that the complementary
  set of four rational curves on $\tilde S$ is even, since the existence of the ramified
  double cover $X \to S$ requires that the sum of $B$ and the six
  rational $(-2)$-curves is divisible by two.

  The three possible pairs of fixed points $D_i \cap D_{i+\kappa}$ provide
  the three required even sets of nodes.
\end{proof}

It should be noted that numerical Campedelli surfaces may be also
constructed using an Enriques surface with invariants
$(\alpha,k) = (1,6)$, giving fundamental groups $\Z/4 \times \Z/2$
or $\Z/2$. We hope to study these families of Enriques surfaces in
detail in a future work.

\subsection{Embeddings of cohomology lattices for double covers}
\label{ssec:index-comp}

In order to compare the period map of Enriques surfaces with the
actual period map of a family of covering Todorov-Enriques surfaces,
we give a formula computing the index of the map $H^2(S,\Z_-) \to
H^2(X,\Z_\kappa)^G_\mathrm{num}$, when $f: X \to S$ is a degree two cover of an Enriques
surface, $G$ is the group generated by the associated involution of
$X$, and $\Z_\kappa = f^\star \Z_-$. Then $f^\ast: H^2(S, \Z_-) \to H^2(X, \Z_\kappa)_\mathrm{num}$
is a morphism of quadratic lattices, which multiplies the intersection
form by two. The notation $H_\mathrm{num}$ denotes the quotient of an
abelian group $H$ by its torsion subgroup (cohomology classes up to
numerical equivalence).

We make the following assumptions: the ramification locus of $f$ is a
disjoint union of $\rho$ smooth curves ($\rho$ being a positive
integer), and any étale double cover $X'$ of $X$ is regular ($q(X')=0$).
The only needed consequence of $S$
being an Enriques surface is that $H^3(S, \Z_-) = 0$.

\index{suite spectrale!de cohomologie équivariante}
The computation is done using spectral sequences for
$G$-equivariant cohomology: here $H^k_G(X, \bullet)$ can
be understood as the $k$-th derived functor of
$\mathcal F \mapsto \Gamma(X, \mathcal F)^G$.
The Borel diagram for the action of $G$ is
\[ \xymatrix{
  EG \ar[d] & X \times EG \ar[l] \ar[r] \ar[d]^{\tilde f} & X \ar[d]^{f}  \\
  BG & [X/G] \ar[l]^{\pi} \ar[r]_\gamma & S \\
} \]
where $BG$ is the classifying space of $G$ and $EG$ the universal
$G$-bundle over $BG$. Then $H^\bullet_G(X, \Z_\kappa) = H^\bullet([X/G], \Z_\kappa)$
can be calculated by the Leray spectral sequences for $\pi$ and $\gamma$.
\[ {E'}_2^{p,q} = H^p(G, H^q(X, \Z_\kappa))
\implies H^{p+q}_G(X, \Z_\kappa) \]
\[ {E"}_2^{p,q} = H^p(S, R^q \gamma_\star \Z_\kappa)
\implies H^{p+q}_G(X, \Z_\kappa) \]

The computation of sheaves $R^q \gamma_\star \Z_\kappa$ corresponds to
cohomology groups of $\Z/2\Z$:
\[ \gamma_\star \Z_\kappa \simeq \Z_- \qquad
R^{2k+1} \gamma_\star \Z_\kappa = 0 \quad \text{and} \quad
R^{2k+2} \gamma_\star \Z_\kappa \simeq (\Z/2\Z)_{\restriction F}
\qquad (k \geq 0) \text. \]

The map $f^\star$ is decomposed as follows:
\[ f^\star: H^2(S, \Z_-) \xrightarrow{\gamma^\star}
H^2_G(X,\Z_\kappa) \xrightarrow{\bar f^\star} H^2(X, \Z_\kappa)^G \]
By eliminating torsion, we obtain two maps of free abelian groups
\[ \begin{aligned}
  H^2(S, \Z_-)_\mathrm{num} &\to H^2_G(X,\Z_\kappa)_\mathrm{num}
  \text{ (of index $2^{N_1}$)} \\
  \text{and } H^2_G(X,\Z_\kappa)_\mathrm{num} &\to
  H^2(X, \Z_\kappa)^G_\mathrm{num} \text{ (of index $2^{N_2}$).} \\
\end{aligned} \]
The index of $f^\star$ is thus expressed as the $(N_1+N_2)$-th power of two.
We are actually going to prove the following formulas:
\begin{prop}
  The integer $N_1$ (see proposition \ref{prop:compute-index-N1}) is equal to
  \begin{equation}
    \label{eqn:compute-N1}
    \rho - \ell_2(H^2_G(X,\Z_\kappa)) + \ell_2(H^2(S,\Z_-)) \text.
  \end{equation}
  The integer $N_2$ (see propositions \ref{prop:compute-index-N2-a} and
  \ref{prop:compute-index-N2-b}) is equal to
  \begin{equation}
    \label{eqn:compute-N2}
    \ell_2(H^2_G(X,\Z_\kappa)) - \ell_2(H^2(X, \Z_\kappa)^G) - \eps
  \end{equation}
  where $\eps$ is 1 if $\Z_\kappa$ is trivial, 0 otherwise.
\end{prop}
Here $\ell_2$ denotes the length of the 2-adic torsion subgroup.

\begin{thm} \label{thm:index-pullback}
  Under the hypotheses above, the map of lattices associated to
  $f^\ast$ has index
  \[ 2^{\ell_2(H^2(S, \Z_-)) - \ell_2(H^2(X, \Z_\kappa)^G) + \rho} \]
  if $\Z_\kappa$ is non trivial,
  \[ 2^{\ell_2(H^2(S, \Z_-)) - \ell_2(H^2(X, \Z_\kappa)^G) + \rho - 1} \]
  otherwise.
\end{thm}

\subsubsection{Second spectral sequence}

Since ${E"}_2^{p,q} = 0$ for any odd $q$, the differential $d_2$
is zero. The remaining differential is $d_3^{(0,2)}: H^0(F, \Z/2\Z)
\to H^3(S, \Z_-) = 0$.
The resulting filtration on $H^2_G(X,\Z_\kappa)$ is graded by
\[ \mathrm{Gr}^2 = H^2(S,\Z_-) \qquad \mathrm{Gr}^1 = 0
\qquad \mathrm{Gr}^0 = H^0(F, \Z/2\Z) \simeq (\Z/2\Z)^\rho \]

\begin{prop}
  \label{prop:compute-index-N1}
  The map $H^2(S,\Z_-) \to H^2_G(X,\Z_\kappa)$ is injective, its cokernel is
  a $\Z/2\Z$-module of rank $\rho$, and
  $N_1 = \rho - \ell_2(H^2_G(X,\Z_\kappa)) + \ell_2(H^2(S,\Z_-))$, as stated
  by formula \ref{eqn:compute-N1}.
\end{prop}

\begin{proof}
  Let $T(S)$, $T_G(X)$ be the torsion groups of $H^2(S,\Z_-)$ and
  $H^2_G(X,\Z_\kappa)$. The natural map $H^2(S,\Z_-)/T(S) \to H^2_G(X,\Z_\kappa)/T(S)$
  has again cokernel $\simeq (\Z/2)^\rho$. This cokernel is an
  extension of the 2-torsion of the target (with length
  $\ell_2(H^2_G(X,\Z_\kappa)) - \ell_2(H^2(S,\Z_\kappa))$), and the contribution
  from the torsion-free part (a group of order $2^\kappa$).
  We get the equation stated above.
\end{proof}

\subsubsection{First spectral sequence}

\begin{lemma}
  Let $\Z_\kappa$ be a $\Z$-local system on an algebraic surface $X$ with no
  irregular étale double cover. Then $H^1(X, \Z_\kappa) = 0$ if $\Z_\kappa = \Z_X$,
  $\Z/2\Z$ if $\Z_\kappa$ is non trivial.
\end{lemma}

\begin{proof}
  If $\Z_\kappa = \Z_X$, then
  $H^1(X, \Z) = \Hom(H_1(X), \Z)$ by the universal coefficient
  theorem, and this is torsion-free, but
  since $H^1(X, \bbC) = 0$, it must be zero.

  If $\Z_\kappa$ is non trivial, then $\Z_\kappa$ is given by a character $\kappa: \pi_1(X)
  \to \Z/2\Z$. There is a canonical double cover $X_\kappa \to X$
  and a short exact sequence:
  \[ 0 \to \Z \xrightarrow 2 \Z \to H^1(X, \Z_\kappa) \to H^1(X_\kappa, \Z)
   = 0 \text. \]
\end{proof}




Suppose that $\Z_\kappa = \Z_X$. Then $H^p(G, H^q(X, \Z_\kappa))$ is zero
for $q=1$, and for $(p,q) = (2p'+1, 0)$. The graded parts $\mathrm{Gr}_I^p$
of $H^2_G(X, \Z)$ are
\[ \begin{aligned}
  \mathrm{Gr}^2 &= H^2(G, \Z) \simeq \Z/2\Z \\
  \mathrm{Gr}^1 &= 0 \\
  \mathrm{Gr}^0 &= H^0(G, H^2(X, \Z)) \\
\end{aligned} \]

\begin{prop}
  \label{prop:compute-index-N2-a}
  When $\Z_\kappa$ is trivial, the map $H^2_G(X,\Z) \to H^2(X, \Z)^G$ is surjective and induces an
  isomorphism between the torsion-free quotients (i.e. $N_2=0$). Moreover the following
  relation holds:
  \[ \ell_2(H^2_G(X,\Z)) = \ell_2(H^2(X,\Z)^G) + 1 \text. \]
\end{prop}

Suppose now that $\Z_\kappa \neq \Z_X$: then $H^0(X,\Z_\kappa)=0$ and $H^1(X,\Z_\kappa) = \Z/2$. The graded
parts $\mathrm{Gr}_I^p$ of $H^2_G(X, \Z_\kappa)$ are
\[ \begin{aligned}
  \mathrm{Gr}^2 &= 0 \\
  \mathrm{Gr}^1 &= H^1(G, \Z/2\Z) \simeq \Z/2\Z \\
  \mathrm{Gr}^0 &= \ker: H^0(G, H^2(X, \Z_\kappa)) \to H^2(G, H^1(X, \Z_\kappa)) \\
\end{aligned} \]
hence there is an exact sequence
\[ 0 \to \Z/2\Z \to H^2_G(X, \Z_\kappa) \to H^2(X, \Z_\kappa)^G \to \Z/2\Z \to 0
\text. \]

\begin{prop}
  \label{prop:compute-index-N2-b}
  When $\Z_\kappa$ is not trivial, $N_2 = \ell_2(H^2_G(X, \Z_\kappa)) - \ell_2(H^2(X, \Z_\kappa)^G)$
  as in formula \ref{eqn:compute-N2}.
\end{prop}

\section{Enriques surfaces with a
  \texorpdfstring{$D_{1,6}$}{D(1,6)} polarisation}

The Enriques surfaces appearing as quotients of Campedelli surfaces
have (in the generic case) six ordinary double points: these surfaces
are exactly the Enriques surfaces $S$ containing a $D_6$ sublattice in
$H^{1,1}(S, \Z_-)$. In other words their six nodes can form three even
sets, as seen in proposition \ref{prop:campedelli-are-todorov-enriques}.

This allows to determine the generic transdencental part of $H^2(S,
\Z_-)$ (section \ref{ssec:tanscendental-for-d6}), which is the
orthogonal complement of the $D_6$ lattice. It is isomorphic to
$\Z^2(2) \oplus \Z^4(-1)$ (proposition \ref{prop:gen-transcendental-lattice}).

These Enriques surfaces can be defined as bidouble covers of the plane
(proposition \ref{prop:Enriques-galois-cover}), whose ramification
locus consists of three pairs of lines. We will see that a Cremona
transformation can map these pairs to another configuration
(which is not projectively equivalent in generic cases),
defining a birationally equivalent Enriques surface
(proposition \ref{prop:cremona-involution}).

\subsection{Linear systems and geometry}
\label{ssec:d16-enriques-geometry}

Let $D_{1,6}$ be the index 2 sublattice of the standard Lorentzian
lattice $\Z^{1,6}$, containing vectors whose sum of coordinates is
even. If $\pair{e_0; e_1, \dots, e_6}$ is the canonical basis of
$\Z^{1,6}$, we distinguish a norm 4 vector $2e_0$ and six mutually
orthogonal $(-2)$-vectors $e_1 \pm e_2$, $e_3 \pm e_4$, $e_5 \pm
e_6$.

\begin{defn}
  A $D_{1,6}$-polarised Enriques surface is an Enriques surface whose
  Picard group contains a primitively embedded copy of $D_{1,6}$
  such that the distinguished vectors described above correspond
  to 6 smooth rational curves $R_i$ and a nef class $H$
  with $H^2 = 4$.

  This is equivalent to the requirement that the $R_i$'s contain three
  even sets of rational curves, and $H+\sum R_i$ is divisible by two.
\end{defn}

In proposition \ref{prop:campedelli-are-todorov-enriques},
we proved essentially that the quotient of a smooth Campedelli surface by the
involution induced by a reflection of the form $s_\kappa$ is such an
Enriques surface: the six nodes are images of the isolated fixed points
of $s_\kappa$ and form three even sets.

In this case, the Enriques surface can be written as the quotient of
a complete intersection of three diagonal quadrics in $\mathbb P^5$
by a group $G \simeq (\Z/2)^3$: the linear system $H$ is generated by
$G$-invariants quadratic forms (the various $x_i^2$).


\begin{prop}
  On a $D_{1,6}$-Enriques surface $S$, the linear system $|H|$ is base
  point free and induces a map $S \to \mathbb P^2$ of degree 4.
\end{prop}

\begin{proof}
  According to \cite{Cossec-picard}, since $H$ is nef, $|H|$ has no
  fixed components (nef divisors with fixed components have
  self-intersection 2). Moreover $|H|$ has base points if and only if
  it has the form $2E+F$ (where $E$ and $F$ are half-elliptic pencils
  with $EF=1$) or $3E+R$ where $R$ is a smooth rational curve and
  $E$ is a half-pencil with $ER = 1$ \cite{Cossec-picard}*{2.12}.

  Following \cite{Morrison}, we note that
  $E \cdot H = 1$, but
  \[ E \cdot (H + \sum R_i) \]
  is an even integer, so $E$ should intersect one of the $R_i$'s,
  let $Q$ be this rational curve.

  \textit{The case $H = 3E+R$.} Then $0 = Q \cdot H = Q \cdot (3E +
  R)$, hence $Q \cdot R < 0$, hence $Q = R$, but then $0 = H \cdot Q =
  (3E+R) \cdot R = 1$, yielding a contradiction.

  \textit{The case $H = 2E+F$.} We note that $Q \cdot (2E+F)$ is
  zero, but since $Q \cdot E > 0$, $Q \cdot F < 0$, which is absurd.

  It is also known that $H^1(\Ocal_S(H)) = 0$, see
  \cite{Cossec-picard}.
\end{proof}

This proposition is a particular case of the following analogue of
\cite{Morrison}*{Lemma 5.1}

\begin{prop}
  Let $S$ be a smooth Enriques surface with $\nu$ disjoint rational
  $(-2)$-curves $R_i$, and $H$ a nef divisor with $H^2 >
  0$ such that $H \cdot R_i = 0$ and $H + \sum R_i$ is divisible by 2
  in $\Pic(S)$.

  Then $|H|$ contracts the curves $R_i$, and is base point free, except
  in the \emph{special} case: $H^2 = 2$ and there is one rational curve
  such that the others form an even set of nodes. In the special case,
  the special rational curve is a fixed component of $|H|$, and is
  not contracted by $|H|$.
\end{prop}

\begin{proof}
  According to \cite{Cossec-picard}, if $H$ has a fixed component or
  base points, it can be written $kE+F$, where $E$ and $F$
  are half-elliptic pencils with $EF = 1$, or $(k+1)E+R$ where $E$ is
  a half-pencil, $R$ a rational nodal curve and $ER = 1$. In both
  cases $H^2 = 2k$, and $HE = 1$. We need to show that we are then
  precisely in the \emph{special} case.

  Since $E$ has an even intersection number with $H+\sum R_i$
  (which is divisible by two), and $E \cdot H = 1$, for some $Q$ among the $R_i$'s,
  $E \cdot Q > 0$.

  Now $Q \cdot H = 0$. In the case $H = (k+1)E+R$, it follows that $QR < 0$,
  hence $Q = R$, and $HR = k+1+R^2 = 0$, which is a contradiction, except when
  $k=1$. If $H = kE+F$, then $HQ = 0$ gives $FQ < 0$, which is also a
  contradiction.

  Thus $k=1$, and it follows that $2E = H-Q$ is divisible by two, but since
  $H + \sum R_i$ is also divisible by two, $\sum_{R_i \neq Q} R_i$
  is an even set of nodal curves, hence we are in the special case.
  In other cases, $H$ is base point free and has no fixed component:
  since $H \cdot R_i = 0$, the curves $R_i$ are contracted by $|H|$.

  Conversely, under the hypotheses of the special case, let $Q$ be the distinguished
  rational curve. Then $H+Q$ and $H-Q$ are divisible by two and
  effective. Let $H = 2E + Q$. Then $h^0(H) \geq h^0(2E) \geq 2$,
  and since $H$ is nef and big, $h^0(H) = \chi(H) = 2$, hence $Q$ is a
  fixed component of $H$.
\end{proof}

\begin{prop}
  Let $S$ be a $D_{1,6}$-polarised Enriques surface, $H$ the
  distinguished positive class in $D_{1,6}$. Then $S \to |H|^\vee$ is
  a surjective morphism of degree 4, contracting the 6 rational curves
  and ramified over six lines.
\end{prop}


\begin{proof}
  As before, $H$ is base point free, has no fixed components, and
  $H \cdot R_i = 0$ so $|H|$ contracts all curves $R_i$.
  Consider the three elliptic pencils $2E_1 = H-R_1-R_2$, $2E_2 =
  H-R_3-R_4$, $2E_3 = H-R_5-R_6$ (by construction they are
  2-divisible). They have a natural interpretation as linear
  subsystems of $|H|$, so they are pulled back from pencils of lines
  in $|H|^\vee$.

  The corresponding half-pencils map to six lines in the
  plane. Each of them is part of the ramification locus, since pulling
  back one of these lines gives a multiplicity 2 divisor.

  Now note that $K_S = K_{\mathbb P^2} + R$ (the formula is both valid
  for $S$ and the nodal surface obtained by contracting the $R_i$'s)
  where $R$ is the ramification locus of $S \to
  \mathbb P^2$.  Then $R \equiv 3H + K_S$ (since $K_{\mathbb P^2} \equiv -3H$),
  and the definition of the pencils implies
  \[ R \equiv \sum_{i=1,2,3} E_i + (E_i+K_S) + \sum_{i=1}^6 R_i
  \text. \]
  Since $R$ is an effective divisor containing all half-pencils
  as well as the $R_i$'s, the
  linear equivalence above is an equality.
\end{proof}


In the proof above, we see that $E_1 \cdot R_1 = E_1 \cdot R_2 = 1$:
this indicates that the image point of $R_1$ and $R_2$ is actually
the base point of the elliptic pencil $|2E_1|$, which is the
intersection of the images of $E_1$ and $E'_1 = E_1 + K_S$ in the
plane.

\begin{prop}
  \label{prop:Enriques-galois-cover}
  The above morphism is a Galois cover with Galois group $(\Z/2)^2$.
\end{prop}

If $X$ is a Campedelli surface, with an involution $s_\kappa$, the
Galois cover $X \to \mathbb P^2$ given by the bicanonical map
(proposition \ref{prop:campedelli-galois-cover}) factors through the
quotient $S_\kappa$, and the Galois group of the cover $S_\kappa \to
\mathbb P^2$ is identified with the order 4 group $\Gamma / \pair{G, s_\kappa}$.

\begin{proof}[Proof of proposition \ref{prop:Enriques-galois-cover}]
  Let $E_i$ ($i=1,2,3$) be the three half-pencils as above, and
  consider the divisor $M_i = E_j + E_k$, where $\set{i,j,k} =
  \set{1,2,3}$

  Then $M_1$ is a genus two linear system, and by \cite{Cossec-models}*{6.1}
  we know that $|2M_1|$ defines a degree two morphism onto a degree 4
  del Pezzo surface. Let $s_1$ be the associated involution. Note that
  $|H|$ is a linear subsystem of $|2M_1|$ by the map $D \mapsto D + L_{23}$
  where $L_{23}$ is the unique element of $|H-R_3-R_4-R_5-R_6|$.
  Hence $X \to |H|^\vee$ is $s_1$-equivariant.

  Such involutions determine the divisor class $|2M_i|$ as the
  ramification divisor of $X \to X/s_i$.

  Hence $s_i$ ($i=1,2,3$) are distinct involutions, and since $X \to
  |H|^\vee$ has degree four, the group $\pair{1,s_1,s_2,s_3}$ is
  exactly the Galois group.
\end{proof}

\index{revêtement bidouble}
A Galois cover as above is called a \emph{bidouble} cover: it is a
special case of abelian covers which are studied in
\cite{Pardini-covers}. Such a cover has three involutions: for each of
them, we consider the image of the (divisorial part of the) fixed
locus in the base variety (here $\mathbb P^2$). Write $f$, $g$, $h$
for sections of line bundles defining them. By Galois theory, $fgh$
can be written as a square $s^2$, where $s$ is a section of some line
bundle.

Then the bidouble cover can be defined by equations of the form
\[ u^2 = f \qquad v^2 = g \qquad w^2 = h \qquad uvw = s \]
(changing signs of $u$, $v$, $w$ gives the other component $uvw=-s$).
In the case of $D_{1,6}$-Enriques surfaces, the divisors of
$f$, $g$, $h$ are $e_i + e'_i + e_j + e'_j$ (denoting by $e_i$ the
image of $E_i$ in the plane). This point of view highlights the fact
that the data of three pairs of lines in the planes
$(e_i,e'_i)$ $(i=1,2,3)$ determines uniquely a surface which is
generically the blow-down of six rational curves on an Enriques
surface.

\subsection{An involution of the moduli space}
\label{ssec:involution}

Consider a $D_{1,6}$-polarised Enriques surface $S$: we will now
choose more symmetric notations. The six rational curves are
denoted by $R_i^+$ and $R_i^-$, the six half-pencils $E_i^+$ and
$E_i^-$, such that $2E_i^+ = 2E_i^- = H - R_i^+ - R_i^-$.
The linear system $|H|$ determined a finite Galois cover
$S_r \to \mathbb P^2$, where $S_r$ is obtained from $S$ by contracting
the $(-2)$-curves $R_i^\pm$.

Assume that $S$ is given by three pairs of lines in general position,
whose vertices are denoted by $a_i$ (which are the images of the
double points of $S$).

Let $\mathcal P$ be the blowup of the three vertices: the pencil of
lines through $a_i$ corresponds to a linear system $|e_i|$ on $\mathcal
P$ such that $e_i^2 = 0$, and $e_i e_j = 1$ if $i \neq j$. The proper
transforms of the three pairs of lines define pairs of divisors
$e_i^\pm$ in each $|e_i|$.

Let $S'$ be the (well-defined) bidouble cover of $\mathcal
P$ ramified over the three pairs $(e_i^\pm)$. The linear system $e_i$
pulls back to the elliptic pencil $|2E_i|$.
Note that $\mathcal P$ has 6 $(-1)$-curves, given by the $a_i$ and
$\ell_i$ (proper transform of the line $a_j a_k$). The inverse image
of $a_i$ in $S'$ is made of two disjoint curves since $a_i$ is
disjoint from $e_j^\pm$ and $e_k^\pm$ ($\set{i,j,k} = \set{1,2,3}$):
it is actually the union of two $(-2)$-curves. It follows that $S$ and
$S'$ are isomorphic, and that the inverse image of $a_i$ is identified
with $R_i^+ + R_i^-$.

\index{involution!transformation de Cremona}
The curves $a_i$ and $\ell_i$ play equivalent roles (they can be exchanged by a
standard quadratic transformation). The intersection numbers on
$\mathcal P$ are $a_i \cdot \ell_j = 1$ if $i \neq j$, zero otherwise.
Using a similar proof, we show that $\ell_i$ pull back to two disjoint
$(-2)$-curves $\Lambda_i^\pm$: $\ell_i$ is also disjoint from $e_j^\pm$ and $e_k^\pm$.
If $h$ denotes the linear system of lines on $\mathbb P^2$ and its
pull-back to $\mathcal P$, $h = e_i + a_i$ and $h = \ell_i + a_j + a_k$.

The $D_{1,6}$-polarisation on $S$ can be described by $3H = \sum \Lambda_i^\pm
+ 2 \sum R_i^\pm$ and the six rational curves $R_i^\pm$ mapping to
$a_i$. Let $h' = e_i + \ell_i = 2h-a_1-a_2-a_3$.
\begin{prop}
  The pull-back $H'$ of $h'$ to $S$ is an effective divisor such
  that ${H'}^2 = 4$, and $\Lambda_i^\pm$ define six disjoint
  $(-2)$-curves such that $H' + \sum \Lambda_i^\pm$ is divisible by
  two. In other words $(H', \Lambda_i^\pm)$ is another
  $D_{1,6}$-polarisation on $S$.
\end{prop}

\begin{proof}
  The pullback of $h'$ is $H' = 2E_i + \Lambda_i^+ + \Lambda_i^-$,
  which is an effective divisor and $(H')^2 = 4$. The formula
  \[ H' + \sum \Lambda_i^\pm = (3H - H) - \sum R_i^\pm + \sum
  \Lambda_i^\pm = 2 \big( \sum \Lambda_i^\pm + \sum R_i^\pm \big) -
  (H + \sum R_i^\pm) \]
  shows that it is also 2-divisible.

  Moreover $H' - \Lambda_i^+ - \Lambda_i^-$ is the pull-back of
  $h'-\ell_i = e_i$, which is the elliptic pencil $2E_i$.
  This proves that the lattice generated by $H'$ and the
  $\Lambda_i^\pm$ is $D_{1,6}$.
\end{proof}

\begin{prop}
  The new $D_{1,6}$-polarisation defines another configuration of
  lines which can be obtained from the initial one by performing a
  Cremona quadratic transformation of the plane, with vertices $a_i$. \qed
\end{prop}


\begin{prop}
  The nonzero intersection numbers of the $R_i$'s and $\Lambda_i$'s are
  $R_i^{\pm} \cdot \Lambda_j^\pm = 1$ for $i \neq j$. \qed
\end{prop}

\begin{prop}
  A $D_{1,6}$-polarised Enriques surface $S$ defines an embedding
  of the root lattice $D_6 \subset H^2(S, \Z_-)$, using the classes
  of the six rational curves $R_i^\pm$ as a standard orthogonal frame.
\end{prop}

\begin{proof}
  The curves $R_i^\pm$ and $\Lambda_i^\pm$ define classes in $H^2(S,
  \Z_-)$ (up to a choice of sign), using their fundamental classes and
  the fact that $\Z_-$ restricts to a trivial local system on rational
  curves.  Their self-intersection is $-2$: this gives an embedding of
  $\Z^6(-2)$ in $H^2(S, \Z_-)$, since the $R_i^\pm$ do not intersect.

  The intersection numbers show that the classes $[R_1^+]$,
  $[\Lambda_3^+]$, $[R_2^+]$, $[\Lambda_1^+]$, $[R_3^\pm]$ (up to a
  choice of signs) form a $D_6$ lattice: the associated curves define
  a dual graph which is the Dynkin diagram of $D_6$ (this is graph is
  simply connected, so there actually exists a suitable choice of
  sign).

  By a small deformation, we can assume that the Hodge structure on
  $H^2(S,\Z_-)$ is generic (using the infinitesimal period map), and
  that the $D_6$ lattice is actually the whole $H^{1,1}(S, \Z_-)$.
  It is then easily checked that $[R_2^+] + [R_2^-] + [R_3^+] + [R_3^-]$
  is divisible by 2: it has even pairing with any of the roots mentioned above.
\end{proof}


Configurations of six lines may be represented by elements of
$\Gr(3,6)$: the coefficients of six equations of lines in $\mathbb
P^2$ may be written in the 6 columns of a $3 \times 6$ matrix, whose
lines define a three-dimensional subspace of $\bbC^6$ if the lines do
not pass through a common point. A given point in $\Gr(3,6)$ defines
the configuration only up to a projective transformation. We now
prove the following fact:
\begin{prop}
  \label{prop:cremona-involution}
  There is an explicit biregular involution $Q$ on the Grassmann
  variety $\Gr(3,6)$, such that for general $x \in \Gr(3,6)$,
  $x$ and $Q(x)$ represent conjugate configurations under the
  transformation described above: $Q(x)$ is equivalent to the
  Cremona transform of $x$ with vertices $x_1 \wedge x_2$,
  $x_3 \wedge x_4$, $x_5 \wedge x_6$.
\end{prop}


\begin{proof}
  Given a configuration of six lines in general position, there exists
  a choice of basis such that $A_1 = [1:0:0]$, etc. Let $M$ be a
  matrix of the configuration given by columns $M_1$, \dots, $M_6$,
  and $N$ be the matrix whose lines are $M_1 \wedge M_2$,
  $M_3 \wedge M_4$, $M_5 \wedge M_6$ (which are coordinates of the
  vertices of each pair of lines).

  Then $N$ is invertible and
  \[ NM = \begin{pmatrix}
    0 & 0 & m_{123} & m_{124} & m_{125} & m_{126} \\
    m_{341} & m_{342} & 0 & 0 & m_{345} & m_{346} \\
    m_{561} & m_{562} & m_{563} & m_{564} & 0 & 0 \\
  \end{pmatrix} \]
  where $m_{ijk} = \det(M_i,M_j,M_k)$.

  The quadratic transformation with centres $A_i$ can be described by
  inverting
  each of the non-zero coefficients of $NM$, or
  equivalently, exchanging nonzero coefficients in each column
  (up to a dilation on columns).

  Let $Q$ be the transformation of $\mathbb P^{19} \supset \Gr(3,6)$
  which exchanges the 6 pairs of coordinates as above
  (e.g. $x_{123}$ and $x_{563}$) and leaves the others untouched.
  Then $Q$ is a well-defined involution of $\Gr(3,6)$, and it acts as
  specified.
\end{proof}

\index{involution!association de Coble}
The description of the action of Coble's association show that it
differs from $Q$ by the transpositions of each pair of columns.
In particular if $t$ is an element of the torus $T$, then
$Q(t·x) = t' Q(x)$, where $t'$ is obtained from $t'^{-1}$
by exchanging the pairs of coefficients. It follows $Q$ also acts on
the ring of $T$-invariants of $\Gr(3,6)$, and descends to an
involution of $\Gr(3,6) \sslash T$.

\subsection{GIT stability for configurations of six lines}

The moduli space we are interested in is the space of configurations
of six lines, which can be described as a GIT quotient
$\bbC^{18} \sslash (GL_3 \times T)$ where $T$ is the diagonal torus in
$GL_6$ and $GL_3 \times T$ acts on $\bbC^3 \otimes \bbC^6
\simeq \bbC^{18}$. This space is also $(\mathbb P^2)^6 \sslash SL_3$
whose coordinate ring is determined by the $\mathfrak S_6$-equivariant
line bundle $\Ocal(1, \dots, 1)$, or
$\Gr(3,6) \sslash T$ where $\Gr(3,6)$ is the Grassmannian, with the linear
action of $T$ on $\bigwedge^3 \bbC^6$ (which is the target of the Plücker
embedding of the Grassmannian).

\index{stabilité!configurations de droites}
Semi-stability and instability can be given the following
description: we say a configuration of points in the plane has type
$ijk$ with respect to a flag $\bbC^3 = V_1 \supset V_2 \supset V_3$ of $\bbC^3$
iff the associated matrix has $i$ (resp. $j$, $k$) columns belonging
to the first (resp. second, third) element of the flag (note that
$i+j+k=6$). The type of a configuration can be read on the Schubert cell
it belongs to: configurations of type $ijk$ corresponding to a Schubert cell
$X_{i,i+j,i+j+k}$. Stability can be decided from the type, using
Seshadri's criterion. In this particular case, the conditions are very explicit:

\begin{prop}
  Given a flag of $\bbC^3$ and a 1-parameter subgroup of $SL_3$
  $t \mapsto (t^a, t^b, t^c)$ where $a>b>c$, a point $p$ of
  $(\mathbb P^2)^6$ is unstable (resp. non stable) w.r.t. $(g_t)$ iff
  $p$ is described by an element of the Schubert cell
  $\mathcal X_{i,i+j,i+j+k}$ where $ia + jb + kc < 0$ (resp. $\leq
  0$).
\end{prop}

\begin{proof}
  In the Segre embedding, basic coordinates of $p$ are given by the
  product of coordinates from each points: the given formula is the
  biggest weight which can be obtained in this way.
\end{proof}

It is enough to check stability against subgroup whose weights are
$(2,-1,-1)$ and $(1,1,-2)$: we need to check the sign of
$2i-(j+k) = 3i-6$ and $(i+j)-2k = 6-3k$.

\begin{prop}
  If a configuration has type $0xx$ in some basis, it is unstable.
\end{prop}

In the following we thus only consider configurations having type
$ijk$ where $i>0$. The shape of configurations of given type
is here described by interpreting them as configurations of lines,
which is more directly related to the shape of corresponding
surfaces.

\begin{prop}
  The unstable types are $1xx$ and $xx3$. The strictly semi-stable
  types are $231$ (4 concurrent lines), $312$ (2 identical lines),
  and $222$.
\end{prop}


The notions of stability, semi-stability and polystability coincide
for the GIT quotients $\Gr(3,6) \sslash T$ and $(\mathbb P^2)^6 \sslash SL_3$.

We obtain the following stratification of the Grassmannian:
\begin{itemize}
\item unstable 141 (codim. 3): five concurrent lines;
\item unstable 213 (codim. 4): three identical lines;
\item polystable 222 (codim. 6): three pairs of identical lines,
  the stabiliser has dimension 2, note there are three possible
  combinatorial types $(AA,BB,CC)$, $(AA,BC,BC)$ or $(AB,BC,CA)$
  (where $A$, $B$, $C$ denote lines from each pair);
\item semistable 222 (codim. 3): four concurrent lines, two of them
  being identical;
\item polystable 231=312 (codim. 4): 2 identical and 4 concurrent;
  there are two ways, up to permutation, in which they can be arranged
  (either two lines in the same pairs coincide, or two lines from
  different pairs coincide), the stabiliser has dimension 1;
\item semistable 231 (codim. 2): 4 concurrent lines;
\item semistable 312 (codim. 2): 2 identical lines;
\item type 321 (codim. 1): only one concurrent triple of lines;
\item type 411 (open): lines linearly in general position
  (intersection of the “big cells”).
\end{itemize}

The description of the stability locus can be summarised as follows:
\begin{prop}
  \label{prop:describe-stable-locus}
  A configuration of lines is stable if and only if it consists of six
  distinct lines and has at worst triple points.
\end{prop}

\subsection{The cohomology of a generic
  \texorpdfstring{D$_{1,6}$}{D(1,6)}-Enriques surfaces}
\label{ssec:tanscendental-for-d6}

\subsubsection{Cohomology lattices of Enriques surfaces}

Let $S$ be a generic smooth Enriques surface of type $D_{1,6}$.
Traditional presentations of the period map of Enriques surfaces
use the anti-invariant cohomology lattice of the double cover
$\pi: T \to S$, which is isomorphic to $H \oplus E_{10}(-2)$.
For some uses, it may be more convenient to study the unique
non-trivial local system of rank one on $S$, which we denote
by $\mathbb Z_-$. A lattice-theoretic version of this choice was
described by Allcock \cite{Allcock}.

\begin{prop}
  The cohomology groups of $\Z_-$ are $H^0 = 0$,
  $H^1 = \Z/2\Z$ and $H^2$ is torsion-free.
\end{prop}

\begin{proof}
  From the exact sequence of sheaves
  \[ 0 \to \Z_- \xrightarrow{\pi^{-1}} \pi_\ast \Z_T
  \xrightarrow{\pi_\ast} \Z_S \to 0 \]
  we deduce
  \[ 0=H^0(S, \Z_-) \to H^0(T, \Z) \to
  H^0(S, \Z) \to H^1(S, \Z_-) \to H^1(T, \Z)=0 \]
  which gives $H^0$ and $H^1$, since $H^0(T, \Z) \to H^0(S, \Z)$
  is integration on fibres.

  Since $H^1(S, \Z) = 0$, the torsion subgroup of $H^2(S, \Z_-)$ maps
  injectively to the torsion subgroup of $H^2(T, \Z)$, which is zero.
\end{proof}

\begin{prop}
  The torsion subgroups of $H^2(S, \Z_-)$ and
  $H^3(S, \Z_-)$ are isomorphic.
\end{prop}

\begin{proof}
  Since $\Z_-$ is Verdier self-dual, Poincaré-Verdier duality gives a
  (degenerate) spectral sequence
  \[ E_2^{p,q} = \Ext^p (H^{n-q}(S, \Z_-), \Z) \implies
  H^{p+q}(S, \Z_-) \]

  Hence there is a short exact sequence
  \[ 0 \to \Ext^1(H^2(S, \Z_-), \Z) \to
  H^3(S, \Z_-) \to \Hom(H^1(S, \Z_-), \Z) = 0 \]
  And similarly
  \[ 0 \to \Ext^1(H^3(S, \Z_-), \Z) \to
  H^2(S, \Z_-) \to \Hom(H^2(S, \Z_-), \Z) \to 0 \]
\end{proof}

\begin{prop}
  The lattice $H^2(S, \Z_-)$ is odd and unimodular, with
  signature $(2,10)$.
\end{prop}

\begin{proof}
  The parity of the lattice follows from Wu's formula: the Steenrod
  square $Sq^2$ coincides with the cup product with the second
  Stiefel-Whitney class of $S$. But the class $w_2(S)$ can be seen as
  the reduction mod~2 of $c_1(S) \in H^2(S, \Z)$, which is nonzero.
  Since the reduction map $H^2(S, \Z_-) \to H^2(S, \Z/2)$ is
  surjective ($H^3(\Z_-) = 0$),
  the intersection form is odd.

  Unimodularity follows from Poincaré duality. Additionally, $H^2(S,
  \R)$ and $H^2(S, \R_-)$ have the same signature, either by using
  index formula, or by identifying $H^2(S, \R_-)$ with the
  anti-invariant subspace of $H^2(T, \R)$, and resorting to standard
  calculations \cite{BHPV}.
\end{proof}

\begin{prop}
  Let $H^2(T, \Z)_-$ be the sublattice of vectors in $H^2(T, \Z)$
  which are anti-invariant under the involutive deck transformation of $T \to S$.
  The pull-back map from the lattice $H^2(S, \Z_-)$ to
  $H^2(T, \Z)_-(1/2)$ is an isometric embedding of index 2.
\end{prop}

\begin{proof}
  This is either deduced from our computations in section \ref{ssec:index-comp},
  or from the unimodularity of $H^2(S, \Z_-)$ and the
  fact that $H^2(T)_-(1/2)$ is isometric to $H(1/2)
  \oplus E_{10}(-1)$ \cite{BHPV} which has discriminant 1/4.
\end{proof}

There is actually a unique odd unimodular sublattice of
$H(1/2) \oplus E_{10}(-1)$, as shown by Allcock in \cite{Allcock}.

\begin{thm}[Allcock]
  There is a unique odd unimodular sublattice of $H^2(T)_-(1/2)$.
  The dual of $H^2(T)_-(1/2)$ is characterised as the lattice of even
  vectors in this unimodular lattice.
\end{thm}

\begin{proof}
  Let $L$ be the lattice $H(1/2) \oplus E_{10}(-1)$ and $L^\vee = H(2)
  \oplus E_{10}(-1)$ be the dual lattice. If $M$ is a unimodular
  sublattice of $L$, then $L^\vee \subset M^\vee = M \subset L$.

  Note that $L / L^\vee = H(1/2) / H(2)$, which is isomorphic to
  $(\Z/2)^2$. The quadratic form on $L/L^\vee$ is well defined mod 2,
  and takes values $0$, $0$, $1$, on the three nonzero vectors.
  So among the three integral lattices between $L$ and $L^\vee$,
  there is only one which is odd.
\end{proof}


\subsubsection{The orthogonal complement of the polarising lattice}

We are interested in the embedding $D_6 \into \Z^{2,10}$ which
represents the type $(1,1)$ part of the lattice $H^2(S, \Z_-)$ for a
generic $D_{1,6}$-Enriques surface $S$. Its orthogonal complement
(the “transcendental part”) will be denoted by $L$.
Note that $D_6 \subset H^2(S, \Z_-)$
is primitive, since its pull-back to $H^2(T, \Z)_-$ is also primitive.

First note that the dual of $D_6 \subset \Z^6$ (the embedding being
the standard embedding of $D_6$ as the even sublattice of $\Z^6$) is the lattice
generated by $\Z^6$ and $(1/2, \dots, 1/2)$. The discriminant group
of $D_6$ is then generated by a basis vector of $\Z^6$ and the
half-integer vector. The matrix of the associated bilinear form
$b: \Sym^2(D_6^\ast/D_6) \to \Q/\Z$ is
\[ \begin{pmatrix}
  0 & 1/2 \\
  1/2 & 1/2 \\
\end{pmatrix} \]
The quadratic form takes the basis vector to $1 \pmod 2$, the
half-integer vector to $-1/2 \pmod 2$ and their sum to $-1/2$.

The following proposition is classical:
\begin{prop}
  Let $M$ be a primitive sublattice of a unimodular lattice $L$.
  The inclusions $M \oplus M^\perp \subset L \subset M^\ast
  \oplus (M^\perp)^\ast$ define maps $L \to M^\ast / M$ and
  $L \to (M^\perp)^\ast / M^\perp$ which induce a bijective
  correspondence between the discriminant groups.


  This correspondence is an isometry for the associated bilinear forms
  $-b_M$ and $b_{M^\perp}$ with values in $\Q/\Z$. If all lattices
  involved are even, it is also an isometry for the quadratic forms
  $-q_M$ and $q_{M^\perp}$, with values in $\Q/2\Z$.

  Conversely, given such an isometry of bilinear forms, the pull back
  of the graph along the map $(M \oplus M^\perp)^\ast \to D_M \oplus
  D_{M^\perp}$ defines a unimodular lattice.
\end{prop}

\begin{prop}
  The orthogonal complement of $D_6(-1)$ in $\Z^{2,10}$ is a lattice
  $L$ of signature $(2,4)$ and discriminant 4. Its discriminant
  bilinear form has matrix
  $\big(\begin{smallmatrix} 0 & 1/2 \\ 1/2 & 1/2 \\ \end{smallmatrix}\big)$.
\end{prop}

\begin{lemma}
  The lattice $L$ is isomorphic to an index two sublattice of
  $\Z^{2,4}$.
\end{lemma}

\begin{proof}
  The values of the discriminant bilinear form show that there is a
  lattice between $L$ and $L^\ast$ which is integral.
  Since it must be unimodular, it is isomorphic to $\Z^{2,4}$.
\end{proof}

Note that any sublattice of index two in $\Z^{2,4}$ has the form
$\Z^{p,q} \oplus D_{r,s}$ where $D_{r,s}$ is the lattice of vectors in
$\Z^{r,s}$ whose coordinates have an even sum: such a sublattice is
necessarily obtained as the kernel of a group homomorphism $\Z^{2,4} \to \Z/2\Z$,
which is the set of vectors whose sum of specified coordinates is
even. The discriminant of $D_{r,s}$ is 4.

\begin{lemma}
  A sublattice $\Z^{p,q} \oplus D_{r,s}$ of $\Z^{2,4}$ has a 2-torsion
  discriminant group if and only if $r+s$ is even. The index 2
  sublattices of $\Z^{2,4}$ with 2-torsion discriminant group are
  \[ \begin{aligned}
    D_{2,4} &\simeq H \oplus H \oplus \Z^2(-2) \\
    D_{0,4} \oplus \Z^2 &\simeq H(2) \oplus \Z^{1,3} \\
    D_{2,2} \oplus \Z^{0,2} &\simeq H(2) \oplus H \oplus \Z^2(-1) \\
    D_{2,0} \oplus \Z^{0,4} &\simeq \Z^2(2) \oplus \Z^4(-1) \\
    D_{0,2} \oplus \Z^{2,2} &\simeq \Z^{2,2} \oplus \Z^2(-2) \\
    D_{1,3} \oplus \Z^{1,1} &\simeq H \oplus \Z^2(-2) \oplus \Z^{1,1} \\
    D_{1,1} \oplus \Z^{1,3} &\simeq H(2) \oplus \Z^{1,3} \\
  \end{aligned} \]
\end{lemma}

A basis for $H \subset D_{1,3}$ is given by $(1;1,0,0)$
and $(1;0,1,0)$. A basis for $\Z^{0,3}$ in $\Z \oplus D_4(-1)$
is given by $(1;1,1,0,0)$, $(1;0,1,1,0)$, $(1;0,0,1,1)$.
Note that
\[ D_{2,2} \oplus \Z^{0,2} \simeq D_{1,1} \oplus \Z^{1,3}
\simeq D_{0,4} \oplus \Z^{2,0} \]
and $D_{2,0} \oplus \Z^{0,4} \simeq  D_{1,3} \oplus \Z^{1,1}
\simeq  D_{0,2} \oplus \Z^{2,2}$.

\begin{lemma}
  The discriminant bilinear forms of $D_6$ and $D_{r,s}$ coincide
  if and only if $r-s \equiv 6 \pmod 4$. This rules out $(r,s) =
  (2,2)$ and $(r,s) = (1,1)$ and $(r,s) = (0,4)$.
\end{lemma}

\begin{lemma}
  The discriminant \emph{quadratic forms}, with values in $\Q/2\Z$, of
  $D_6$ and $D_{2,4}$ coincide. In particular, any unimodular lattice
  containing orthogonal primitive copies of $-D_6$ and $D_{2,4}$
  should be even.
\end{lemma}

\begin{proof}
  Any isomorphism between the discriminant groups which is an
  isometry between the bilinear forms is automatically an isometry for
  the quadratic forms whenever they coincide.
\end{proof}

There is only one case left.
\begin{thm}
  The lattice $L$ is isomorphic to $D_{2,0} \oplus \Z^{0,4} \simeq
  D_{1,3} \oplus \Z^{1,1} \simeq  D_{0,2} \oplus \Z^{2,2}$. \qed
\end{thm}

\begin{corol}
  All primitive embeddings $D_6(-1) \subset \Z^{2,10}$ are conjugate.

  Any primitive embedding $D_6(-1) \subset \Z^{2,10}$ can be factored
  as $D_6(-1) \to \Z^{2,6} \to \Z^{2,6} \oplus \Z^{0,4}$ or
  $D_6(-1) \to E_8(-1) \to \Z^{2,2} \oplus E_8(-1)$.
\end{corol}

\begin{proof}
  For any primitive embedding $D_6(-1) \subset \Z^{2,10}$, $D_\perp$
  contains a copy of $\Z^{0,4}$, hence the image of $D_6(-1)$ is
  contained in the orthogonal complement of $\Z^{0,4}$, which is
  isomorphic to $\Z^{2,6}$. Similarly, $D_6(-1)$ and $D_{0,2}$
  should span a copy of $E_8(-1)$.

  The fact that all primitive embeddings are conjugate follows from
  the construction of the unimodular lattice: suppose we are given two
  sublattices of the form $D_6(-1) \oplus D_\perp$ inside $\Z^{2,10}$,
  and an isomorphism between the copies of $D_6(-1)$. Then there is a
  choice of isomorphism between the copies of $D_\perp$ which is
  compatible with the correspondences between discriminant groups
  induced by $\Z^{2,10}$: it thus extends to an automorphism of
  $\Z^{2,10}$.
\end{proof}

A particular embedding of $D_6(-1)$ in $\Z^{2,10}$ is given by
the embedding of Dynkin diagrams between $D_6$ and $E_8$, and
the fact that $\Z^{2,10} \simeq \Z^{2,2} \oplus E_8(-1)$.

The Torelli property of Enriques surfaces can be propagated to a
smaller lattice provided that it embeds uniquely in the standard
Enriques lattice and any automorphism of $L$ extends to $\Z^{2,10}$.


\begin{corol}
  Two $D_{1,6}$-Enriques surfaces are isomorphic if and only if their
  periods with values in $L$ are equivalent.
\end{corol}

\begin{prop}
  \label{prop:gen-transcendental-lattice}
  The generic transcendental lattice of a $D_{1,6}$-Enriques surface,
  in the Hodge structure $H^2(S, \Z_-)$, is isomorphic to the lattice
  $\Z^2(2) \oplus \Z^4(-1)$.
\end{prop}



\section{The period map}

In this section we study the period map of $D_{1,6}$-polarised
Enriques surfaces, which are parametrised by the space of line configurations in
the plane: we have seen in section \ref{ssec:d16-enriques-geometry}
how three pairs of lines gave rise to an Enriques surface as a
bidouble cover of the plane.

The period map behaves well on the stable locus of the parameter space
(which can be identified with $\bbC^{3 \times 6}$): since extra
singularities appear when configurations of lines acquire triple
points, the discriminant locus consists of the vanishing locus of
Plücker coordinates (which map $\bbC^{3 \times 6}$ to the Grassmann
variety $\Gr(3,6) \subset \mathbb P^{19}$. This is divisor with normal
crossings (proposition \ref{prop:plucker-normal-crossings}), and the
finite cover obtained by adding square roots to these coordinates
produces a local uniformisation for the period map, with values in the
space $\mathcal X_L = \mathcal D(L) / O(L)$.

The goal is to actually prove that denoting by $\Mgit$ the quotient
$\Gr(3,6) \sslash T$ (where the stable locus parametrises
configurations having at worst triple points), the
quotient $\Mgit[s] / (W_3 \times \langle Q \rangle)$
(where $W_3$ is the wreath product $\Z/2 \wr \mathfrak S_3$
and $Q$ is the Cremona transformation of proposition
\ref{prop:cremona-involution}) is actually isomorphic to $\mathcal
X_L$ via the period map. Theorem \ref{thm:period-is-local-iso} states
that a suitably chosen point has only one preimage, and that the
period map is a local isomorphism around it (using section
\ref{ssec:period-double-pt}).

A criterion of Looijenga and Swierstra, described in section
\ref{sssec:boundary-to-boundary}, can be used to prove that the
complement of $\Mgit[s]$ in $\Mgit$ is mapped to the complement
of $\mathcal X_L$ in $\mathcal X_L^\mathrm{BB}$, i.e. that the
map $\Mgit[s] \to \mathcal X_L$ is proper, and the previous remark
show that it has degree one. It suffices to note that this map has
finite fibres (since fibres correspond to a set of [almost]
polarisations on some K3 surface), to obtained the required isomorphism.

The extension to the compactified moduli spaces follows from a
Hartogs-type argument. However, we give an explicit calculation in
section \ref{ssec:boundary-period-chi}, which details the geometry of
the singularities in the neighbourhood of one of the boundary strata.

\subsection{The period mapping and its extension to the stable locus}

Let $S$ be a $D_{1,6}$-polarised Enriques surface. Then $H^2(S, \Z_-)$
contains a distinguished $D_6$ sublattice as previously explained,
and the line of (twisted) 2-forms $\bbC \omega$ lies in the complex
vector space spanned by the orthogonal complement of the distinguished
sublattice. If $\phi$ is a choice of isometry between $D_6^\perp
\subset H^2(S, \Z_-)$ and $L$, $(S,\phi)$ determines a point in
$\mathbb P(L \otimes \bbC)$ which is the image of the line $\bbC
\omega$. Let $q$ be the quadratic form on $L$: for example, we can
choose integral coordinates $(x_1, x_2; y_1, y_2, y_3, y_4)$ such that
\[ q = 2 x_1^2 + 2x_2^2 - y_1^2 - y_2^2 - y_3^2 - y_4^2 \text. \]
The associated symmetric bilinear form is denoted by $q(a,b)$.

Since $H^2(S, \Z_-)$ is a polarised Hodge structure, a representative
$\omega$ must satisfy $q(\omega) = 0$ and $q(\omega,\bar \omega) > 0$.
The \emph{period point} of $(S,\phi)$ is the element $[\omega]$ of
the Hermitian symmetric domain associated to the lattice $L$ computed above:
\[ \mathcal D_L = \set{[\omega] \in \mathbb P(L \otimes \bbC)
  \text{ such that } \pair{\omega, \bar \omega} > 0 \text{ and }
  \pair{\omega, \omega} = 0} \]

\index{domaine de périodes}
The period domain $\mathcal D_L$ contains two connected components:
the equations of $\mathcal D_L$ induce the constraints
\[ \begin{aligned}
  2x_1^2 + 2x_2^2 &= y_1^2 + y_2^2 + y_3^2 + y_4^2 \\
  2|x_1|^2 + 2|x_2|^2 &> |y_1|^2 + |y_2|^2 + |y_3|^2 + |y_4|^2 \\
\end{aligned} \]
which imply, for example, that $|x_1^2 + x_2^2| < |x_1|^2 + |x_2|^2$.
It follows that a point of $\mathcal D_L$ with coordinates $(x_i; y_j)$
cannot be such that $x_1/x_2$ is real. The connected components of
$\mathcal D_L$ are distinguished by the sign of $\Im(x_2/x_1)$.

The correspondence between Enriques surfaces with a $D_{1,6}$
polarisation and their period points is holomorphic.
\begin{thm}[Griffiths \cite{Voisin}*{chapitre 10}]
  Let $\mathcal S \to B$ be a holomorphic family of $D_{1,6}$-polarised
  Enriques surfaces, parametrised by a base $B$, equipped with a
  continuous family of isomorphisms $\phi(b)$ between the orthogonal
  complement of the $D_6$ sublattice of $H^2(S_b, \Z_-)$ and~$L$.

  Then the map $B \to \mathcal D_L$ mapping a point $b$ to the period
  point of $(S_b, \phi(b))$ is holomorphic.
\end{thm}

Let $\Gamma = O(L)$ and $\mathcal X_L = \mathcal D_L / \Gamma$: note
that elements of $\Gamma$ can exchange the connected components of
$\mathcal D_L$. Our goal is to define a period map from the stable
locus $\Mgit[s]$ to $\mathcal X_L$, and then to extend it to a
morphism from the full GIT quotient $\Mgit$ to the Baily-Borel
compactification $\mathcal X_L^\mathrm{BB}$ of $\mathcal X_L$.

Let $\mathcal U^\mathrm{sm}$ and $\mathcal U^{\mathrm s}$
be the open subsets in $\bbC^{18}$ parametrising configurations of lines
without triple points (resp. with at worst triple points).
\begin{prop}
  There is a well-defined holomorphic map $\mathcal U^\mathrm{sm} \to
  \mathcal X_L$, which is locally liftable to $\mathcal D_L$.
  The automorphic line bundle over $\mathcal X_L$ lifts to the
  linearised line bundle $\Ocal(3)$ on $\mathbb P(\mathcal U^\mathrm{sm})$.
\end{prop}

\begin{proof}
  Let $\mathcal T \to \mathcal U^\mathrm{sm}$ be the family of K3
  surfaces with double points, defined in $\mathbb P^5$ by the
  equations $u^2 = f_t(x,y,z)$, $v^2 = g_t(x,y,z)$ and $w^2 =
  h_t(x,y,z)$, where $f_t$, $g_t$, $h_t$ are products of two linear
  forms, depending on the parameter $t \in \mathcal U^\mathrm{sm}$.
  There is an associated family $\widetilde{\mathcal T}$ of
  \emph{smooth} K3 surfaces over $\mathcal U^\mathrm{sm}$: it is
  obtained as the corresponding cover of $\mathcal P$ (a family of
  rational surfaces such that $\mathcal P_t$ is the plane blown up at
  the vertices of $f_t$, $g_t$, and $h_t$).

  Let $dF$, $dG$, $dH$ be the differentials of these equations,
  which can be naturally interpreted as elements of
  $H^0(T_t, N_{T_t/\mathbb P^5}^\ast(2))$, and $\tau$ be a fixed
  section of $(\det T\mathbb P^5)(-6)$ on $\mathbb P^5$
  (which is unique up to a scalar).

  Then the pairing of $\tau$ with $dF \wedge dG \wedge dH$ is a
  well-defined, nowhere degenerate bivector field $w$
  (section of the dual of $\Omega^2_{T_t}$) over the regular part
  of $T_t$, which defines a symplectic form $\omega = w^{-1}$
  on $\widetilde T_t$.

  Since $\widetilde{\mathcal T}$ is a locally trivial fibration,
  there is on the universal cover of $\mathcal U^\mathrm{sm}$
  a uniform choice of basis of twisted homology classes $\gamma_i$
  ($i = 1...6$), giving a period map from the universal cover of
  $\mathcal U^\mathrm{sm}$ to $\mathcal D_L$. This gives local lifts
  for the quotient map:
  \[ \mathcal U^\mathrm{sm} = \frac{\widetilde{\mathcal
      U^\mathrm{sm}}}{\pi_1(\mathcal U^\mathrm{sm})} \to
  \frac{\mathcal D_L}{\Gamma} = \mathcal X_L \]

  Let $T_{ks}$ be the K3 surface obtained from a surface $T_s$ by
  multiplying the matrix $s$ of linear forms by $k$. An isomorphism
  $m_k: T_s \to T_{ks}$ is obtained by multiplying
  $u$, $v$, and $w$, by $k$. The differential forms transform as
  $(m_k)^\star (dF dG dH)_{ks} = k^6 (dF dG dH)_s$, and $\tau$ pulls back as
  $(m_k^{-1})_\star \tau = k^{-3} \tau$. It results
  that $w_{ks} = k^3 w_s$ under the identification $m_k: T_s \simeq T_{ks}$
  and that the periods of $T_{ks}$ are $k^{-3}$ times the periods
  of $T_s$. This proves that the automorphic line bundle pulls back as
  $\Ocal(3)$ on $\mathbb P(\mathcal U^\mathrm{sm})$.
\end{proof}

The period map is of course equivariant under action of $GL_3$ and
$T$. Note that the 2-form $\omega$ defined above can be written
with the simpler formula:
\[ \omega = \frac{z dx \wedge dy + x dy \wedge dz + y dz \wedge dx}
{u v w} \]
since $\omega dF dG dH$ is a multiple of the standard 5-form on
$\mathbb P^5$ with values in $\Ocal(6)$ (replace $dF$ by $u du$, etc.).

\begin{prop}
  The period map from the universal cover of $\mathcal U^{\mathrm{sm}}$
  to $\mathcal D_L$ is submersive. Consequently, the map $\Mgit[sm]
  \to \mathcal X_L$ can be locally lifted to an étale map with target
  $\mathcal D_L$.
\end{prop}

\begin{proof}
  Since the two spaces have the same dimension, it is enough to prove
  that it is a submersion. But if $\mathcal V^{\mathrm{sm}}$ is the
  subspace of $\bbC^{4 \times 7}$ which parametrises smooth Campedelli
  surfaces, the forgetful map $\mathcal V^{\mathrm{sm}} \to \mathcal
  U^{\mathrm{sm}}$ is itself submersive.

  Since the local period map $\mathcal V^{\mathrm{sm}} \to \mathcal
  D_L$ is submersive (section \ref{sec:campedelli-jacobian}),
  the result follows.
\end{proof}

We are going to use the following theorem of Borel
\begin{thm}[Borel \cite{Borel-extension}*{Thm. A}]
  Let $U$ be a polydisc in $\bbC^n$, and $U^\ast$ the complement
  of a standard normal crossing divisor in $U$. Let $\mathcal D$ be a
  bounded symmetric domain, $\Gamma$ an arithmetic subgroup of the
  associated group $\mathcal G$, and let $X = \mathcal D / \Gamma$.

  Suppose we are given a holomorphic map $f: U^\ast \to X$. Then if $f$
  is locally liftable to $\mathcal D$, then $f$ extends to a
  holomorphic map $U \to X^\mathrm{BB}$ where $X^\mathrm{BB}$ is the
  Baily-Borel compactification of $X$.
\end{thm}

The extension to the whole stable locus can be carried out using
this theorem and the following property:
\begin{prop}
  \label{prop:plucker-normal-crossings}
  The complement of $\Mgit[sm]$ in $\Mgit[s]$ is made
  of coordinate hyperplane sections, with normal crossings.
\end{prop}

The geometry of the variety $\Mgit$ is very well known and studied
thoroughly, for example in \cite{MSY} in a similar context. We
nevertheless provide a proof for the needed properties. Note the
statement here concerns configurations of six ordered lines.

\begin{proof}
  First observe that a stable configuration of six lines, having at
  worst standard triple points, can have at most four such points
  (proposition \ref{prop:describe-stable-locus}).
  Indeed, a given line cannot go through three triple points, hence
  five triple points would involve at least 8 lines.
  If there are four triple points, then \emph{any} line goes through two of them.
  The graph having triple points as vertices and lines as edges is
  made of four vertices of valence three, which is only possible if
  it is a complete graph. This means that the combinatorics of the
  corresponding arrangement of lines (up to permutation in $\mathfrak
  S_6$) are fully determined (they form a complete quadrangle).

  Then note that a stable configuration of lines always contains four
  lines which constitute a projective frame. If this were not the
  case, the 15 sets of four lines would contribute 5 triplets of
  concurrent lines (because each triplet is part of at most 3 sets of
  four lines), which is impossible by the previous remark.

  Again up to a permutation, we can assume the stable configuration is
  described by a matrix of the form
  \[ \begin{pmatrix}
    1 & 0 & 0 & 1 & a & b \\
    0 & 1 & 0 & 1 & c & d \\
    0 & 0 & 1 & 1 & e & f \\
  \end{pmatrix} \]

  $\bullet$ For one triple point, we can assume that the vanishing minor is
  $m_{125} = e$, since a projective frame can be made out of every
  4-tuple of non concurring lines. Neither $a$ nor $b$ can vanish,
  hence can be set to $1$. The local structure of $\Mgit$ around this
  point is then described by the “étale slice”, and the locus of
  triple points is a smooth divisor $e=0$.

  In the case of two triple points, choosing two lines from each
  triple makes a projective frame: the vanishing minors are
  \begin{itemize}
  \item either $m_{125} = m_{346} = e = b-d = 0$; then $a$ and $b$
    cannot vanish again, hence we set them to 1, the locus of triple
    points is the union of $e=0$, $d=1$ which meet transversely;
  \item either $m_{125} = m_{345} = e = a-c = 0$; then $a$ and $f$
    cannot both vanish, the locus of triple points is a union of
    $e=0$, $c=1$.
  \end{itemize}

  $\bullet$ In the case of three triple points, two triples have a common line,
  say $125$ and $345$, up to permutation the third can be chosen to be $136$.
  The corresponding equations are $e=a-c=d=0$. Again $a$ and $b$
  cannot vanish, hence can be chosen to be $1$, and in the étale
  slice, the locus of triple points is defined by $ed(c-1) = 0$,
  around $e=d=c-1=0$.

  $\bullet$ In the case of four triple points, as before we can assume that two
  triples are $125$ and $345$, then the others can be chosen to be
  $136$ and $246$, corresponding to the vanishing of $e$, $d$, $b-f$
  and $a-c$. As before $a$, $b$ can be chosen to be 1, and the
  locus of triple points is then the union of $e=0$,
  $d=0$, $c=1$, $f=1$ which is again a normal crossing divisor.
\end{proof}

Triple points in the configuration of lines correspond to the
appearance of rational double points in the K3 surfaces $T_t$: each
triple point creates a local monodromy or order two. These local
monodromies are eliminated by the following process: let
$\widetilde{\mathcal U^s} = \mathrm{sq}^{-1}(\mathcal U^s)$ be the
finite cover obtained as the inverse image under the morphism
$\mathrm{sq}: \mathbb P^{19} \to \mathbb P^{19}$ which takes
each Plücker coordinate to its square. By the previous proposition,
$\widetilde{\mathcal U^s}$ is the smooth locally closed subscheme of
$\mathbb P^{19}$, and the action of the diagonal torus $T$ admits a
lift to $\widetilde{\mathcal U^s}$ (and actually also to the inverse
image of the Grassmannian, which is a singular variety).

\begin{prop}
  Let $\widetilde{\Mgit[s]}$ be the GIT quotient $\widetilde{\mathcal
    U^s} \sslash T$, which is a Galois cover of $\Mgit[s]$. Then the period map
  $\widetilde{\Mgit[s]} \to \mathcal X_L$ has local lifts to $\mathcal
  D_L$.
\end{prop}

\begin{proof}
  Using the previous remark and the calculation carried in section
  \ref{ssec:period-double-pt}, the finite cover $\widetilde{\Mgit[s]}
  \to \Mgit[s]$ trivialises the monodromy: the period map extends
  regularly and has local lifts to $\mathcal D_L$.
\end{proof}

\subsection{Local structure around the one-dimensional $\chi$ stratum}
\label{ssec:boundary-period-chi}

Consider a polystable configuration in the one-dimensional stratum of
$\Mgit$. It is made of four concurrent lines (whose cross-ratio is
denoted by $t$), and two identical lines. We focus on the $\chi$ case,
where the identical lines belong to the same pair.
Up to some permutation, we can assume that the configuration is given by the
matrix
\[ \begin{pmatrix}
  1 & 0 & 1 & 1 & 0 & 0 \\
  0 & 1 & 1 & t & 0 & 0 \\
  0 & 0 & 0 & 0 & 1 & 1 \\
\end{pmatrix} \]
with $t \neq 0$. Its stabiliser has dimension 1, it acts
by multiplying lines by $(\lambda^{-1}, \lambda^{-1}, \lambda^2)$
and columns by
$(\lambda, \lambda, \lambda, \lambda, \lambda^{-2}, \lambda^{-2})$.

An “étale slice” of $\Mgit$ around this configuration given by
the five parameter family:
\[ N(t,a,b,c,d) = \begin{pmatrix}
  1 & 0 & 1 & 1 & 0 & c \\
  0 & 1 & 1 & t & 0 & d \\
  0 & 0 & a & b & 1 & 1 \\
\end{pmatrix} \]
where $\mathbb G_m$ acts by $\lambda \cdot (t,a,b,c,d) \mapsto
(t, \lambda^3 a, \lambda^3 b, \lambda^{-3} c, \lambda^{-3} d)$
on the parameters (this action is equivalent to the action
on lines and columns above). The étale slice itself is
the quotient $\bbC \pair{t,a,b,c,d} \sslash \mathbb G_m$
whose coordinate ring is
$\bbC[t,ac,ad,bc,bd]$ (a cylinder over the quadratic cone with equation $(ac)(bd) =
(ad)(bc)$).

Also recall that for a generic configuration $\gamma$, the quadratic
transformation associated to the vertices of the three pairs of lines
defines a new configuration $Q(\gamma)$ with isomorphic associated
Enriques surface. Then $Q$ defines a biregular involution of $\Mgit$,
commuting with the action of $(\Z/2 \wr \mathfrak S_3)$.

\begin{prop}
  $Q$ lifts to a biregular involution of the étale slice mentioned
  above, given by the composite of $(a,b) \leftrightarrow (c,d)$
  and reversing the order of the four first columns. Therefore
  the period map is invariant under exchange of $(a,b)$ and $(c,d)$.
\end{prop}

\begin{proof}
  Recall that $Q$ can be defined as follows on matrices:
  \[ \begin{pmatrix}
    0 & 0 & a_2 & b_2 & c_3 & d_3 \\
    a_1 & b_1 & c_2 & d_2 & 0 & 0 \\
    c_1 & d_1 & 0 & 0 & a_3 & b_3 \\
  \end{pmatrix} \mapsto
  \begin{pmatrix}
    0 & 0 & c_2 & d_2 & a_3 & b_3 \\
    c_1 & d_1 & a_2 & b_2 & 0 & 0 \\
    a_1 & b_1 & 0 & 0 & c_3 & d_3 \\
  \end{pmatrix} \]
  Given a configuration with matrix:
  \[ N = \begin{pmatrix}
    1 & 0 & 1 & 1 & 0 & c \\
    0 & 1 & 1 & t & 0 & d \\
    0 & 0 & a & b & 1 & 1 \\
  \end{pmatrix} \]
  its three vertices are given by the columns of the matrix
  \[ M = \begin{pmatrix}
    b-at & -d & 0 \\
    a-b & c & 0 \\
    0 & 0 & 1 \\
  \end{pmatrix} \]
  and the associated matrix in this new basis is
  \[ M^T N = \begin{pmatrix}
    b-at & a-b & 0 & 0 & t-1 & bc+ad-bd-act+t-1 \\
    -d & c & c-d & ct-d & 0 & 0 \\
    0 & 0 & a & b & 1 & 1 \\
  \end{pmatrix} \]
  (write $f(a,b,c,d) = bc+ad-bd-act+t-1 = f(c,d,a,b)$).
  The transformed configuration $Q(M^T N)$ has matrix
  \[ \begin{aligned} &\begin{pmatrix}
    -d & c & 0 & 0 & 1 & 1 \\
    b-at & a-b & a & b & 0 & 0 \\
    0 & 0 & c-d & ct-d & t -1 & f(a,b,c,d) \\
  \end{pmatrix} \\
  \equiv &\begin{pmatrix}
    ct-d & c-d & 0 & 0 & t -1 & f(c,d,a,b) \\
    b & a & a-b & b-at & 0 & 0 \\
    0 & 0 & c & -d & 1 & 1 \\
  \end{pmatrix} \end{aligned} \]
  (the second matrix is obtained by reversing the first four columns
  as well as the lines). Now applying the substitution
  \[ (a,b,c,d) \mapsto (c,d,a,b) \]
  gives back (up to changes of signs in columns) the matrix $M^T N$.
\end{proof}

Note that the lift to the étale slice is equivariant under the action
of $\bbC^\times$. It is useful to define a similar étale slice for the
ramified cover $\widetilde{\Mgit}$: for this we need to formally add
square roots to all minors of $N(t,a,b,c,d)$.

If $a$, $b$, $c$, $d$ are in a small neighbourhood of zero, and $t$
lies in a small disc $U$ avoiding $0$ and $1$, then the only minors of
the configuration matrix which possibly vanish involve two columns
from the first four and one of the last two. More explicitly: three
lines become concurrent when either $a$, $b$, $a-b$, $at-b$, $c$, $d$,
$c-d$, $ct-d$ vanishes.

Let $\mathcal E_i = (E_t)$ be two identical families of cones of elliptic curves
with equation $u_i^2 = x_i^2 - y_i^2$, $v_i^2 = t x_i^2 - y_i^2$. Then
$\mathcal E_1 \times_U \mathcal E_2$ maps to a ramified cover of
$\bbC^5$ (with coordinates $t \in U$, $a$, $b$, $c$, $d$):
\[ a = x_1^2 \qquad b = y_1^2 \qquad c = x_2^2 \qquad d = y_2^2 \]
Each $\mathcal E_i$ is equipped with an action of $\mathbb G_m$,
$(u_i,v_i,x_i,y_i) \mapsto (\lambda u_i,\lambda v_i,\lambda x_i,
\lambda y_i)$.

\begin{prop}
  The action of $\mathbb G_m$ on the affine space lifts to the action
  of $\mathbb G_m$, $\mathcal E_1 \times_U \mathcal E_2$ with weights
  $(6,-6)$, and $\mathcal E_1 \times_U \mathcal E_2 \sslash \mathbb G_m$
  is a local model for $\widetilde{\Mgit}$.
\end{prop}

\subsubsection{The boundary period map}
\label{sssec:boundary-to-boundary}

To study the behaviour the period map, we will not work on $\Mgit$
itself, but on the étale slice of $\widetilde{\Mgit}$ defined earlier.
We denote by $Z_i$ the central curve of $\mathcal E_i$ (defined
by $u_i=v_i=x_i=y_i=0$, and by $\widetilde{\mathcal E_i}$
the blow-up of $\mathcal E_i$ along $Z_i$.

We already know that a period map is well-defined and locally liftable
over $\widetilde{\Mgit[s]}$:
\begin{prop}
  The period map from $\mathcal E_1 \times_U \mathcal E_2$
  to $\mathcal X_L$ is well-defined and locally liftable to $\mathcal D_L$
  outside of $Z_1 \times \mathcal E_2 \cup \mathcal E_1 \times Z_2$.
  \qed
\end{prop}

This observation is compatible with a classical description: if
$\Gamma$ is a neat arithmetic subgroup of $O(L)$, $\mathcal
D_L/\Gamma$ locally looks like a family over a disc $U$ of cones over
a product of elliptic curves, near a point of a 1-dim. boundary
component.

\begin{prop}
  The period map with source domain $\widetilde{\mathcal E_1} \times_U
  \widetilde{\mathcal E_2}$ is defined on the complement of a normal
  crossing divisor, which is the union of $\widetilde Z_1 \times
  \widetilde{\mathcal E_2}$ and $\widetilde{\mathcal E_2} \times
  \widetilde Z_2$ (here $\widetilde Z_i$ is the exceptional divisor in
  $\widetilde{\mathcal E_i}$).
\end{prop}

By Borel's extension theorem, the period map has a unique extension to
this desingularisation, which takes values in the Baily-Borel
compactification $\mathcal X_L^\mathrm{BB}$.

We now want to prove that for a fixed $t$, the exceptional divisor is
contracted to a point in $\mathcal X_L^\mathrm{BB}$, which would prove
that in some small neighbourhood of the semistable point, the period
map $\Mgit \to \mathcal X_L^\mathrm{BB}$ is well-defined. Since $Q$ lifts
to the ramified cover, and exchanges the two factors, it is enough to
prove it for one of the components: we choose to study the case where
$c=d=0$ (double line).

\begin{prop}
  Let $\Delta \to \widetilde{\mathcal E_1} \times_U
  \widetilde{\mathcal E_2}$ be a one-parameter family such that
  $c=d=0$ at the central point. Then for generic values of $a$ and
  $b$, the limit of the period map only depends on $t$.
\end{prop}

The proof of this proposition is delayed to the next section:
for the other boundary strata, we rely on a weaker result.
We are going to prove that the rational map given by the period map
$\Mgit \to \mathcal X_L^\mathrm{BB}$ maps the special semistable strata
(0-dimensional strata) to the boundary as well.

Consider a one-parameter degeneration with special fibre
a configuration of three double lines: the three ramification
divisors degenerate to $x^2$, $y^2$, $z^2$ (for given
projective plane coordinates $[x:y:z]$). We are using the same
criterion as Looijenga in \cite{Looijenga-4fold} and
\cite{LooijengaSwierstra}*{sec. 3}.

\begin{prop}[see \cites{LooijengaSwierstra,Looijenga-4fold}]
  Let $T_t$ be a one-parameter family of K3 surfaces
  whose special fibre is singular, and $\omega$ be a choice of
  twisted 2-forms on $(T_t)$ such that $\omega_{t=0}$ is a nonzero
  cohomology class in $H^2((T_0)_\mathrm{reg}, \bbC)$, but
  the $L^2$-norm of $\omega$ goes to infinity as $t \to 0$.

  Then the limit of the period map of the family as $t \to 0$
  belongs to $\mathbb P(\ker: H^2(T_t, \bbC) \to
  H^2((T_0)_\mathrm{reg}, \bbC))$. Moreover, the limit of the period
  map belongs to the boundary of the period domain.
\end{prop}

On the special fibre, the chosen 2-form can be written as
\[ \omega = \frac{x dy dz + y dz dx + z dx dy}{xyz} \]
which has non zero periods on the cycle
\[ \gamma: (\theta_1, \theta_2) \mapsto
   [\exp(i\theta_1) : \exp(i\theta_2): 1] \]
since for $z=1$, $\omega=dx/x \wedge dy/y$: which has period $(2i\pi)^2$.

Moreover $\omega \wedge \bar \omega = dx d\bar x dy d\bar y / |xy|^2$ is
not integrable over $\bbC^2$: it follows that for any smoothing of
this central fibre, the monodromy cannot have finite order. Now using
Borel's extension theorem, we know that some birational model of
$\Mgit$, which is an isomorphism outside the non-stable locus,
can be mapped regularly to $\mathcal X_L^\mathrm{BB}$. The discussion above
implies that the fibre above the special semistable point is mapped to
the boundary. The discussion is similar for the remaining boundary
strata, which correspond to similar line configurations, with various
labellings.

\subsection{The most special configuration and its isotropy group}
\label{sec:most-special-config}

We are interested in a very special point of $\mathcal X_L$.
\begin{prop}
  There exists a unique period point $[\omega_0]$ such that
  $\omega_0^\perp$ in $L$ is isomorphic to $\Z^4(-1)$.
\end{prop}
Since there are $(-1)$-classes in $\omega_0^\perp$, this $\omega_0$ is
of course not the period point of an Enriques surface (but rather the
one of a rational degeneration of Enriques surfaces).

\begin{proof}
  Such period points are in bijection with isometry classes of
  decompositions $L \simeq \Z^2(2) \oplus \Z^4(-1)$. Indeed,
  any positive plane $T$ such that $T^\perp \simeq \Z^4(-1)$
  gives such a decomposition. It is clear that any two such
  decompositions are conjugate under some isometry.
\end{proof}

We want to prove the following theorem: for convenience, $W_3$ denotes
the wreath product $(\Z/2 \wr \mathfrak S_3)$.
\begin{thm}
  \label{thm:period-is-local-iso}
  There is only one point which is mapped to $[\omega_0]$ under the
  period map $\mathcal P: \Mgit / (W_3 \times \pair{Q}) \to \mathcal
  X_L^\mathrm{BB}$, and $\mathcal P$ is a local isomorphism around these points.
\end{thm}
In order to prove the theorem, we will describe how to uniformise the
period map for each of these spaces: we will show that in a
neighbourhood of the special point, uniformisation is obtained by
performing a finite cover of degree $24 \times 32$, both at the source and
target, and that the uniformised period map is étale.

The stabiliser of the special period point is the group of automorphisms of
$\Z^2(2) \oplus \Z^4(-1)$, preserving the direct sum decomposition
(and the connected components of $\mathcal D_L$): it is the quotient of
$SO_2(\Z) \times O_4(\Z)$ by $\pm\id$, which is a group of order $24
\times 32$. This settles the statement for the target space.

\begin{prop}
  A stable configuration of lines has period point $\omega_0$ if and only if
  it is a complete quadrangle and the three pairs of lines are the
  opposite sides of the quadrangle.
\end{prop}

\begin{proof}
  The associated surface should be a K3 surface with involution, with
  4 disjoint rational $(-2)$-curves orthogonal to the polarisation
  $H$. Each one defines a fixed rational double point, and the only
  stable configurations having these should have exactly four distinct triple points
  through which goes a line from each pair.

  Of course no line can go through three triple points, and counting
  vertices and edges in the graph drawn by the lines, we see that
  each line goes through exactly two triple points,
  giving a complete quadrilateral.
\end{proof}

\begin{prop}
  The GIT quotient $(\mathbb P^2)^6 \sslash PSL_3$ contains exactly two
  points representing complete quadrangles, whose triple points lie
  exactly on one line from each pair $L_1 L_2$, $L_3 L_4$, $L_5 L_6$.
\end{prop}

\begin{proof}
  First note that a complete quadrangle has trivial stabiliser in
  $PGL_3$. Fixing one of the lines $\ell_1$ the combinatorics of the
  arrangement are fully determined by the two pairs of lines meeting
  on $\ell_1$: there are two possible inequivalent configurations
  (of six ordered lines).
\end{proof}

A point representing a complete quadrangle lies at the intersection of
four transverse hypersurfaces (one for each concurring triplet of
lines).

The group $(\Z_2)^3 \rtimes \mathfrak S_3$ acts on these two points
by the exact sequence
\[ 1 \to \mathfrak S_4 \to (\Z/2) \wr \mathfrak S_3 \to \Z/2 \to
1 \]
where the first map is the action of $\mathfrak S_4$ on the six pairs
of vertices ($\mathfrak S_4$ is a semi-direct product $(\Z/2)^2
\rtimes \mathfrak S_3$), and the second one the signature map.

As before let $\widetilde{\mathcal U^s}$ be the ramified cover which adds a square
root to each Plücker coordinate. The action of the torus $T$ lifts
to $\widetilde{\mathcal U^s}$ and defines a quotient $\widetilde{\Mgit[s]}$.
\begin{thm}
  The period map $\widetilde{\mathcal U^s} \to \mathcal X_L$ is locally
  liftable to $\mathcal D_L$ and the local lifts are submersive around the special
  points.
\end{thm}

\begin{proof}
  Consider the 4-parameter family of configurations:
  \[ \begin{pmatrix}
    1 &  1 &  1 &  1 & a+b &  1  \\
    1 &  1 & -1 & -1 &  1  & c+d \\
    1 & -1 &  1 & -1 & a-b & c-d \\
  \end{pmatrix} \]
  corresponding to the configuration of lines
  \[ (x+y=\pm 1, x-y = \pm 1, y + a(x+1) + b(x-1) = 0,
  x + c(y+1) + d(y-1) = 0 \]

  The relevant minors or the matrix are
  \[ m_{135} = 4b, m_{245} = -4a, m_{146} = 4d, m_{236} = -4c \]
  which define divisors with normal crossings. Hence adding a square
  root to each of the variables $a$, $b$, $c$, $d$ allows to lift the
  period map to $\widetilde{U} \to \mathcal D_L$ where $U$ is a small
  neighbourhood of the origin. Write $a = \alpha^2$, $b = \beta^2$,
  $c = \gamma^2$, $d = \delta^2$.

  Similarly, write $S_a$, $S_b$, $S_c$, $S_d$ for the four
  $(-1)$-classes which are associated to the vanishing of $a$, $b$,
  $c$, $d$. They provide a basis of $H^{1,1}(\widetilde S,
  \Z_-)$ at the special point, and they also give local coordinates for
  $\mathcal D_L$ around the associated period point.

  We now need to compute the partial derivatives of $\int_{S_\bullet}
  \omega$, with respect to the four variables: it is clear that the
  Jacobian matrix is diagonal, so it is enough to prove that
  $\partial/\partial \alpha(\int_{S_a} \omega)$ is nonzero.
  But this computation was done in section
  \ref{ssec:period-double-pt}.

  Note that since we are dealing with the lift of the period map,
  the integral structure is no longer relevant, and we can
  again work with Campedelli surfaces (and their universal cover),
  instead of Enriques surfaces, to determine the local structure of
  the period map.
\end{proof}

The action of $W_3 \times \Z_2$ lifts to $\widetilde{\Mgit[s]}$ in
the following way: $W_3$ acts by permutation matrices on
$\bigwedge^3(\bbC^6)$ (with coefficients $\pm 1$). Then its natural
lift $\widetilde W_3 \simeq (\pm 1)^{19} \rtimes W_3$
acts on the finite cover $\widetilde{\Mgit[s]} \to \Mgit[s]$
by changing signs in the square roots of the Plücker coordinates
and by the action of $W_3$ on the base.

\begin{prop}
  The stabiliser of a special configuration in $\widetilde{\Mgit[s]}$
  inside $\widetilde W_3$ is an order $24 \times 16$ group.
\end{prop}

\begin{proof}
  Remember that the local structure of $\widetilde{\Mgit[s]} \to \Mgit[s]$
  around a special configuration is the same as the double cover of $\bbC^4$
  given by the formula $(a,b,c,d) \mapsto (a^2,b^2,c^2,d^2)$.
  The number $24 \times 16$ then accounts for the actions of $(\Z/2)^4$
  and $\mathfrak S_4$ on this space, which stabilise the origin.
\end{proof}

\begin{prop}
  The action of $Q$ lifts to the ramified cover. As a consequence,
  the uniformisation $\widetilde{\Mgit[s]} \to \Mgit[s] / W_3 \times
  \pair{Q}$ has local degree
  $24 \times 32$ on small neighbourhoods of the special points.
\end{prop}

This proves the main result of this section:
\begin{thm}
  The period map $\Mgit/(W_3 \times \Z/2) \dashrightarrow \mathcal X_L^\mathrm{BB}$
  is birational.
\end{thm}

\begin{proof}
  The preimage is $\mathcal X_L$ under the period map is exactly the
  stable locus, hence the period map $\Mgit[s]
  \to \mathcal X_L$ is proper, and finite (it lifts locally to an
  étale map).

  The map $\Mgit[s] / (W_3 \times \Z/2) \to \mathcal X_L$ is finite
  and the preimage of $[\omega_0]$ consists of a single point $P$.
  Considering the diagram of germs around the special points
  \[ \xymatrix{
    (\widetilde{\Mgit[s]}, \ast) \ar[r] \ar[d] & (\mathcal D_L,
    [\omega_0]) \ar[d] \\
    (\Mgit[s] / W_3 \times \pair{Q}, \ast) \ar[r] & (\mathcal X_L,
    [\omega_0]) \\
  } \]
  we check that the period map is a local isomorphism around this point. But since the map
  is finite, the preimage of a small contractible neighbourhood of
  $[\omega_0]$ can be chosen to be a small contractible neighbourhood
  of $P$: this proves that a generic fibre consists of a single
  element, and that $\Mgit[s] / (W_3 \times \Z/2) \to \mathcal X_L$
  is an isomorphism.
\end{proof}

The previous property can be refined in a genuine isomorphism.
\begin{corol}
  The rational map of the previous theorem is an isomorphism.
\end{corol}

\begin{proof}
  We already know that the natural map $\Mgit[s] / (W_3 \times
  \Z/2) \to \mathcal X_L$ is an isomorphism, and that the ample line
  bundles $\Ocal(1)$ over $\Mgit[s]$ and the automorphic line bundle
  $\mathcal L$ over $\mathcal X_L$ are identified under this isomorphism, using
  the formula for the universal twisted 2-form. In particular, the
  sections of these line bundles on these Zariski open subsets can be
  identified.

  The projective coordinate ring of the Baily-Borel compactification
  is by definition $\bigoplus H^0(\mathcal X_L, \mathcal L^k)$ \cite{BailyBorel}
  while $\Mgit / (W_3 \times \Z/2)$ can be obtained from the ring of invariants
  \[ \Proj \bigoplus H^0(\Mgit[s], \Ocal(3k))
  ^{SL_3 \times T \rtimes W_3 \times \Z/2} \text. \]
  But sections
  over $\Mgit[s]$ are nothing more than invariant sections over the
  stable locus in the Grassmannian, which is the
  complement of a codimension 2 union of Schubert varieties: they
  extend to the whole Grassmannian, where invariant sections
  give, by definition, the coordinate ring of the GIT quotient.
\end{proof}

\section{Mixed Hodge structures and boundary configurations}

As an appendix to section \ref{ssec:boundary-period-chi}, we determine
more precisely how the Hodge structure on the twisted cohomology of
Enriques surfaces $H^2(S, \Z_-)$ degenerates when it approaches the
boundary component “of type $\chi$” which was studied there.

The asymptotic behaviour of the periods is characterised by an
isotropic sublattice $I$ in the lattice $\Z^{2,10}$ (which is
isomorphic to the generic $H^2(S, \Z_-)$). The main result is
proposition \ref{prop:mixed-hodge-on-type-chi}, which states that
$I^\perp / I$ is isomorphic to $E_8(-1)$. A consequence is that the
analogous computation in the smaller lattice $L$ (which is used for
our period map) gives a quadratic lattice $I^\perp / I \subset L / I$
isomorphic to $\Z^2(-2)$. This fact is used to determine the complete
correspondence between boundary strata of the GIT moduli space and
strata of the Baily-Borel compactification $\mathcal X_L^\mathrm{BB}$
(see section \ref{sec:boundary-details}).

\subsection{Geometric setup}

The type $\chi$ degeneration is a 1-parameter degeneration of Enriques
surfaces, with central fibre a rational surface self-intersecting
along an elliptic curve. It happens when two lines belonging to the same
ramification divisor coincide.

The type $\chi$ degeneration admits the following description: consider the
blowup $\mathcal P$ of $\mathbb P^2$ at three points $A$, $B$, $C$,
and let $L_A$, $L_B$, $L_C$ be the $(-1)$-curves arising as the proper
transforms of the sides $BC$, $CA$, $AB$. Then the linear systems $H_A
= H - A$, $H_B$, $H_C$ are pencils of rational curves with
self-intersection zero. Note that $K_{\mathcal P} = -3H+A+B+C = -H_A-H_B-H_C$.

\begin{figure}[ht]
  \begin{center}
    \includegraphics[width=.5\textwidth]{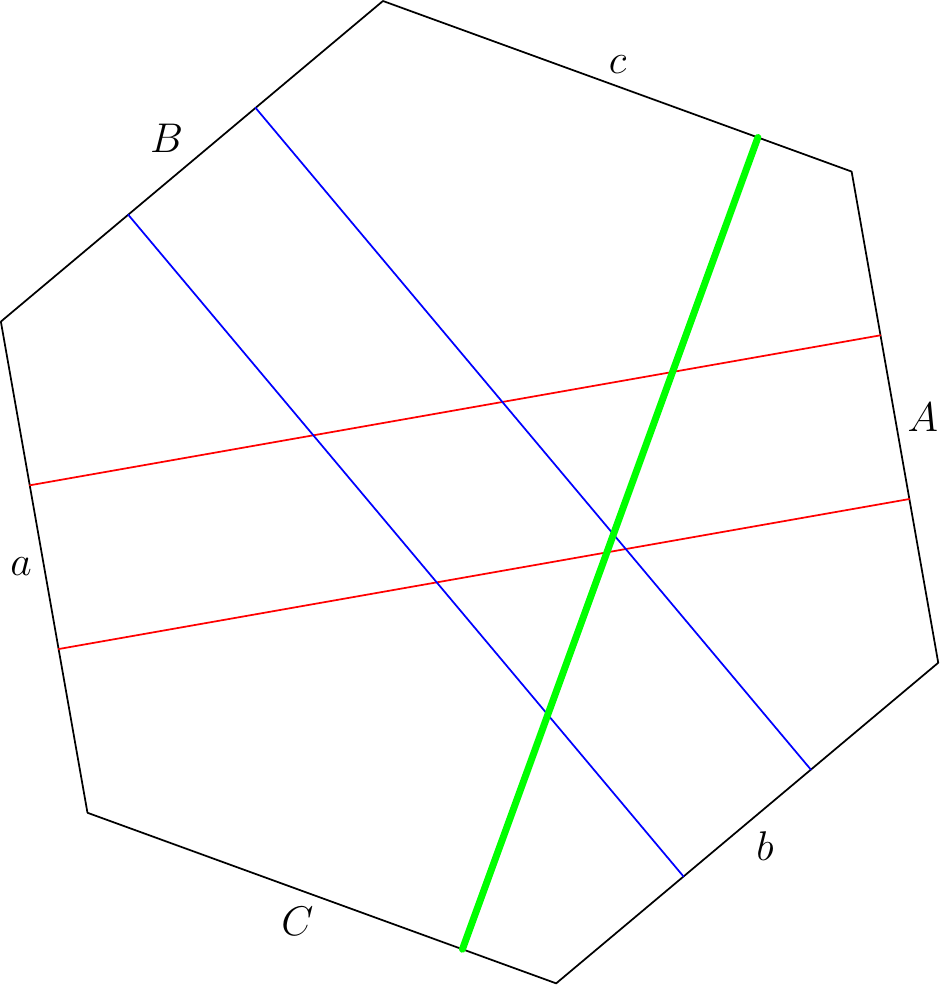}
    \caption{Type $\chi$ degeneration}
    \label{fig:degenchi}
  \end{center}
\end{figure}

A $D_{1,6}$-Enriques surface can be constructed as the bidouble cover
of $\mathcal P$ with ramification divisors of the form $R_A = D_A +
D'_A$ where $D_A$ and $D'_A$ are elements of $|H_A|$. A type I
degeneration is obtained by making $R_C$ into a double line, and we
additionally suppose that $R_A$ and $R_B$ remain generic (see figure
\ref{fig:degenchi}).

\begin{prop}
  The central fibre of a type $\chi$ degeneration of the K3 surfaces is a
  union of two isomorphic rational surfaces $V_i$, obtained as bidouble
  covers ramified over $R_A$ and $R_B$, and $K_V = -H_C$.

  The central fibre of the corresponding degeneration of Enriques
  surfaces is a rational surface whose normalisation has $K^2=0$,
  having everywhere normal crossings along a smooth elliptic curve
  $E$.
\end{prop}

Note that this degeneration can be obtained from the classical model
of 6 lines in the plane by blowing up the four points lying over the
vertex of the two collapsing lines.

\begin{prop}[see Friedman, or details in the next section]
  A type $\chi$ degeneration is a semistable type II degeneration in
  the sense of Kulikov. The limit mixed Hodge structure of a type $\chi$
  degeneration has $W_1 = H^1(E)$, $W_2/W_1 \simeq E^\perp/E$
  (where $E$ is an element of $H^2(V_1 \sqcup V_2)$).
\end{prop}

Let $V$ be the normalisation of the degenerate surface, which is a Del
Pezzo surface of degree 4 blown-up at 4 points. Note that $E = -K_V$,
hence the lattice $E^\perp$ is isomorphic to $-E_9 \simeq \mathbb Z E
\oplus -E_8$.

\begin{prop}
  \label{prop:mixed-hodge-on-type-chi}
  The mixed Hodge structure of a type $\chi$ degeneration is
  characterised by a lattice $E^\perp/E \simeq E_8(-1)$, and only
  depends (up to isomorphism) on the isomorphism class of the
  elliptic curve $E$.
\end{prop}

Since the orthogonal complement of $D_6$ in $E_8$ is isomorphic to
$\Z^2(2)$, the corresponding boundary
component of the period space corresponds to isotropic sublattices
of $L$ such that $I^\perp / I \simeq \Z^2(-2)$.

\subsection{The Clemens-Schmid exact sequence for cohomology with local
  coefficients}

Let $\mathcal S \to \Delta$ be a type $\chi$ degeneration of Enriques
surfaces, and write $\mathcal S = \mathcal T / \iota$ where
$\iota$ is an involution without fixed points. As usual, we denote by
$\mathbb Z_-$ the nontrivial local system of integers over $\mathcal
S$. The fibre of $\mathcal S$ over $t \in \Delta$ is denoted by $S_t$
and $\mathcal S^\ast$ is the fibre product $\mathcal S \times_\Delta
\Delta^\ast$.

\begin{prop}
  The homology groups of $S_0$ with coefficients in $\Z_-$ are
  $H_0 = \Z/2\Z$, $H_1 = 0$, $H_2 = \Z^{11}$,
  $H_3 = 0$, $H_4 = \Z$.
\end{prop}

\begin{proof}
  Remember that $S_0$ is obtained from its normalisation $V_0$ by
  replacing the double curve $\widetilde E$ by the singular locus $E$.
  Note that $\Z_-$ pulls back to $\Z$ on $V_0$. Then by a
  Mayer-Vietoris type exact sequence
  \[ H_n(E, \Z) \to H_n(V_0, \Z) \to H_n(S_0, \Z_-) \to -1 \]

  We get $H_0 = \Z/2$, $H_1 = 0$, $H_2$ is an extension
  of $H_1(E)$ by $H_2(V_0, \Z) / [E]$, $H_3 = 0$, $H_4 = \Z$.
\end{proof}

\begin{prop}
  The cohomology groups of $S_0$ with coefficients in $\Z_-$ are
  $H^0 = 0$, $H^1 = \Z/2$, $H_2 = \Z^{11}$,
  $H^3 = 0$, $H^4 = \Z$.
\end{prop}

\begin{proof}
  From the sequence of sheaves
  \[ 0 \to (\Z_-)_{S_0} \to \nu_\star \Z_{V_0} \to \Z_E \to 0 \]
  we derive the long exact sequence
  \[ H^n(S_0, \Z_-) \to H^n(V_0, \Z) \to H^n(E, \Z) \to +1 \]

  We get $H^0 = 0$, $H^1 = \Z/2$, $H^2$ is an extension
  of $H^2(V_0)$ by $H^1(E)$, $H^3 = 0$, $H^4 = \Z$.
\end{proof}

\begin{prop}
  \index{suite exacte de Clemens-Schmid}
  The Clemens-Schmid exact sequence
  \[ H_4(S_0, \frac 1 2 \Z_-) = \Z E \to H^2(S_0, \Z_-) \to
  H^2_\mathrm{lim}(S_t, \Z_-) \xrightarrow N H^2_\mathrm{lim} (S_t,
  \mathbb Z_-) \]
  is exact over the integers.
\end{prop}

\begin{proof}
  Proceeding as in \cite{Friedman}, we decompose the Clemens-Schmid
  exact sequence into a part of Wang's long exact sequence
  \[ H^1(S_t, \Z_-) \to H^2(S^\ast, \Z_-) \to H^2(S_t, \Z_-)
  \xrightarrow{T-1} H^2(S_t, \Z_-) \]
  and the relative cohomology exact sequence
  \[ H^1(\mathcal S^\ast, \Z_-) \to H^2(\mathcal S, \mathcal S^\ast;
  \Z_-) \to H^2(\mathcal S, \Z_-) \to H^2(\mathcal S^\ast, \Z_-) \to
  H^3(\mathcal S, \mathcal S^\ast, \Z_-) \]

  By Poincaré duality, $H^3(\mathcal S, \mathcal S^\ast, \Z_-) \simeq
  H_3(S_0, \Z_-) = 0$. Hence $H^2(\mathcal S^\ast, \Z_-)$ is the
  cokernel of $H^2(\mathcal S, \mathcal S^\ast; \Z_-) \to H^2(S,
  \Z_-)$ which is also $H_4(S_0, \Z_-) \to H^2(\mathcal S, \Z_-)$.

  Moreover, since the image of a generator of $H_4(S_0, \Z_-)$ is
  $\pm 2E$, the 2-torsion element in $H^2(\mathcal S^\ast, \Z_-)$ is $[E]$.

  Now $S_t$ is a honest Enriques surface, hence $H^2(S_t, \Z_-)$
  is torsion-free, and $H^1(S_t, \Z_-)$ is isomorphic to $\Z/2$, hence the kernel
  of $T-1$, or equivalently the monodromy operator $N$, is the
  quotient $H^2(S^\ast, \Z_-) / [E]$.
\end{proof}

As a conclusion, $Gr^W_2 H^2_\mathrm{lim}(S_t, \Z_-)$ is isomorphic
to $H^2(S_0, \Z_-) / [E], H^1(E, \Z)$. But
\[ H^2(S_0, \Z_-)/H^1(E, \Z) \simeq E^\perp \subset H^2(V_0, \Z)
\text.\]

\section{The monodromy of marked Campedelli surfaces}

The results we obtained concerning Enriques surfaces can be summarised
by the theorem stated in the introduction:
\begin{thm*}
  The structure of the period map of Campedelli surfaces can be
  described by the following diagram
  \[ \xymatrix{
    \Gr(4,7) \ar@{-->}[r] \ar[d] & \Gr(3,6) \ar[d] & \\
    \Gr(4,7) / \mathfrak S_4 \ar@{-->}[r] &
    \Gr(3,6) / (\Z/2)^2 \rtimes \mathfrak S_3 \ar[r] &
    \frac{\Gr(3,6)}{(\Z/2 \wr \mathfrak S_3) \times \pair{Q}}
    \simeq \mathcal X_L^\mathrm{BB} \\
  } \]
  where $\mathcal X_L^\mathrm{BB}$ is the Baily-Borel compactification
  of $\mathcal X_L$.
\end{thm*}

In particular, the fibres of the period map we defined have
generically four connected components. It is natural to ask whether
the same statement is true when considering the period map
of the integral Hodge structure $H^2(X, \Z_\kappa)$ (we call
\emph{Campedelli lattice} the underlying integral quadratic lattice),
instead of $H^2(S, \Z_-)$ (whose lattice is isomorphic to $L$.

The main result of this section is corollary
\ref{corol:auto-enriques-campedelli-lats}: the embedding of the
lattice $L$ in $H^2(X, \Z_\kappa)$ is invariant under isometries of
$L$, meaning that any isometry of $L$ induces an isometry of the
overlattice $H^2(X, \Z_\kappa)$ (it is not true, however, that $L$ is
invariant under any isometry of $H^2(X, \Z_\kappa)$. In particular,
the isomorphism class of the Hodge structure $H^2(S, \Z_-)$ fully
determines that of $H^2(X, \Z_\kappa)$, meaning that the
true period map of marked Campedelli surfaces also has disconnected fibres.

\subsection{The cohomology lattice of Campedelli surfaces}

Suppose $X$ is a smooth Campedelli surface, with a chosen involution $s_\kappa$,
and $X/s_\kappa = S$. By blowing up the six fixed points, we get a
double cover $\tilde X \to \tilde S$ ramified over six $(-2)$-curves
and a genus 3 curve, where $S$ is the minimal resolution of $S_\kappa$
and $\tilde X$ is $X$ blown up at the six isolated fixed points of $s_\kappa$.
The results of section \ref{ssec:index-comp} apply, and
the index of $H^2(\tilde S, \Z_-)(2) \subset H^2(\tilde X,
\Z_\kappa)^{s_\kappa}_\mathrm{num}$ is $2^{0-2+7} = 2^5$. The notation
$H_\mathrm{num}$ still denotes the quotient of $H$ by its torsion subgroup.
We will see in the next section that
$H^2(\tilde S, \Z_-)$ contains $L \oplus D_6(-1)$ as an index 4
subgroup: consider the following maps
\[ L(2) \oplus D_6(-2) \to H^2(\tilde S, \Z_-)(2) \to H^2(X,
\Z_\kappa)^{s_\kappa}_\mathrm{num} \oplus \Z^6(-1) \]
where the first one has index four, the second one index $2^5$.

Then the equation
\[ 4 \cdot 2^5 = [\Z^6(-1):D_6(-2)] \times
[H^2(X, \Z_\kappa)^{s_\kappa}_\mathrm{num}:L(2)] \]
gives that $L(2) \subset H^2(X,\Z_\kappa)^{s_\kappa}_\mathrm{num}$ has
index $2^3$, and the lattice $L_0 =
H^2(X,\Z_\kappa)^{s_\kappa}_\mathrm{num}$ has also discriminant 4.
Our goal is now to characterise $L_0$ as an integral overlattice of
$L(2)$: as such, it is completely determined by the image $\Lambda_0$
of $L_0$ in $L(2)^\ast/L(2) = \Lambda$.

\begin{lemma}
  Let $\Lambda = \Z/4(4)^2 \oplus \Z/2(-2)^4$ be the discriminant group
  of $L(2)$, with its fractional bilinear form
  \[ b(x,y) = \frac{x_1 y_1 + x_2 y_2}{4} - \frac
  {x_3 y_3 + x_4 y_4 + x_5 y_5 + x_6 y_6}{2} \pmod 1 \text. \]

  An index $2^3$ integral overlattice of $L(2)$ correspond to an
  isotropic subgroup of order 8 of $\Lambda$. Elements with
  integral norm form the subgroup $\Lambda_i$ of order 64
  \[ \set{(a,b,c,d,e,f) \text{ such that }
    a \equiv b \pmod 2, s+c+d+e+f \equiv 0 \pmod 2} \]
  where $s$ is the common parity of $a$ and $b$.
\end{lemma}

\begin{proof}
  The decomposition $L(2) = \Z^2(4) \oplus \Z^4(-2)$ induces the
  description of $\Lambda$ as a product of cyclic groups, and the
  description of isotropic vectors results from the fact
  that $(x^2 + y^2)$ is even if and only if $x \equiv y \pmod 2$.
\end{proof}

The determination of $\Lambda_0$ uses a strong symmetry property of
the lattice $L_0$: from this point, we identify $L$ with the (limit of
the) lattice $H^{1,1}(\Z_-)$ at the special point $[\omega_0]$
representing a configuration with four triple points.  Consider a
uniformising disk $\mathbb D$ around $[\omega_0]$ in the space
$\widetilde{\Mgit[s]}$. This disk supports a global identification
of the twisted Picard lattice of a general member with $L$.
The $\mathfrak S_4$-action on the configuration,
extends to an action on the disk (which at first order is equivalent
to the linear action of $\mathfrak S_4$ on four coordinates).
As in section \ref{sec:most-special-config}, the group $\mathfrak
S_4$ acts on the six lines of the configuration by its embedding
$\mathfrak S_4 \subset (\Z/2) \wr \mathfrak S_3$: this action
corresponds to the permutation of coordinates which we used to define
Campedelli surfaces and their quotients (via the isomorphism
$\mathfrak S_4 \simeq GA_2(\mathbb F_2)$ with the group of affine
transformations of $\mathbb A^2(\mathbb F_2)$).

\begin{prop}
  The $\mathfrak S_4$-action on $\mathbb D$ lifts to the family of
  Enriques surfaces supported by the complement of the triple point
  locus, hence to a well-defined automorphism of $L$. \qed
\end{prop}

The action of $\mathfrak S_4$ on the four triple points of the
configuration gives the following property of the action
on the cohomology:
\begin{prop}
  The action of $\mathfrak S_4$ on the lattice induces the natural
  action of $\mathfrak S_4$ on the $\Z^4(-1)$ summand (up to choices
  of signs).
\end{prop}

We are then looking for the particular isotropic subgroup $\Lambda_0$
corresponding to the overlattice of Campedelli surfaces, which is also
globally defined on $\mathbb D$: there exists a family of Campedelli
surfaces, corresponding to the choices of a seventh line, whose base maps to
$\mathbb D$ with connected fibres, and to which the $\mathfrak
S_4$-action can be lifted, using the same action on coordinates.
This family defines a uniform choice of
$\Lambda_0$ in $\Lambda$ over $\mathbb D$.
\begin{prop}
  The subgroup $\Lambda_0$ is invariant under the $\mathfrak S_4$
  action. \qed
\end{prop}

\begin{prop}
  The projection of $\Lambda_0$ of the second summand $(\Z/2)^4$
  consists of even vectors (i.e. the sum of their coordinates is even).
\end{prop}

\begin{proof}
  A $\mathfrak S_4$-invariant subspace of $(\Z/2)^4$ contains all vectors of a
  given weight (number of nonzero coefficients). Since the projection
  of $\Lambda_0$ may only contain at most 8 elements, there cannot be
  a weight one vector (basis vector), and there cannot be a weight
  three vector either (the sum of three different weight 3 vectors is
  a weight 1 vector). The subgroup of even vectors contains exactly
  eight elements.
\end{proof}

In particular, $\Lambda_0$ only contains elements $(a,b,c,d,e,f)$
such that $c+d+e+f$ is even, then by the characterisation of isotropic
vectors given above, $a$ and $b$ are both even. But since the pairing
of $(2a,2b,1,1,0,0)$ and $(2a',2b',0,1,1,0)$ is $1/2$ (thus not an
integer), $\Lambda_0$ cannot actually contain even vectors except
for the \emph{characteristic element} $(1,1,1,1)$: it is the only
non zero element of $\Lambda$ which is orthogonal to the 2-torsion subgroup
of $\Lambda$.

\begin{prop}
  The group $\Lambda_0$ is $\pair{2a, 2b, c+d+e+f}$ (where the letters
  stand for the standard generators of $\Lambda$).
\end{prop}

\begin{proof}
  The discussion has shown that $\Lambda_0$ was a subgroup of
  $\pair{2a, 2b, c+d+e+f}$, which has order 8.
\end{proof}

This implies that $L_0$ is absolute in the sense that it is invariant
under the action of $O(L)$, in other words any element of $O(L)$
naturally induces an element of $O(L_0)$. This is because $\Lambda_0$
is generated by the 2-divisible elements of $\Lambda$, and by the
characteristic vector $c+d+e+f$. This leads to the following
description:
\begin{thm}
  The cohomology lattice of a linear system $\Z_\kappa$ on a
  Campedelli surface splits as an invariant part of signature $(2,4)$
  and an anti-invariant part which is negative definite.

  The invariant part is isomorphic to $\Z^2 \oplus D_4(-1)$,
  is which $L(2) \simeq \Z^2(4) \oplus \Z^4(-2)$ is embedded
  in the most natural way (as computed above).
\end{thm}

\begin{corol}
  \label{corol:auto-enriques-campedelli-lats}
  The automorphism group $\Aut(L)$ is naturally embedded in the
  orthogonal group of the Campedelli lattice, i.e. any automorphism
  of $L$ extends to an isometry of the Campedelli lattice. Two marked Campedelli
  surfaces with the same associated Enriques surface have the same
  invariant Hodge structure $H^2(\Z_\kappa)$.

  In particular marked Campedelli surfaces with the same periods form
  a disconnected family.
\end{corol}

\section{Details on the boundary period map}
\label{sec:boundary-details}

The list of boundary components in both the GIT moduli space and the
Baily-Borel compactification of the period space are listed in figure
\ref{fig:tab-boundary}. The description of the strata is given by the
lines which coincide, the lines of each pair being denoted by $L_A$
and $L'_A$ for example. The entries in the table result from the
discussion below: this section is dedicated to the classification
of primitive isotropic sublattices of $L$ up to isometries.

\begin{figure}[ht]
\begin{center}
\begin{tabular}{|l|c|l|l|l|}
  \hline
  GIT component & Dimension & Isotropic $\ell$ & $\ell^\perp / \ell$ \\
  \hline
  $L_A = L_A'$ (type $\chi$) & 1 & Odd plane & $\Z^2(-2)$ \\
  $L_A = L_B$ & 1 & Even plane & $\Z^2(-1)$ \\
  $L_A = L_A'$, $L_B=L_B'$ & 0 & Odd of type 2 & $H \oplus \Z^2(-2)$ \\
  $L_A=L_A'$, $L_B=L_C$ & 0 & Odd of type 1 & $\Z^{1,1} \oplus \Z^2(-2)
  = \Z^{1,1}(2) \oplus \Z^{0,2}$ \\
  $L_A=L'_B$, $L_B=L'_C$... & 0 & Even & $\Z^{1,3}$ \\
  \hline
\end{tabular}
\caption{List of boundary components in GIT and Baily-Borel compactifications}
\label{fig:tab-boundary}
\end{center}
\end{figure}

\subsection{Boundary components and local structure}

We give here an example of how the one-dimensional boundary components
can be explicitly identified.

\begin{prop}
  The boundary period map $\mathcal M_\chi / (\Z/2 \wr \mathfrak S_3)
  \to X(2A_1)$ is birational.
\end{prop}

We know from the study of type $\chi$ degenerations that $\mathcal
M_\chi$ is mapped to $\mathcal X(2A_1)$. In the next section, we
prove that $\mathcal X(2A_1)$ is the quotient of the upper half-plane
under the action of an index three subgroup of $PSL_2(\Z)$. An
elementary argument can be used to reprove this using almost only
lattice arithmetic.

\begin{proof}
  For generic $\chi$ configurations, $\mathcal M_\chi$ is isomorphic
  to $\mathbb P^1$, the isomorphism being given by the cross-ratio
  $t$ of the ordered four lines crossing the double line. The singular
  locus is an elliptic curve whose ramification points can be chosen
  to be $(\pm 1, \pm \sqrt t)$: its $j$-invariant is a degree 6
  rational function of $t$.

  The group $(\Z/2 \wr \mathfrak S_2)$ acts on the cross-ratio with an
  order 4 kernel. Hence the $j$-invariants $\mathcal M_\chi / (\Z/2
  \wr \mathfrak S_3) \to \mathbb P^1$ and $X(2A_1) \to \mathbb P^1$
  both have degree 3, hence the period map has degree one.
\end{proof}

\subsection{Isotropic sublattices of the cohomology of
  \texorpdfstring{D$_6$}{D6} Enriques surfaces}

Denote by $L$ the lattice $\Z^2(2) \oplus \Z^4(-1)$. Remember that $L$
embeds as an index two sublattice of $\Z^{2,4}$.

\subsubsection{Isotropic vectors}

\begin{lemma}
  If $\ell$ is a primitive isotropic vector in $L$, the set of
  $\pair{\ell,x}$ for $x \in L$ is either $\Z$, either $2\Z$.
  According to this distinction, we say $\ell$ is odd or even.
\end{lemma}

\begin{proof}
  Let $d$ be the gcd of all $\pair{\ell, x}$. If
  $\ell = (a_1, a_2; b_1, b_2, b_3, b_4)$ in the coordinates
  where $L = \Z^2(2) \oplus \Z^4(-1)$, then $d$ is the gcd
  of $2a_i$ and $b_i$.

  Since the gcd of $a_i$ and $b_i$ is one, the conclusion follows.
\end{proof}

\begin{lemma}
  If $\ell$ is an even vector, there exists a decomposition
  $L = \Z^{1,1}(2) \oplus \Z^{1,3}$ such that $\ell$ belongs to the
  first summand. Up to automorphism of $L$, there is only one
  even primitive isotropic vector.
\end{lemma}

\begin{lemma}
  If $\ell$ is an odd vector, then $\ell$ belongs to a unimodular
  sublattice of $L$ of signature $(1,1)$. The corresponding
  possible decompositions of $L$ are:
  \[ \Z^{1,1} \oplus (\Z^{1,1} \oplus \Z^2(-2)) \]
  \[ H \oplus (\Z^{1,1} \oplus \Z^2(-2)) \]
  \[ \Z^{1,1} \oplus D_{1,3} \]
  giving two orbits (the first and the second decompositions occur for
  the same vectors). We will say the first case is \emph{type 1},
  the last case being \emph{type 2}. The orbits are characterised
  by $\ell^\perp / \ell$, which is either $\Z^{1,1} \oplus \Z^2(-2)$
  or $H \oplus \Z^2(-2)$.
\end{lemma}

\begin{proof}
  Let $(\ell, m_1, \dots, m_4)$ be a basis of $\ell^\perp$ and
  $\nu$ be such that $\pair{\ell, \nu}$ is minimal. In the basis,
  $(\ell, \nu, m_\bullet)$ the quadratic form has matrix
  \[ \begin{pmatrix}
    0 & \eps & 0\\
    \eps & \nu^2 & N \\
    0 & N & A \\
  \end{pmatrix} \]
  where $A$ is the matrix of the quadratic form on $\bigoplus
  \Z m_\bullet$ and $N_\bullet$ is $\pair{\nu, m_\bullet}$.

  Its determinant is $4 = - \eps^2 \det A$. If $\eps = 2$,
  i.e. $\ell$ is even, $A$ is unimodular, hence $m_\bullet$
  span a unimodular lattice whose orthogonal complement
  contains $\ell$ and has discriminant 4: its matrix
  should have the form $\big(\begin{smallmatrix}
    0 & 2 \\ 2 & n \\ \end{smallmatrix} \big)$. Since it has the same
  discriminant form as $L$, it must be an index two sublattice of
  $H$ or $\Z^{1,1}$. Then it is $\Z^{1,1}(2)$ (since $H(2)$ has
  two isotropic elements for its discriminant bilinear form).

  If $\eps = 1$, then $\Z \ell \oplus \Z \nu$ is unimodular, and
  $(m_\bullet)$ should be a index 2 lattice in $\Z^{1,3}$,
  hence the decompositions mentioned in the statement.
\end{proof}

\begin{prop}
  The lattice $L$ contains three orbits of primitive isotropic vectors.
  The corresponding vectors are even, odd of type 1, odd of type 2.
  They are distinguished by the fact that $\ell^\perp/\ell$
  is isomorphic to $\Z^{1,3}$ (resp. $\Z^{1,1} \oplus \Z^2(-2)$,
  $D_{1,3} = H \oplus \Z^2(-2)$).
\end{prop}

\subsubsection{Isotropic planes}

\begin{prop}
  The lattice $L$ contains two orbits of maximal isotropic
  sublattices, distinguished by the fact that they contain (or not)
  even isotropic vectors.

  More precisely, if $\lambda$ is an odd isotropic lattice of rank 2,
  then there exists a decomposition $L = \Z^{2,2} \oplus \Z^2(-2)$
  such that $\lambda$ is a maximal isotropic lattice of the first
  summand. Primitive isotropic vectors in $\Z^{2,2}$ are either
  odd of type 1, e.g. $(1,0; 1, 0)$, or odd of type 2, e.g.
  $(1,1; 1,1)$ (whose orthogonal complement is even).

  If $\lambda$ is even (i.e. contains even isotropic vectors),
  it is maximal isotropic in a $\Z^2 \oplus
  \Z^2(-2)$ summand of $L$, and $\lambda^\perp/\lambda$ is isomorphic
  to $\Z^2(-1)$. It contains primitive even vectors and odd vectors of
  type 1.
\end{prop}

\begin{proof}
  Let $\lambda$ be a rank 2 isotropic sublattice of $L$. As before,
  choose a basis
  \[ (\lambda_1, \lambda_2, \nu_1, \nu_2, \mu_1, \mu_2) \]
  such that $\nu_\bullet$ span $\lambda^\perp/\lambda$. Then the
  matrix of the quadratic form of $L$ can be written with
  $(2\times2)$-sized blocks:
  \[ \begin{pmatrix}
    0 & 0 & A \\
    0 & B & C \\
    A^T & C^T & D \\
  \end{pmatrix} \]
  and the discriminant is $4 = \det B (\det A)^2$. Then either
  $\det B = 1$ or $\det A = 1$.

  If $\det A = 1$, then $\lambda$ is embedded in the unimodular
  lattice spanned by $\lambda$ and $\mu$, which is indefinite hence
  isomorphic to $H^2$ or $\Z^{2,2}$. But the case $H^2$ is
  impossible, since it would imply that $L$ is even. In this case
  $L$ contains no primitive even isotropic vectors (it would imply
  that $\det A$ is even).

  If $\det B = 1$, then $\lambda \oplus \nu$ has signature $(2,2)$
  and the same discriminant form as $L$, and $B$ is definite
  of rank 2, and unimodular, hence isomorphic to $\Z^2$.
  Since $\lambda \oplus \nu$ cannot be isomorphic to $D_{2,2}$,
  which has the wrong bilinear form, it must be isomorphic to
  $\Z^2 \oplus \Z^2(-2)$.
\end{proof}

The study of isotropic sublattices of $L$ allows to determine the
structure of the boundary components of the Baily-Borel-Satake
compactification of the period space.

\begin{prop}
  The boundary of $\mathcal D_L$ consists of 3 distinguished points
  $p_\mathrm{even}$, $q_1$, $q_2$ and two rational curves $C_\mathrm{even}$
  going through $p_\mathrm{even}$ and $q_1$, $C_\mathrm{odd}$ going
  through $q_1$ and $q_2$.
\end{prop}

\begin{proof}
  The identification of dimension 1 strata results from the study of
  type $\chi$ degenerations. The identification between strata of
  dimension zero results from the incidence relations with the
  one-dimensional strata.
\end{proof}

\subsection{Modular curves at the boundary of the period space}

In the following, $\gamma$ will denote the index three subgroup of
$GL_2(\Z)$ of matrices $\big( \begin{smallmatrix} a & b \\ c & d \\
\end{smallmatrix} \big)$ such that $a \equiv d$ and $b \equiv c \pmod 2$.

\subsubsection{Boundary curve for odd isotropic lattices}

\begin{prop}
  Let $I$ be a primitive isotropic plane in $L \simeq \Z^{2,2} \oplus
  \Z^2(-2)$, such that $I^\perp/I \simeq \Z^2(-2)$. Then up to
  isometry, we can assume that $I \subset \Z^{2,2}$.

  The image of the stabiliser of $I$ in $\Aut(L)$ is an index three
  subgroup $\gamma(I)$ of $GL_2(I)$, which is identified with $\gamma$
  for a suitable choice of basis.
\end{prop}

Before explaining the proof, we choose a basis of $L$ which gives
bases for $I$, $I^\perp/I$, and $L/I^\perp \simeq I^\vee$ (and assume
that the intersection matrix of $I$ and $I^\vee$ is given by
the identity matrix). Then an isometry of $L$ preserving $I$ has a
matrix of the following form
\[ \begin{pmatrix}
  H & A & B \\
  0 & G & C \\
  0 & 0 & H^\dagger \\
\end{pmatrix} \]
where $H^\dagger$ is the inverse transpose of $H$, and $G \in O_2(\Z)$.
In order to define an isometry, the matrix has to satisfy additionally
\begin{itemize}
\item $A^T H^\dagger - 2 G^T C = 0$;
\item $(H^\dagger)^T H^\dagger + B^T H^\dagger + H^{-1} B - 2 C^T C = I$.
\end{itemize}
where $I$ is the identity matrix.

\begin{proof}
  Choose a basis $(e_1, e_2; f_1, f_2)$ of $\Z^{2,2}$ such that $I$ is
  generated by $(e_1+f_1, e_2+f_2)$. Then choosing any orthogonal
  basis $(g_1,g_2)$ for $\Z^2(-2)$ and $(e_1, e_2)$ for a basis of $I^\vee$ gives
  a basis
  \[ (e_1+f_1,e_2+f_2,g_1,g_2,e_1,e_2) \]
  of $L$ satisfying the axioms above.

  With the previous notations,
  \[ (H^\dagger)^T H^\dagger = I + (H^{-1} B) + (H^{-1} B)^T \pmod 2 \]
  thus is either the identity matrix or $\big(\begin{smallmatrix} 1 & 1 \\ 1 &
    1 \\ \end{smallmatrix} \big)$. Hence $H$ must map to an element of
  the orthogonal group $O_2(\mathbb F_2)$ which has index three in
  $GL_2(\mathbb F^2)$ ($H$ is an element of $\gamma$).

  Conversely, if $H \pmod 2$ is in $O_2(\mathbb F_2)$, then write
  $(H^\dagger)^T (H^\dagger) = I - 2N$ where $N$ is an integral
  symmetric matrix. Then $A=0$, $B=HN$ and $C=0$, $G=I$
  defines an isometry of $L$, whose restriction to $I$ is given by $H$.
\end{proof}

Let $\Gamma_I$ be the stabiliser of $I$ in $\Aut L$
(notations are borrowed from \cite{Looijenga-compact2}) and let
\[ 1 \to N_I \to \Gamma_I \to \gamma_I \times O_2(\Z) \to 1 \]
be its Levi decomposition. Note also that $N_I$ is a central extension
\[ 1 \to \Z \to N_I \to \Hom(I^\vee, I^\perp/I) \to 0 \]
where $\Z$ is the subgroup of matrices
\[ \begin{pmatrix}
  1 & 0 & B \\
  0 & 1 & 0 \\
  0 & 0 & 1 \\
\end{pmatrix} \]
where $B$ is integral and skew-symmetric (which is central in $N_I$).

Let $\mathcal D_L$ be the period domain associated to the lattice $L$,
and $\mathcal D_L^\mathrm{BB}$ the topological space defined by Satake,
such that the Baily-Borel compactification $\mathcal X_L^\mathrm{BB}$ is
$\mathbb D_L^\mathrm{BB} / \Aut(L)$. Then the boundary curve $\mathcal
X_L(2A_1)$ corresponding to odd isotropic planes is isomorphic to
$\mathfrak H / \gamma$, where $\gamma$ is the index three subgroup of
$SL_2(\mathbb Z)$ as above.

\begin{prop}
  Let $\Gamma^I$ be the kernel of $\Gamma_I \to \gamma_I$. Then in a
  neighbourhood of a generic point of $\mathcal X_L(2A_1)$, $\mathcal
  X_L^\mathrm{BB}$ is locally isomorphic to $\mathcal D_L^\mathrm{BB} /
  \Gamma^I$.
\end{prop}

Following Looijenga \cite{Looijenga-compact2}, let $\pi_W$ be the
projection $L \to L/W$, and consider the morphisms
\[ \mathcal D_L \to \pi_I \mathcal D_L \to \pi_{I^\perp} \mathcal D_L
\]
and the filtration
\[ 0 \to \Z \to N_I \to \Gamma^I \to \Gamma_I \]
with quotients $\Z$, $\Hom(I^\perp/I, I)$, $O_2(\Z)$, $\gamma_I$.

Note that $\mathcal D_L \to \pi_I \mathcal D_L$ is a bundle of upper
half-planes, acted on by the subgroup acting as the identity on $L/I$,
which is $\Z \subset N_I$, acting by translations.

\begin{prop}
  $\mathcal D_L / \Z \to \pi_I \mathcal D_L$ is a bundle
  of punctured disks.
\end{prop}

Now consider $\pi_I \mathcal D_L \to \pi_{I^\perp} \mathcal D_L$.
Note that $\pi_I \mathcal D_L$ is the set of images of positive planes
inside $\Gr(2, L/I \otimes \R)$ or $\Gr(2, I^\perp_\R)$, which is a
principal bundle over $\pi_{I^\perp} \mathcal D_L$ under action of
$\Hom_\R(I, I^\perp / I)$ (which is the stabiliser of $I$ in
$GL(I^\perp)$, which acts by isometries).

\begin{prop}
  $\pi_I \mathcal D_L / N_I \to \pi_{I^\perp} \mathcal D_L$
  is a bundle of abelian surfaces.
\end{prop}

Here $\pi_{I^\perp} \mathcal D_L$ is the upper half plane associated
to $L / I^\perp \simeq I^\vee$ which is acted on by $\gamma_I$, it
maps to the boundary curve $\mathcal X_L(2A_1)$.

\begin{prop}
  There are natural projections
  \[ \mathcal D_L / N_I \to \pi_I \mathcal D_L / N_I
  \to \pi_{I^\perp} \mathcal D_L \]
  such that the second morphism is a bundle of abelian surfaces,
  and the first morphism is a bundle of punctured disks, and $O_2(\Z)$
  acts properly on the first two spaces.

  The local structure of the Baily-Borel compactification is obtained
  by extending $\mathcal D_L / \Gamma_I^o$ into a disk bundle, and
  contracting central fibres, then taking the quotient under
  $O_2(\Z)$.
\end{prop}

\subsubsection{Boundary curve for even isotropic lattices}

\begin{prop}
  Let $I$ be a primitive isotropic plane in $L \simeq \Z^2 \oplus
  \Z^2(-2) \oplus \Z^2(-1)$, such that $I^\perp/I \simeq \Z^2(-1)$.
  Then up to isometry, we can assume that $I \subset \Z^2 \oplus
  \Z^2(-2)$.

  The image of the stabiliser of $I$ in $\Aut(L)$ has
  index three in $SL_2(I)$.
\end{prop}

Note that a primitive isotropic plane in $\Z^2 \oplus \Z^2(-2)$
(in a canonical basis $(e_1, e_2; f_1, f_2)$) is given by
$I = \pair{e_1+e_2+f_1, e_1-e_2+f_2}$. Let $g_1, g_2$ be a orthonormal
basis of the remaining summand $\Z^2(-1)$, and consider the basis of
$L$ given by
\[ e_1+e_2+f_1, e_1-e_2+f_2, g_1, g_2, e_1, e_2 \]
The matrix of the quadratic form is then
\[ \begin{pmatrix}
  0 & 0 & M \\
  0 & -1 & 0 \\
  M & 0 & 1 \\
\end{pmatrix} \]
where $M = \big(\begin{smallmatrix} 1 & 1 \\ 1 & -1 \\ \end{smallmatrix}
\big)$. A matrix
\[ g = \begin{pmatrix}
  H & A & B \\
  0 & K & C \\
  0 & 0 & L \\
\end{pmatrix} \]
then defines an isometry if and only if
\begin{itemize}
\item $K \in O_2(\Z)$;
\item $H^T M L = M$;
\item $A^T M L - K^T C = 0$;
\item $B^T M L + L^T M B + L^T L - 2 C^T C = I$.
\end{itemize}

\begin{lemma}
  An explicit computation shows that if $X = MYM/2$, and $X$, $Y$ are
  elements of $GL_2(\Z)$, then $X$ and $Y$ belong to $\gamma$.
\end{lemma}

\begin{proof}
  Write $Y = \big(\begin{smallmatrix} a & b \\ c & d
    \\ \end{smallmatrix} \big)$. Then
  \[ X = \frac 1 2 \begin{pmatrix}
    a+b+c+d & a-b+c-d \\
    a+b-c-d & a-b-c+d \\
  \end{pmatrix} \]
  which implies that $a+b+c+d$ is even, that is $Y \in \gamma$.
  It is then easy to check that $X \in \gamma$.
\end{proof}

\begin{proof}[Proof of the proposition]
  From the above equations we get $H^T = M L^{-1} M / 2$,
  hence $H$ is an element of $\gamma$.

  Conversely, given any matrix $H$ mapping to $O_2(\mathbb F_2)$,
  i.e. $H \in \gamma$, there exists an isometry
  \[ g = \begin{pmatrix}
    M H^\dagger M / 2 & 0 & B \\
    0 & 1 & 0 \\
    0 & 0 & H \\
  \end{pmatrix} \]
  where $B^T M H + H^T M B + H^T H = I$. Since $I - H^T H = 2N$ for
  some integral symmetric matrix $N$. Note that diagonal coefficients
  of $N$ are even, hence either $N$ or $N' = N + \big(\begin{smallmatrix} 0
    & 1 \\ -1 & 0 \\ \end{smallmatrix} \big)$ is even.

  Setting $B = M H^\dagger N / 2$ or $M H^\dagger N' / 2$ gives $g$.
\end{proof}



\begin{bibsection}
\begin{biblist}

\bib{AlexeevPardini}{article}{
  author={Alexeev, V.},
  author={Pardini, R.},
  title={Explicit compactifications of moduli spaces of Campedelli and
    Burniat surfaces},
  date={2009},
  eprint={\arXiv{0901.4431}},
}

\bib{Allcock}{article}{
   author={Allcock, Daniel},
   title={The period lattice for Enriques surfaces},
   journal={Math. Ann.},
   volume={317},
   date={2000},
   number={3},
   pages={483--488},
   issn={0025-5831},
   review={\MR{1776113 (2002a:14040)}},
   eprint={\arXiv{math/9905166}},
}

\bib{Atiyah-resolution}{article}{
   author={Atiyah, M. F.},
   title={On analytic surfaces with double points},
   journal={Proc. Roy. Soc. London. Ser. A},
   volume={247},
   date={1958},
   pages={237--244},
   issn={0962-8444},
   review={\MR{0095974 (20 \#2472)}},
}

\bib{AtiyahSinger}{article}{
   author={Atiyah, M. F.},
   author={Singer, I. M.},
   title={The index of elliptic operators. III},
   journal={Ann. of Math. (2)},
   volume={87},
   date={1968},
   pages={546--604},
   issn={0003-486X},
   review={\MR{0236952 (38 \#5245)}},
}

\bib{BailyBorel}{article}{
  author={Baily, W. L. Jr.},
  author={Borel, Armand},
  title={Compactification of arithmetic quotients of bounded symmetric domains},
  journal={Ann. of Math. (2)},
  volume={84},
  date={1966},
  pages={442--528},
  issn={0003-486X},
}

\bib{BarthPeters}{article}{
   author={Barth, W.},
   author={Peters, C.},
   title={Automorphisms of Enriques surfaces},
   journal={Invent. Math.},
   volume={73},
   date={1983},
   number={3},
   pages={383--411},
   issn={0020-9910},
   review={\MR{718937 (85g:14052)}},
   doi={10.1007/BF01388435},
}

\bib{BHPV}{book}{
   author={Barth, Wolf P.},
   author={Hulek, Klaus},
   author={Peters, Chris A. M.},
   author={Van de Ven, Antonius},
   title={Compact complex surfaces},
   series={Ergebnisse der Mathematik und ihrer Grenzgebiete. 3. Folge. A
   Series of Modern Surveys in Mathematics},
   volume={4},
   edition={2},
   publisher={Springer-Verlag},
   place={Berlin},
   date={2004},
   pages={xii+436},
   isbn={3-540-00832-2},
   review={\MR{2030225 (2004m:14070)}},
}

\bib{Borel-extension}{article}{
   author={Borel, Armand},
   title={Some metric properties of arithmetic quotients of symmetric spaces
   and an extension theorem},
   note={Collection of articles dedicated to S. S. Chern and D. C. Spencer
   on their sixtieth birthdays},
   journal={J. Differential Geometry},
   volume={6},
   date={1972},
   pages={543--560},
   issn={0022-040X},
   review={\MR{0338456 (49 \#3220)}},
}

\bib{Campedelli}{article}{
  author={Campedelli, Luigi},
  title={Sopra alcuni piani doppi notevoli con curva di diramazione
    del decimo ordine},
  language={Italian},
  journal={Atti Accad. Naz. Lincei, Rend., VI. Ser.},
  volume={15},
  pages={536--542},
  date={1932},
}

\bib{Pardini-involution}{article}{
   author={Calabri, A.},
   author={Mendes Lopes, M.},
   author={Pardini, R.},
   title={Involutions on numerical Campedelli surfaces},
   journal={Tohoku Math. J. (2)},
   volume={60},
   date={2008},
   number={1},
   pages={1--22},
   issn={0040-8735},
   review={\MR{2419034 (2009d:14045)}},
   eprint={\arXiv{math/0511391}},
}

\bib{Cossec-models}{article}{
   author={Cossec, François R.},
   title={Projective models of Enriques surfaces},
   journal={Math. Ann.},
   volume={265},
   date={1983},
   number={3},
   pages={283--334},
   issn={0025-5831},
   review={\MR{721398 (86d:14035)}},
   doi={10.1007/BF01456021},
}

\bib{Cossec-picard}{article}{
   author={Cossec, François R.},
   title={On the Picard group of Enriques surfaces},
   journal={Math. Ann.},
   volume={271},
   date={1985},
   number={4},
   pages={577--600},
   issn={0025-5831},
   review={\MR{790116 (86k:14027)}},
   doi={10.1007/BF01456135},
}

\bib{CossecDolgachev}{book}{
   author={Cossec, François R.},
   author={Dolgachev, Igor V.},
   title={Enriques surfaces. I},
   series={Progress in Mathematics},
   volume={76},
   publisher={Birkhäuser Boston Inc.},
   place={Boston, MA},
   date={1989},
   pages={x+397},
   isbn={0-8176-3417-7},
   review={\MR{986969 (90h:14052)}},
}

\bib{Dolgachev-aut}{article}{
   author={Dolgachev, I.},
   title={On automorphisms of Enriques surfaces},
   journal={Invent. Math.},
   volume={76},
   date={1984},
   number={1},
   pages={163--177},
   issn={0020-9910},
   review={\MR{739632 (85j:14076)}},
   doi={10.1007/BF01388499},
}

\bib{Dolgachev-refl}{article}{
   author={Dolgachev, Igor V.},
   title={Reflection groups in algebraic geometry},
   journal={Bull. Amer. Math. Soc. (N.S.)},
   volume={45},
   date={2008},
   number={1},
   pages={1--60},
   issn={0273-0979},
   review={\MR{2358376 (2009h:14001)}},
   doi={10.1090/S0273-0979-07-01190-1},
}

\bib{DGK}{article}{
   author={Dolgachev, I.},
   author={van Geemen, B.},
   author={Kond\=o, S.},
   title={A complex ball uniformization of the moduli space of cubic
   surfaces via periods of $K3$ surfaces},
   journal={J. Reine Angew. Math.},
   volume={588},
   date={2005},
   pages={99--148},
   issn={0075-4102},
   review={\MR{2196731 (2006h:14051)}},
   doi={10.1515/crll.2005.2005.588.99},
}

\bib{Friedman}{article}{
   author={Friedman, Robert},
   title={A new proof of the global Torelli theorem for $K3$ surfaces},
   journal={Ann. of Math. (2)},
   volume={120},
   date={1984},
   number={2},
   pages={237--269},
   issn={0003-486X},
   review={\MR{763907 (86k:14028)}},
   doi={10.2307/2006942},
}

\bib{Gerstein}{article}{
   author={Gerstein, Larry J.},
   title={Integral decomposition of hermitian forms},
   journal={Amer. J. Math.},
   volume={92},
   date={1970},
   pages={398--418},
   issn={0002-9327},
   review={\MR{0269592 (42 \#4487)}},
}

\bib{Green}{article}{
   author={Green, Mark L.},
   title={The period map for hypersurface sections of high degree of
     an
   arbitrary variety},
   journal={Compositio Math.},
   volume={55},
   date={1985},
   number={2},
   pages={135--156},
   issn={0010-437X},
   review={\MR{795711 (87b:32038)}},
}

\bib{Godeaux}{article}{
  author={Godeaux, L.},
  title={Sur une surface algébrique de genre zéro et de bigenre deux},
  journal={Atti Acad. Naz. Lincei},
  volume={14},
  date={1931},
  pages={479--481},
}

\bib{Griffiths}{article}{
   author={Griffiths, Phillip A.},
   title={On the periods of certain rational integrals. I, II},
   journal={Ann. of Math. (2) 90 (1969), 460-495; ibid. (2)},
   volume={90},
   date={1969},
   pages={496--541},
   issn={0003-486X},
   review={\MR{0260733 (41 \#5357)}},
}

\bib{Horikawa1}{article}{
   author={Horikawa, Eiji},
   title={On the periods of Enriques surfaces. I},
   journal={Math. Ann.},
   volume={234},
   date={1978},
   number={1},
   pages={73--88},
   issn={0025-5831},
   review={\MR{0491725 (58 \#10927a)}},
}

\bib{Horikawa2}{article}{
   author={Horikawa, Eiji},
   title={On the periods of Enriques surfaces. II},
   journal={Math. Ann.},
   volume={235},
   date={1978},
   number={3},
   pages={217--246},
   issn={0025-5831},
   review={\MR{0491726 (58 \#10927b)}},
}

\bib{Konno}{article}{
   author={Konno, Kazuhiro},
   title={On the variational Torelli problem for complete intersections},
   journal={Compositio Math.},
   volume={78},
   date={1991},
   number={3},
   pages={271--296},
   issn={0010-437X},
   review={\MR{1106298 (92f:14009)}},
}

\bib{Looijenga-compact2}{article}{
   author={Looijenga, Eduard},
   title={Compactifications defined by arrangements. II. Locally symmetric
   varieties of type IV},
   journal={Duke Math. J.},
   volume={119},
   date={2003},
   number={3},
   pages={527--588},
   issn={0012-7094},
   review={\MR{2003125 (2004i:14042b)}},
   doi={10.1215/S0012-7094-03-11933-X},
}

\bib{LooijengaSwierstra}{article}{
   author={Looijenga, Eduard},
   author={Swierstra, Rogier},
   title={On period maps that are open embeddings},
   journal={J. Reine Angew. Math.},
   volume={617},
   date={2008},
   pages={169--192},
   issn={0075-4102},
   review={\MR{2400994 (2010a:32030)}},
   doi={10.1515/CRELLE.2008.029},
}

\bib{Looijenga-4fold}{article}{
   author={Looijenga, Eduard},
   title={The period map for cubic fourfolds},
   journal={Invent. Math.},
   volume={177},
   date={2009},
   number={1},
   pages={213--233},
   issn={0020-9910},
   review={\MR{2507640}},
   doi={10.1007/s00222-009-0178-6},
}

\bib{Lusztig}{article}{
   author={Lusztig, Gheorghe},
   title={Novikov's higher signature and families of elliptic
     operators},
   journal={J. Differential Geometry},
   volume={7},
   date={1972},
   pages={229--256},
   issn={0022-040X},
   review={\MR{0322889 (48 \#1250)}},
}

\bib{MSY}{article}{
   author={Matsumoto, Keiji},
   author={Sasaki, Takeshi},
   author={Yoshida, Masaaki},
   title={The monodromy of the period map of a $4$-parameter family of $K3$
   surfaces and the hypergeometric function of type $(3,6)$},
   journal={Internat. J. Math.},
   volume={3},
   date={1992},
   number={1},
   pages={164},
   issn={0129-167X},
   review={\MR{1136204 (93a:33029)}},
   doi={10.1142/S0129167X92000023},
}

\bib{MLP}{article}{
   author={Mendes Lopes, Margarida},
   author={Pardini, Rita},
   title={A new family of surfaces with $p_g=0$ and $K^2=3$},
   language={English, with English and French summaries},
   journal={Ann. Sci. École Norm. Sup. (4)},
   volume={37},
   date={2004},
   number={4},
   pages={507--531},
   issn={0012-9593},
   review={\MR{2097891 (2005h:14095)}},
   doi={10.1016/j.ansens.2004.04.001},
}

\bib{MirandaMorrison}{article}{
   author={Miranda, Rick},
   author={Morrison, David R.},
   title={The number of embeddings of integral quadratic forms. I},
   journal={Proc. Japan Acad. Ser. A Math. Sci.},
   volume={61},
   date={1985},
   number={10},
   pages={317--320},
   issn={0386-2194},
   review={\MR{834537 (87j:11031a)}},
}

\bib{Miyaoka}{article}{
   author={Miyaoka, Yoichi},
   title={Tricanonical maps of numerical Godeaux surfaces},
   journal={Invent. Math.},
   volume={34},
   date={1976},
   number={2},
   pages={99--111},
   issn={0020-9910},
   review={\MR{0409481 (53 \#13236)}},
}

\bib{Miyaoka-campedelli}{article}{
   author={Miyaoka, Y.},
   title={On numerical Campedelli surfaces},
   conference={
      title={Complex analysis and algebraic geometry},
   },
   book={
      publisher={Iwanami Shoten},
      place={Tokyo},
   },
   date={1977},
   pages={113--118},
   review={\MR{0447258 (56 \#5573)}},
}

\bib{Morrison}{article}{
   author={Morrison, David R.},
   title={On the moduli of Todorov surfaces},
   conference={
      title={Algebraic geometry and commutative algebra, Vol.\ I},
   },
   book={
      publisher={Kinokuniya},
      place={Tokyo},
   },
   date={1988},
   pages={313--355},
   review={\MR{977767 (90a:14051)}},
}

\bib{Naie}{article}{
   author={Naie, Daniel},
   title={Surfaces d'Enriques et une construction de surfaces de type
   général avec $p\sb g=0$},
   language={French},
   journal={Math. Z.},
   volume={215},
   date={1994},
   number={2},
   pages={269--280},
   issn={0025-5874},
   review={\MR{1259462 (94m:14055)}},
}

\bib{Namikawa}{article}{
   author={Namikawa, Yukihiko},
   title={Periods of Enriques surfaces},
   journal={Math. Ann.},
   volume={270},
   date={1985},
   number={2},
   pages={201--222},
   issn={0025-5831},
   review={\MR{771979 (86j:14035)}},
}

\bib{Nikulin}{article}{
   author={Nikulin, V. V.},
   title={Integer symmetric bilinear forms and some of their geometric
   applications},
   language={Russian},
   journal={Izv. Akad. Nauk SSSR Ser. Mat.},
   volume={43},
   date={1979},
   number={1},
   pages={111--177, 238},
   issn={0373-2436},
   review={\MR{525944 (80j:10031)}},
}

\bib{Pardini-covers}{article}{
   author={Pardini, Rita},
   title={Abelian covers of algebraic varieties},
   journal={J. Reine Angew. Math.},
   volume={417},
   date={1991},
   pages={191--213},
   issn={0075-4102},
   review={\MR{1103912 (92g:14012)}},
   doi={10.1515/crll.1991.417.191},
}

\bib{Reid-k22}{article}{
  author={Reid, Miles},
  title={Surfaces with $p_g=0$, $K^2=2$},
  journal={unpublished},
  eprint={\href{http://www.warwick.ac.uk/~masda/surf/K2=2.pdf}
    {\texttt{http\string:// www.warwick.ac.uk/$\sim$masda/surf}}},
}

\bib{Serre}{book}{
   author={Serre, Jean-Pierre},
   title={Cours d'arithmétique},
   language={French},
   note={Deuxième édition revue et corrigée;
   Le Mathématicien, No. 2},
   publisher={Presses Universitaires de France},
   place={Paris},
   date={1977},
   pages={188},
   review={\MR{0498338 (58 \#16473)}},
}

\bib{Todorov}{article}{
   author={Todorov, Andrei N.},
   title={A construction of surfaces with $p_{g}=1$, $q=0$ and $2\leq
   (K^{2})\leq 8$. Counterexamples of the global Torelli theorem},
   journal={Invent. Math.},
   volume={63},
   date={1981},
   number={2},
   pages={287--304},
   issn={0020-9910},
   review={\MR{610540 (82k:14034)}},
   doi={10.1007/BF01393879},
}

\bib{Voisin-4fold}{article}{
   author={Voisin, Claire},
   title={Théorème de Torelli pour les cubiques de $\mathbf P^5$},
   language={French},
   journal={Invent. Math.},
   volume={86},
   date={1986},
   number={3},
   pages={577--601},
   issn={0020-9910},
   review={\MR{860684 (88g:14006)}},
   doi={10.1007/BF01389270},
}

\bib{Voisin}{book}{
   author={Voisin, Claire},
   title={Théorie de Hodge et géométrie algébrique complexe},
   language={French},
   series={Cours Spécialisés},
   volume={10},
   publisher={Société Mathématique de France},
   place={Paris},
   date={2002},
   pages={viii+595},
   isbn={2-85629-129-5},
   review={\MR{1988456 (2005c:32024a)}},
}

\bib{Wall-ortho}{article}{
   author={Wall, C. T. C.},
   title={On the orthogonal groups of unimodular quadratic forms},
   journal={Math. Ann.},
   volume={147},
   date={1962},
   pages={328--338},
   issn={0025-5831},
   review={\MR{0138565 (25 \#2009)}},
}

\end{biblist}
\end{bibsection}
\end{document}